\documentclass[10pt]{article}
\usepackage[utf8]{inputenc}
\usepackage[english]{babel}

\usepackage{eso-pic}
\usepackage{amsmath} 
\usepackage{amssymb}
\usepackage{mathrsfs}
\usepackage{amsthm}
\usepackage{tikz}
\usetikzlibrary{calc} 
\usetikzlibrary{math}
\usepackage{tikz-cd} 
\usepackage{graphicx}
\usetikzlibrary{shapes.misc}
\usepackage{appendix}
\usepackage{xcolor}
\usepackage{graphicx}
\usepackage{lineno}
\usepackage{titlesec}
\usepackage{scalerel}[2016/12/29]
\usepackage[nottoc, notlof, notlot]{tocbibind}
\usepackage{hyperref}
\usepackage{subfiles}
\usepackage[all,cmtip]{xy}

\usepackage{setspace}

\usepackage{mathtools}
\usepackage{moreverb}
\usepackage{geometry}

\newcommand{\N}{\mathbb{N}}

\newcommand{\U}{\mathcal{U}}
\newcommand{\A}{\mathcal{A}}

\newcommand{\F}{\mathbb{F}_2}
\newcommand{\Fp}{\mathbb{F}_p}
\newcommand{\Fc}{\mathcal{F}}
\newcommand{\Fo}{\mathcal{F}_\omega}
\newcommand{\K}{\mathcal{K}}

\newcommand{\E}{\mathcal{V}^f}
\newcommand{\Ef}{\mathcal{V}}

\newcommand{\Hom}{\text{Hom}}

\newcommand{\HomK}{\text{Hom}_\K}
\newcommand{\HomF}{\text{Hom}_{\F}}
\newcommand{\HomU}{\text{Hom}_\U}
\newcommand{\HomKf}{\text{Hom}_{\K\text{f.g.}}}

\newcommand{\HomC}{\text{Hom}_{\mathcal{VI}}}

\newcommand{\Nil}{\mathcal{N}il}

\newcommand{\reg}{\text{reg}}
\newcommand{\id}{\text{id}}

\newcommand{\Ima}{\text{Im}}
\newcommand{\Gal}{\text{Gal}}
\newcommand{\coker}{\text{coker}}

\newcommand{\SVF}{\mathcal{S}\text{et}^{(\mathcal{V}^f)^\text{op}}}
\newcommand{\SVI}{\mathcal{S}\text{et}^{(\mathcal{VI})^\text{op}}}

\newcommand{\FVF}{\mathcal{F}\text{in}^{(\mathcal{V}^f)^\text{op}}}
\newcommand{\FVI}{\mathcal{F}\text{in}^{(\mathcal{VI})^\text{op}}}

\newcommand{\PFVF}{\mathcal{P}\text{fin}^{(\mathcal{V}^f)^\text{op}}}

\def\commutatif{\ar@{}[rd]|{\circlearrowleft}}

\newtheorem{theorem}{Theorem}[section]

\newtheorem{proposition}[theorem]{Proposition}
\newtheorem{lemma}[theorem]{Lemma}
\newtheorem{corollary}[theorem]{Corollary}
\newtheorem{innercustomthm}{Theorem}
\newenvironment{customthm}[1]
{\renewcommand\theinnercustomthm{#1}\innercustomthm}
{\endinnercustomthm}

\theoremstyle{definition}

\newtheorem{definition}[theorem]{Definition}
\newtheorem{remark}[theorem]{Remark}

\theoremstyle{definition}
\newtheorem{example}[theorem]{Example}

\theoremstyle{definition}
\newtheorem{notation}[theorem]{Notation}

\title{PRESHEAVES ON $\mathcal{VI}$, $nil$-CLOSED UNSTABLE ALGEBRAS AND THEIR CENTRES}
\author{OURIEL BLOEDE}

\date{}
\begin{document}

\maketitle

\begin{abstract}
A $nil$-closed, noetherian, unstable algebra $K$ over the Steenrod Algebra is determined, up to isomorphism, by the functor $\HomKf(K,H^*(\_))$, which is a presheaf on the category $\mathcal{VI}$ of finite dimensional vector spaces and injections, by the theory of Henn-Lannes-Schwartz. In this article, we use this theory to study the centre, in the sense of Heard, of a $nil$-closed noetherian unstable algebra.

For $F$ a presheaf on $\mathcal{VI}$,  we construct a groupoid $\mathcal{G}_F$ which encodes $F$. Then, taking $F:=\HomKf(K,H^*(\_))$, we show how the centre of $K$ is determined by the associated groupoid. We also give a generalisation of the second theorem of Adams-Wilkerson, defining sub-algebras $H^*(W)^\mathcal{G}$ of $H^*(W)$ for appropriate groupoids $\mathcal{G}$.

There is a $H^*(C)$-comodule structure on $K$ that is associated with the centre. For $K$ integral, we explain how the algebra of primitive elements of this $H^*(C)$-comodule structure is also determined by the groupoid associated with $\HomKf(K,H^*(\_))$. Along the way, we prove that this algebra of primitive elements is also noetherian.
\end{abstract}

\section{Introduction}

\subsection{The two theorems of Adams-Wilkerson}

For $p$ a prime number, $W$ a finite dimensional $\Fp$-vector space, and for $BW$ the classifying space of $W$, $W\mapsto H^*(W):=H^*(BW;\Fp)$ defines a functor from $\E$, the category of finite dimensional vector spaces, to $\K$, the category of unstable algebra over the Steenrod algebra over $\Fp$. Then, for $G$ a sub-group of $\text{Gl}(W)$, $G$ acts on $H^*(W)$ and we can define $H^*(W)^G\in\K$, the unstable algebra of invariant elements of $H^*(W)$ under the action of $G$. For $K$ an integral, noetherian, unstable algebra of transcendence degree $\dim(W)$, the first theorem of Adams-Wilkerson states that there always exists an injection $\phi$ from $K$ to $H^*(W)$ and that this injection induces a structure of finitely generated $K$-module on $H^*(W)$. Then, for $\Gal(\phi)$ the sub group of $\text{Gl}(W)$ whose elements are the automorphisms $\alpha$ such that $\alpha^*\phi=\phi$ with $\alpha^*\phi$ the composition of $\phi$ with the isomorphism induced by $\alpha$ on $H^*(W)$, the image of $\phi$ is a sub-algebra of $H^*(W)^{\Gal(\phi)}$. $H^*(W)^{\Gal(\phi)}$ is a first approximation of $K$, but usually $\phi$ does not define an isomorphism between $K$ and $H^*(W)^{\Gal(\phi)}$.\\

The category of unstable modules over the Steenrod algebra admits a localizing subcategory $\Nil$ (cf \cite{S1}), whose objects are called nilpotent unstable modules. We call an unstable module $nil$-closed if its localization away from $\Nil$ is an isomorphism. The algebras $H^*(W)^G$ are $nil$-closed, hence $\phi$ cannot be an isomorphism when $K$ is not $nil$-closed. In this case, and when $p=2$ the second theorem of Adams-Wilkerson gives a condition for $\phi$ to be an isomorphism, when $K$ is $nil$-closed. Namely, if $K$ is \begin{enumerate}
    \item noetherian,
    \item integral,
    \item $nil$-closed,
    \item integrally closed in its field of fractions,
\end{enumerate} the morphism $\phi$ of the first theorem of Adams-Wilkerson is an isomorphism between $K$ and $H^*(W)^{\Gal(\phi)}$. (There is a similar statement when $p$ is odd.)\\

The first objective of this article is the following, for $\mathcal{G}$ a groupoid whose objects are the sub-spaces of $W$ and whose morphisms satisfy a restriction property, we define an object $H^*(W)^\mathcal{G}$ that generalises the algebra of invariants $H^*(W)^G$. We will show that the $H^*(W)^\mathcal{G}$ give a complete list of the $nil$-closed, noetherian, unstable sub-algebras of $H^*(W)$.

\subsection{The centre of an unstable algebra}

In \cite{DW1}, Dwyer and Wilkerson introduced the notion of a central element of an unstable algebra, this notion allowed them to exhibit the only exotic finite loop space at prime 2 in \cite{DW12}. In the case where $K$ is noetherian and connected, the set of central elements of $K$ coincides with the set of pairs $(V,\phi)$ such that
\begin{enumerate}
    \item $\phi\in\HomK(K,H^*(V))$,
    \item $K$ admits a structure $\kappa$ of $H^*(V)$-comodule in $\K$, such that the following diagram commutes:
    $$\xymatrix{ 
K\ar[rr]^\kappa\ar[rrd]_\phi & & K\otimes H^*(V)\ar[d]^-{\epsilon_K\otimes\id} & \\
& & H^*(V),}$$ where $\epsilon_K$ denotes the augmentation of $K$ (which is uniquely defined because of the connectedness of $K).$
\end{enumerate}

In \cite{heard2020topological}, Heard showed that for $K$ noetherian, $K$ admits a unique (up to isomorphism) central element $(C,\gamma)$ such that $\gamma$ induces a structure of finitely generated $K$-module on $H^*(C)$ and $\dim(C)$ is maximal among such central elements. Heard called this central element the centre of $K$. The centre of an unstable algebra have been shown to be an important invariant. In \cite{K11} and \cite{K12}, Kuhn used it to approximate the depth of $K$ as well as invariants $d_0(K)$ and $d_1(K)$ introduced by Henn, Lannes and Schwartz in \cite{HLS1}, in the case where $K$ is the cohomology of a group. Heard generalised those results for $K$ noetherian in \cite{heard2019depth} and \cite{heard2020topological}.\\

For $K$ noetherian, since the centre of $K$ is associated with a $H^*(C)$-comodule structure on $K$, it gives rise to a second invariant: the sub-algebra of primitive elements of $K$ under this $H^*(C)$-comodule structure. The second objective of this article is to explain how the centre of $H^*(W)^\mathcal{G}$ and its sub-algebra of primitive elements are determined by $\mathcal{G}$. We will then be able to classify $nil$-closed, noetherian sub-algebras of $H^*(W)$ with a given centre and algebra of primitive elements.

\subsection{The category $\SVI$}

We consider $\K/\Nil$ the localization of $\K$ with respect to morphisms whose kernels and cokernels are nilpotent unstable modules. In \cite{HLS2}, Henn, Lannes and Schwartz proved that the functor which sends an unstable algebra $K$ to $\HomK(K,H^*(\_))$, the functor which maps $V\in\E$ to $\HomK(K,H^*(V))$, induces an equivalence of category between $\K/\Nil$ and a certain sub-category of $\SVF$ the category of contravariant functors from $\E$ to $\mathcal{S}et$. Furthermore, they used abundantly the fact that, when $K$ is noetherian, the functor $\HomK(K,H^*(\_))$ is fully determined by the $\HomKf(K,H^*(V))$ for $V$ running through $\E$, where $\HomKf(K,H^*(V))$ is the subset of $\HomK(K,H^*(V))$ whose objects are the morphisms from $K$ to $H^*(V)$ which turn $H^*(V)$ into a finitely generated $K$-module.\\

When $K$ is noetherian, $\HomKf(K,H^*(\_))$ defines a contravariant functor from $\mathcal{VI}$ to $\mathcal{S}et$, where $\mathcal{VI}$ is the wide sub-category of $\E$ whose morphisms are the injective morphisms. The category $\SVI$ is a lot easier to understand than the category $\SVF$, and that is going to be our main tool.\\

The two first sections of this article consist of recollections about the equivalences of categories constructed in \cite{HLS2} and the definition of central elements of an unstable algebra. In the third section, we will introduce the category $\SVI$ and its connections with noetherian unstable algebras. Furthermore, we will define a notion of central elements for functors in $\SVI$. For $F=\HomKf(K,H^*(\_))$, with $K$ a noetherian $nil$-closed unstable algebra, the central elements of $F$ will coincide with the central elements $(V,\phi)$ of $K$, such that $\phi$ turns $H^*(V)$ into a finitely generated $K$-module.  
\begin{customthm}{\ref{principal0}}
For $K$ a noetherian unstable algebra, $F=\HomKf(K,H^*(\_))$ and $\phi\in F(V)$, $(V,\phi)$ is central for $K$ if and only if it is central for $F$.
\end{customthm}

\subsection{The groupoid $\mathcal{G}_F$}

An injection $\phi$ from a noetherian unstable algebra $K$ to some $\prod\limits_{i\in I}^{\wedge}H^*(W_i)$, where $\prod\limits^{\wedge}$ denotes the product in the category of connected unstable algebras and such that $\phi$ turns each $H^*(W_i)$ into a finitely generated $K$-module, induces a surjection from $\bigsqcup\limits_{i\in I}\HomC(\_,W_i)$ to $\HomKf(K,H^*(\_))$. For $I$ a set, we consider $_{(W_i)_{i\in I}}\SVI$, the category whose objects are pairs $(F,q_F)$ with $F\in\SVI$ and $q_F$ a natural surjection from $\bigsqcup\limits_{i\in I}\HomC(\_,W_i)$ to $F$.\\ 

In the fourth section, we define an application which sends an object $(F,q_F)\in\ _{(W_i)_{i\in I}}\SVI$ to a groupoid $\mathcal{G}_{(F,q_F)}$ whose set of objects is the disjoint union of the sub-spaces of the $W_i$. This groupoid satisfies a property called the restriction property. The first main result of this article is that the isomophism classes of objects in $_{(W_i)_{i\in I}}\SVI$ are in one-to-one correspondence with such groupoids. This is stated in the following theorem, where $\sim_\mathcal{G}$ is an equivalence relation on $\bigsqcup\limits_{i\in I}\HomC(\_,W_i)$ characterised by the groupoid $\mathcal{G}$.

\begin{customthm}{\ref{principal1}}
\begin{enumerate}
    \item For $\mathcal{G}$ a groupoid whose objects are the sub-spaces of the $W_i$ and whose morphisms are isomorphisms of vector spaces such that $\mathcal{G}$ has the restriction property, $$\mathcal{G}_{(\bigsqcup\limits_{i\in I}\HomC(\_,W_i)/\sim_\mathcal{G},q)}=\mathcal{G},$$ for $q$ the canonical surjection from $\bigsqcup\limits_{i\in I}\HomC(\_,W_i)$ to $\bigsqcup\limits_{i\in I}\HomC(\_,W_i)/\sim_\mathcal{G}$.
    \item Conversely, let $F\in\ _{(W_i)_{i\in I}}\SVI$. Then, $F$ is isomorphic to $\bigsqcup\limits_{i\in I}\HomC(\_,W_i)/\sim_{\mathcal{G}_F}$. 
\end{enumerate}
\end{customthm}

We also show, in Theorem \ref{principal2}, how the central elements of $\bigsqcup\limits_{i\in I}\HomC(\_,W_i)/\sim_\mathcal{G}$ are determined by the groupoid $\mathcal{G}$.\\

Now, for $K$ a noetherian, $nil$-closed, unstable sub algebra of $H^*(W)$ with transcendence degree $\dim(W)$, the inclusion of $K$ in $H^*(W)$ turns $\HomKf(K,H^*(\_))$ into an object of $_W\SVI$. This implies the following: 
\begin{customthm}{\ref{principal4}}
For all $W\in\E$, there is a one-to-one correspondence between the set of $nil$-closed and noetherian sub-algebras of $H^*(W)$ whose transcendence degree is $\dim(W)$ and the set of groupoids with the restriction property, whose objects are the sub-vector spaces of $W$.
\end{customthm}

Then, for $\mathcal{G}$ a groupoid with the restriction property whose objects are sub-spaces of $W$, we will denote by $H^*(W)^\mathcal{G}$ the noetherian, $nil$-closed sub-algebra of $H^*(W)$ such that the groupoid associated with the surjection $\HomC(\_,W)\twoheadrightarrow \HomKf(H^*(W)^\mathcal{G},H^*(\_))$ is $\mathcal{G}$.\\

By Theorems \ref{principal0} and \ref{principal2}, the central elements of $H^*(W)^\mathcal{G}$ are determined by $\mathcal{G}$. In the fifth section we will prove the following, where $P(K,x)$ denotes the algebra of primitive elements of $K$ for the $H^*(V)$-comodule structure induced by a central element $(V,x)$.
\begin{customthm}{\ref{principal7}}
Let $K$ be a noetherian unstable sub algebra of $H^*(W)$ of finite transcendence degree $\dim(W)$ such that $(V,\delta^*\phi)$ is central, for $\phi$ the inclusion of $K$ in $H^*(W)$ and $\delta$ some morphism from $H^*(W)$ to $H^*(V)$. Then, $P(K,\delta^*\phi)$ is $nil$-closed and noetherian.
\end{customthm}
 In this context, we will identify naturally $P(K,\delta^*\phi)$ with a sub algebra of $H^*(W/\Ima(\delta))$ of transcendence degree $\dim(W/\Ima(\delta))$. By Theorem \ref{principal4}, $P(K,\delta^*\phi)$ identifies with some $H^*(W/\Ima(\delta))^{\mathcal{G}'}$. In the case where $K=H^*(W)^\mathcal{G}$, we explain in Theorem \ref{principal6} how to determine $\mathcal{G}'$ from $\mathcal{G}$.\\

We conclude this article by giving examples on how to use those constructions to determine $nil$-closed, noetherian, integral, unstable algebras whose transcendence degree is fixed, with a $H^*(V)$-comodule structure whose primitive elements are isomorphic to some $H^*(U)^{\mathcal{G}'}$ with $U$ in $\E$ and $\mathcal{G}'$ a groupoid whith the restriction property and whose objects are the sub spaces of $U$. For example:

\begin{customthm}{\ref{principal8}}
Let $K$ be a noetherian, $nil$-closed, integral, unstable algebra of transcendence degree $d$. We assume that the centre of $K$ is of dimension $d-1$. Then, there exists $G$, a sub-group of $\text{Gl}(W)$, such that $K$ is isomorphic to the algebra of invariant elements $H^*(W)^G$ with $\dim(W)=d$. Furthermore, $G$ satisfies that the set of element $x\in W$ such that $g(x)=x$ for all $g\in G$ is a sub-vector space of $W$ of dimension $d-1$.

\end{customthm}

This Theorem implies that for $K$ noetherian, $nil$-closed, integral of transcendence degree $d$ whose centre is of dimension $d-1$ and for $\phi$ the injection of Adams-Wilkerson, $\phi$ is an isomorphism from $K$ to $H^*(W)^{\Gal(\phi)}$.\\

\textbf{Acknowledgements: }I am thankful to Geoffrey Powell for its careful proofreading. I also want to thank Antoine Boivin for his help in computing the $H^*(W)^\mathcal{G}$ of the last section. This work was partially supported by the ANR Project {\em ChroK}, {\tt ANR-16-CE40-0003}.

\section{Unstable algebras over the Steenrod algebra and the functor $f$}

In this section, we recall some known facts about Lannes' $T$ functor as well as results from \cite{HLS2} about the localization of the categories of unstable modules and algebras modulo away from nilpotent objects. Recollections about unstable algebras, unstable modules and nilpotent objects can be found in \cite{S1}. In the following, $\A$ denotes the Steenrod algebra over $\Fp$ with $p$ a prime number, $\U$ and $\K$ denote the category of unstable modules and unstable algebras over $\A$ and $\Nil$ denotes the class of nilpotent objects in $\U$.

\subsection{The $T$ functor}

Let us recall the definition of Lannes' $T$ functor.

\begin{theorem}\cite[Proposition 2.1]{L1}
For $V$, a finite dimensional vector space the functor $- \otimes H^*(V)$ has a left adjoint $T_V$.
\end{theorem}

\begin{proposition}\cite[Proposition 3.8.4]{S1}\label{Tvgebre}
\begin{enumerate}
    \item Let $K$ be an unstable algebra, then $T_V(K)$ is given in a natural way a structure of unstable algebra.
    \item Then, $T_V$ defines a functor from $\K$ to $\K$ which is left adjoint to the tensor product with $H^*(V)$ in $\K$.
\end{enumerate}

\end{proposition}

\begin{example}\cite[3.9.1]{S1}\label{hetv}
For $V$ and $W$ two finite dimensional $\Fp$-vector spaces, there is an isomorphism of unstable algebras, natural in both $V$ and $W$, $T_V(H^*(W))\cong H^*(W)\otimes\Fp^{\Hom(V,W)}$. By the adjonction property, we get that $(\Fp^{\Hom(V,W)})^\sharp\cong\HomU(H^*(W),H^*(V)).$ In other words, $\Fp\left[\Hom(V,W)\right]\cong\HomU(H^*(W),H^*(V))$ which is a theorem first proved by Adams, Gunawardena and Miller.
\end{example}

\subsection{$\Nil$-localisation of unstable modules}

 The class of nilpotent modules is a Serre class in $\U$, we recall the existence of an equivalence of categories between $\U/\Nil$ (defined as in \cite{G1}) and a category of functors. The proofs can be found in \cite{HLS2}.

\begin{theorem}\cite[Part I.4]{HLS2}
There is an adjunction of functors:
$$\xymatrix{r_1 : \U\ar@<2pt>[r] & \ar@<2pt>[l] \U/\Nil: s_1,}$$
such that, for $\phi$ a morphism of unstable modules, $r_1(\phi)$ is an isomorphism if and only if $\ker(\phi)$ and $\coker(\phi)$ are objects in $\Nil$. 

Then, $\U/\Nil$ satisfies the following universal property: for $A$ an abelian category and $F\ :\ \U\rightarrow A$ an exact functor such that for all $M\in\Nil$, $F(M)=0$, there exists a unique $G\ :\ \U/\Nil\rightarrow A$ such that $F=G\circ r_1$.
\end{theorem}

\begin{definition}
For $M$ an unstable module, $l_1(M):=s_1\circ r_1(M)$ is the $nil$-localisation of $M$, and we say that $M$ is $nil$-closed if the unit of the adjunction $M\rightarrow l_1(M)$ is an isomorphism.
\end{definition}

For $\Fc$ the category of functors from the category $\E$ of finite dimensional $\Fp$-vector spaces to the category $\Ef$ of $\Fp$-vector spaces, we consider $f\ :\ \U\rightarrow \Fc$ the functor which assigns to $M$ in $\U$, $$f(M)\ :\ V\mapsto T_V(M)^0.$$\\

The class $\Nil$ satisfies that $f(M)=0$ if and only if $M\in\Nil$. Then, $f$ induces a functor $f'$ from $\U/\Nil$ to $\Fc$ such that $f=f'\circ r_1$.\\

We will denote by $\Fo$ the essential image of $f$ in $\Fc$.

\begin{theorem}\cite[Theorem 5.2.6]{S1}\label{eqnil1}
The functor $f'$ induces an equivalence of categories between $\U/\Nil$ and the category $\Fo$.
\end{theorem}

In \cite{HLS2}, the authors exhibit a right adjoint to $f'$, $m'$ such that the restriction of $m'$ to $\Fo$ is an inverse of the equivalence of categories induced by $f'$.

\begin{definition}\label{mmm}

Let $m\ :\ \Fc\rightarrow \U$ be the composition of $m'$ with $s_1$.
\end{definition}

\begin{definition}
For $F\in\Fc$ and $V\in\E$, let $\Delta_VF$ denote the object of $\Fc$ such that $\Delta_VF(W)=F(V\oplus W)$ and $\Delta_VF(\alpha)=F(\id_V\oplus\alpha)$ for all $W\in\E$ and for all morphism $\alpha$ in $\E$.
\end{definition}

We recall the following results from \cite{K1}

\begin{proposition}\label{deltav}
For $M\in\U$ and $V\in\E$, $f(T_V(M))=\Delta_V(f(M))$.\\

For $F$ an object of $\mathcal{F}$, $T_V(m(F))\cong m(\Delta_V(F))$.
\end{proposition}

\begin{corollary}\label{tutut}
For $M\in \U$ $nil$-closed, $T_V(M)$ is also $nil$-closed.
\end{corollary}

\subsection{$\Nil$-localisation of unstable algebras}

Since $\K$ is not abelian, one cannot define a localized category of $K$ in the sense of \cite{G1}. In \cite{HLS2}, Henn, Lannes and Schwartz constructed a localized category $\K/\Nil$ with respect to the morphisms whose kernels and cokernels are in $\Nil$, in the sense of \cite{KS}. Then, the functor $f$ restricted to $\K$ factorises through a functor from $\K/\Nil$ to $\Fc$. The authors of \cite{HLS2} identified the essential image of $f$ restricted to $K$ and they deduced an equivalence of category between $\K/\Nil$ and a category of contravariant functors from $\E$ to the category of profinite sets. \\

\begin{definition}
A $p$-boolean algebra, is an algebra $B$ over $\Fp$, such that, for all $x\in B$, $x^p=x$.
\end{definition}

Since $T_V(K)$ is an unstable algebra, $T_V(K)^0$ is a $p$-boolean algebra. We can then use standard results on $p$-boolean algebras to study $f(K)$.\\

For $\mathcal{B}$ the category of $p$-boolean algebras and $B$ a $p$-boolean algebra, we consider $\Hom_{\mathcal{B}}(B,\Fp)$ the set of morphisms of $\Fp$-algebras from $B$ to $\Fp$. Then $B$ is the direct limit of its finite dimensional subalgebras $B_\alpha$. Therefore, $\Hom_{\mathcal{B}}(B,\Fp)$ is the inverse limit of the $\Hom_{\mathcal{B}}(B_\alpha,\Fp)$ which are finite. $\Hom_{\mathcal{B}}(B,\Fp)$ inherits a structure of profinite set.

 \begin{proposition}\cite{HLS2}\label{booleprofin}
For $\mathcal{P}fin$ the category of profinite sets, the functor $\text{spec}\ :\ \mathcal{B}^{\text{op}}\rightarrow\mathcal{P}fin$, where $\text{spec}(B):=\text{Hom}_{\mathcal{B}}(B,\mathbb{F}_p)$, is an equivalence of categories whose inverse is the functor which sends $S$ to the algebra of continuous maps from $S$ to $\mathbb{F}_p$, $\mathbb{F}_p^S$.
\end{proposition}

In particular, for $K\in\K$, by adjunction $\Hom_{\mathcal{B}}(T_V(K)^0,\Fp)\cong\HomK(K,H^*(V))$ and this isomorphism is an isomorphism of profinite sets, where the structure of profinite set on $\HomK(K,H^*(V))$ comes from the fact that $K$ is the direct limit of the unstable sub-algebra of $K$ which are finitely generated as $\A$-algebras. Then, $T_V(K)^0$ is isomorphic as a $p$-boolean algebra to $\Fp^{\HomK(K,H^*(V))}$. 

\begin{definition}
\begin{enumerate}
    \item Let $\PFVF$ be the category of functors from $(\E)^{\text{op}}$ to $\mathcal{P}fin$,
    \item let $\mathcal{L}$ be Lannes' linearization functor from $(\PFVF)^{\text{op}}$ to $\Fc$ defined by $\mathcal{L}(F)(V):=\Fp^{F(V)}$,
    \item let $g\ :\ \K\rightarrow(\PFVF)^{\text{op}}$ be the functor which sends $K$ to the functor $g(K)\ :\ V\mapsto\HomK(K,H^*(V))$.
\end{enumerate}
\end{definition}

We have a commutative diagram of functors:
$$\xymatrix{\K\ar[r]^-{g}\ar[d] & (\PFVF)^{\text{op}}\ar[d]^{\mathcal{L}}\\
\U \ar[r]^f &\Fc,}$$
where the functor from $\K$ to $\U$ is the forgetful functor. We denote by $\PFVF_\omega$ the full subcategory of $\PFVF$, whose objects are those whose image under $\mathcal{L}$ are in $\Fo$.

The functor $g$ has a unique factorisation of the following form:
$$\K\rightarrow \K/\Nil\rightarrow \PFVF_\omega\rightarrow \PFVF.$$

\begin{theorem}\cite[Theorem 1.5 of Part II]{HLS2}
The functor from $\K/\Nil$ to $\PFVF_\omega$ induced by $g$ is an equivalence of categories.
\end{theorem}

The following lemma will be of importance in the following.

\begin{lemma}\label{coco limite}
The functor $g$ turns injections into surjections and finite inverse limits into direct limits.
\end{lemma}

\begin{proof}
Since $f$ is exact, $f$ sends injections into injections and commutes with finite inverse limits. The result is then a consequence of the isomorphism $f(K)\cong \Fp^{g(K)}$.
\end{proof} 

\section{Connected components of $T_V(K)$ and central elements of an unstable algebra}

In this section, we recall the definition of central elements of an unstable algebra first introduced by Dwyer and Wilkerson in \cite{DW1}.

\subsection{Connected components of $T_V(K)$}

For $K$ an unstable algebra, we recall the definition of the connected components of $T_V(K)$ which is exposed in \cite{heard2019depth} and \cite{heard2020topological}. Such a decomposition exists for any graded algebra over a $p$-boolean algebra.

\begin{lemma}
For $K$ an unstable algebra which is finitely generated as an algebra over $\A$,\newline $\text{Hom}_\mathcal{K}(K,H^*(V))$ is finite. 
\end{lemma}

\begin{definition}
For $V\in\E$ and $\phi \in \text{Hom}_\mathcal{K}(K,H^*(V))$, let $T_{(V,\phi)}(K):=T_V(K)\otimes_{T_V(K)^0}\mathbb{F}_p(\phi)$, where the structure of $T_V(K)^0$-module over $\mathbb{F}_p(\phi)$ is induced by the morphism from $T_V(K)^0$ to $\mathbb{F}_p$ adjoint to $\phi$.

\end{definition}

\begin{proposition}\cite[Equation (2.6)]{heard2019depth}\label{prumi}
For $K$ an unstable algebra finitely generated as an algebra over $\A$ and $V\in\E$, we have the following natural isomorphism of unstable algebra $$T_V(K)\cong \prod\limits_{\phi\in \text{Hom}_\mathcal{K}(K,H^*(V))}T_{(V,\phi)}(K).$$

\end{proposition}

\begin{lemma}\label{upussu}
For $K$ a $nil$-closed unstable algebra, $V\in\E$ and $\phi\in\HomK(K,H^*(V))$, $T_{(V,\phi)}(K)$ is $nil$-closed.
\end{lemma}

\begin{proof}
By corollary \ref{tutut} $T_V(K)$ is $nil$-closed. By the isomorphism of proposition \ref{prumi}, $T_{(V,\phi)}(K)$ is the kernel of the morphism from $T_V(K)$ to $\bigoplus\limits_{\phi\neq\psi}T_V(K)$ which sends $x$ to the direct sum of the components of $x$ in each $T_{(V,\psi)}(K)$ with $\phi\neq\psi$. Since $ \bigoplus\limits_{\phi\neq\psi}T_V(K)$ is $nil$-closed and $l_1$ is left exact, $l_1(T_{(V,\phi)}(K))$ is the kernel of the same morphism, thus, it is isomorphic to $T_{(V,\phi)}(K)$.
\end{proof}

\subsection{Central elements of an unstable algebra}

The notion of a central element of an unstable algebra $K$ is defined by Dwyer and Wilkerson in \cite{DW1} and they used it in \cite{DW12} to exhibit the only exotic finite loop space at the prime $2$. The centre of $K$ has been studied in details in \cite{heard2019depth} and \cite{heard2020topological}.

The aim of this subsection is to recall some known facts about central elements of an unstable algebra.\\
 
 \begin{notation}
 Let $M$ be an unstable module, $V$ a finite dimensional vector space, $K$ an unstable algebra and $\phi\in\HomK(K,H^*(V))$. Denote by:
\begin{itemize}  
    \item  $\eta_{M,V}\ :\ M\rightarrow T_V(M)\otimes H^*(V)$ the unit of the adjunction between $T_V$ and $-\otimes H^*(V)$;
    \item  $\rho_{M,V}$ the following composition $$\rho_{M,V}\ :\ M\overset{\eta_{M,V}}{\longrightarrow} T_V(M)\otimes H^*(V)\overset{id\otimes \epsilon_{V}}{\longrightarrow} T_V(M), $$
where $\epsilon_V$ denote the augmentation of $H^*(V)$;
     \item  $\rho_{K,(V,\phi)}$ the composition of $\rho_{K,V}$ with the projection onto $T_{(V,\phi)}(K)$.
\end{itemize}
 \end{notation} 

\begin{remark}
The morphism $\rho_{M,V}$ identifies with $T_{\iota_0^V}(M)\ :\ M\cong T_0(M)\rightarrow T_V(M)$, the morphism induced by naturality of $T_V(M)$ with respect to $V$ by $\iota_0^V$, the injection from $0$ to $V$.
\end{remark}

\begin{definition}\label{centri}
Let $K$ be an unstable algebra and $\phi\in\HomK(K,H^*(V))$. Then, the pair $(V,\phi)$ is said to be central if $\rho_{K,(V,\phi)}\ :\ K\rightarrow T_{(V,\phi)}(K)$ is an isomorphism.

Let $\textbf{C}(K)$ be the set of central elements of $K$.

\end{definition}

The classical example and first motivation for studying the centre of an unstable algebra is the example of $H^*(G)$, the cohomology of a group $G$. The details can be found in \cite{Henn2001}.

\begin{example}
For $G$ a discrete group or a compact Lie group, $\text{Hom}_\mathcal{K}(H^*(G),H^*(V))\cong\mathbb{F}_p\left[\text{Rep}(V,G)\right]$, where $\text{Rep}(V,G)$ denote the conjugacy classes of morphisms from $V$ to $G$.\\

Let $\rho$ represent a conjugacy class in $\text{Rep}(V,G)$. We consider the morphism $V\times C_G(\rho)\rightarrow G$, where $C_G(\rho)$ denote the centraliser in $G$ of the image of $\rho$, which sends $(v,g)$ to $\rho(v)\cdot g$. It induces a morphism from $H^*(G)\rightarrow H^*(V)\otimes H^*(C_G(\rho))$. By adjunction, it gives us a morphism $T_V(H^*(G))\rightarrow H^*(C_G(\rho))$ which depends only on the conjugacy class of $\rho$. This morphism induce an isomorphism between $T_{(V,\rho)}(H^*(G))$ and $H^*(C_G(\rho))$, and $$\rho_{H^*(G),(V,\rho)}\ :\ H^*(G)\rightarrow T_{(V,\rho)}(H^*(G))\cong H^*(C_G(\rho))$$ is the morphism induced by the injection $C_G(\rho)\hookrightarrow G$.\\

Hence, $(V,\rho)$ is central if and only if the injection $C_G(\rho)\hookrightarrow G$ induces an isomorphism in cohomology.
\end{example}

\begin{definition}
Let $K$ be an unstable algebra, $K$ is connected if $K$ has an augmentation $\epsilon_K\ :\ K\rightarrow\Fp$ which induces an isomorphism $K^0\overset{\cong}{\rightarrow}\Fp$.
\end{definition}

\begin{remark}

The connected components $T_{(V,\phi)}(K)$ are connected, hence, if $K$ is not connected, $\textbf{C}(K)=\emptyset$.

\end{remark}

\begin{example}
The functor $T_0$ is the identity, hence, if $K$ is connected, $T_{(0,\epsilon_K)}(K)\cong K$, for $\epsilon_K\ :\ K\rightarrow \mathbb{F}_p$ the unique morphism of unstable algebra from $K$ to $\Fp$. Hence $(0,\epsilon_K)$ is central.
\end{example}

\begin{notation}
Let $\epsilon_{K,V}$, be the composition of $\epsilon_K$ with the injection from $\Fp$ to $H^*(V)$.
\end{notation}

For $K$ a connected unstable algebra, $I(K)$ denotes the augmentation ideal of $K$. Then the module of indecomposable elements of $K$ is defined by $\mathcal{Q}(K):=I(K)/I(K)^2$. An unstable module $M$ is said to be locally finite if, for all $x\in M$, $\A x$ is finite.

In \cite{AST_1990__191__97_0}, Dwyer and Wilkerson exhibit how, when $\mathcal{Q}(K)$ is locally finite, central elements of $K$ are related to $H^*(V)$-comodule structures on $K$. 

\begin{proposition} \cite[Proof of Theorem 3.2]{AST_1990__191__97_0}\label{indecompi}

Let $K$ be a connected unstable algebra such that $\mathcal{Q}(K)$ is locally finite as an unstable module, then  $(V,\epsilon_{K,V})$ is central for all $V\in\E$.

\end{proposition}

In particular, if $K$ is a connected, noetherian, unstable algebra, then $(V,\epsilon_{K,V})$ is central for all vector space $V$.

We recall the following results of \cite{DW1}.

\begin{proposition} \cite[Proposition 3.4]{DW1}\label{DWC}

Let $K$ be a connected unstable algebra such that $\mathcal{Q}(K)$ is locally finite. Then, for\newline $\phi\in\HomK(K,H^*(V))$, $(V,\phi)$ is central if and only if there exists a morphism from $K$ to $K\otimes H^*(V)$ such that the following diagram commutes:
$$\xymatrix{ & & K\\
K\ar[rru]^{\text{id}}\ar[rr]\ar[rrd]_\phi & & K\otimes H^*(V)\ar[u]_-{\id\otimes\epsilon_{H^*(V)}}\ar[d]^-{\epsilon_K\otimes\id} & \\
& & H^*(V).}$$

\end{proposition}

\begin{corollary}\cite{DW1}\label{DWCT}
Let $K$ be a connected unstable algebra such that $\mathcal{Q}(K)$ is locally finite. For $\phi\in\HomK(K,H^*(V))$, $(V,\phi)$ is central if and only if $K$ has a structure of $H^*(V)$-comodule $\kappa$ in $\K$, such that the following diagram commutes:
$$\xymatrix{ 
K\ar[rr]^\kappa\ar[rrd]_\phi & & K\otimes H^*(V)\ar[d]^-{\epsilon_K\otimes\id} & \\
& & H^*(V).}$$
\end{corollary}

In particular, this implies:

\begin{proposition}\label{incl}
Let $K$ be an unstable algebra such that $\mathcal{Q}(K)$ is locally finite, then for  $\phi\in\textbf{C}(K)$ and $\alpha\ :\ V\rightarrow E$ a morphism in $\E$, $(V,\alpha^*\circ \phi)\in\textbf{C}(K)$.
\end{proposition}

\begin{example}
For $W\in\E$, the addition in $W$, $\nabla_W$, induces on $H^*(W)$ a coalgebra structure in $\K$. Then, for every morphism of unstable modules $\phi$ from $H^*(W)$ to $ H^*(V)$, one can take the composition of $\nabla_W^*$ with $\id_{H^*(W)}\otimes f$ to define a $H^*(V)$-comodule structure on $H^*(W)$ satisfying the hypothesis of corollary \ref{DWCT}. Therefore $(V,\phi)$ is central.
\end{example}

\section{The category $\SVI$}

In \cite{RECTOR1984191}, Rector used the fact that for a noetherian unstable algebra $K$, the functor which maps $V$ to $\HomK(K,H^*(V))$ is fully determined by the functor which maps $V$ to $\HomKf(K,H^*(V))$, where $\HomKf(K,H^*(V))$ is the set of morphism from $K$ to $H^*(V)$ which makes $H^*(V)$ a finitely generated $K$-module. The functor $V\mapsto \HomKf(K,H^*(V))$ is defined on the category $\mathcal{VI}$, whose objects are finite dimensional vector spaces on $\Fp$ and whose morphisms are injective morphisms. In this section, we start by recalling the notion of a regular element of a functor in $\PFVF$, introduced in \cite{HLS2}. This allows us to understand the passage from $\HomK(K,H^*(\_))$ to $\HomKf(K,H^*(\_))$ as a construction on functors in $\PFVF$.

Then, we define a shift functor for contravariant functors on the category $\mathcal{VI}$. This will allow us to define a notion of centrality for objects of $\SVI$, which coincides for noetherian unstable algebras, with the notion of centrality "away from $\Nil$" for pairs $(V,\phi)$ such that $\phi\in\HomKf(K,H^*(V))$. 

\subsection{The functor $\HomKf(K,H^*(\_))$}

\begin{definition}
Let $\mathcal{VI}$ be the wide sub-category of $\E$ with only injective morphisms.
\end{definition}

In the following, we will navigate between different presheaf categories from $\E$ or $\mathcal{VI}$ to variants of $\mathcal{S}et$. Let us define all those categories.

\begin{definition}
\begin{enumerate}
    \item Let $\mathcal{S}\text{et}^{(\mathcal{V}^f)^\text{op}}$ and $\mathcal{S}\text{et}^{(\mathcal{VI})^\text{op}}$ be the categories of contravariant functors from $\mathcal{V}^f$ and $\mathcal{VI}$ to the category of sets.
    \item Let $\mathcal{F}\text{in}^{(\mathcal{V}^f)^\text{op}}$ and $\mathcal{F}\text{in}^{(\mathcal{VI})^\text{op}}$ be the categories of contravariant functors from $\mathcal{V}^f$ and $\mathcal{VI}$ to the category of finite sets.
    \item Let $\mathcal{P}\text{fin}^{(\mathcal{V}^f)^\text{op}}$ and $\mathcal{P}\text{fin}^{(\mathcal{VI})^\text{op}}$ be the categories of contravariant functors from $\mathcal{V}^f$ and $\mathcal{VI}$ to the category of profinite sets.
\end{enumerate}

\end{definition}

\begin{definition}
For $K$ in $\K$, we define $\text{Hom}_{\mathcal{K}\text{f.g.}}(K,H^*(\_))$ in $\mathcal{S}\text{et}^{(\mathcal{VI})^\text{op}}$ which maps $V\in\mathcal{VI}$ to the set of morphisms $\phi$ from $K$ to $H^*(V)$ such that $\phi$ makes $H^*(V)$ a finitely generated $K$-module.
\end{definition}

\begin{remark}
If $K$ is finitely generated as an algebra over $\A$, then, for all finite dimensional vector spaces $V$, $\text{Hom}_\mathcal{K}(K,H^*(V))$ is finite. Thus $\text{Hom}_\mathcal{K}(K,H^*(\_))$ and $\text{Hom}_{\mathcal{K}\text{f.g.}}(K,H^*(\_))$ are respectively in $\mathcal{F}\text{in}^{(\mathcal{V}^f)^\text{op}}$ and $\mathcal{F}\text{in}^{(\mathcal{VI})^\text{op}}$.\\
In particular, this is the case for $K$ noetherian.
\end{remark}

\subsection{Regular elements of an object in $\SVF$}

In this sub-section we recall the notion of a regular element of an object of $\SVF$. The interest of this notion is that, for $F$ in $\SVF$ satisfying a noetherianity condition introduced in \cite{HLS2}, the regular elements of $F$ define an object in the much simpler category $\SVI$. Moreover, if $F\cong\HomK(K,H^*(\_))$ for $K$ a noetherian unstable algebra, the regular elements of $F$ are given by $\HomKf(H,H^*(\_))$ and we can show that $F$ is fully determined by its regular elements. This will allow us to make the connection between the study of the $\Nil$-localisation of noetherian unstable algebras, and the study of functors in $\FVI$.

%\begin{definition}\cite[Partie II.2]{HLS2}
%We will say that an unstable algebra $K$ is of transcendence degree $d\in\N\cup\{\infty\}$ if $d$ is the superior borne of cardinals of finite sets of homogeneous elements of $K$ which are algebraically independents.
%\end{definition}

%\begin{remark}
%In particular, if $K$ is noetherian, $K$ is of finite transcendence degree.
%\end{remark}

%We will see that for $K$ of transcendence degree $d$, $\HomKf(K,H^*(V))$ is trivial for $\dim(V)>d$.

\begin{proposition}\cite[Proposition-Definition 5.1]{HLS2}\label{imp}
Let $G\in\mathcal{S}\text{et}^{(\mathcal{V}^f)^\text{op}}$, $V\in\E$ and $s\in G(V)$. Then, there exists a unique sub-vector space $U$ of $V$, denoted by $\text{ker}(s)$, such that: 
\begin{enumerate}
    \item For all $t\in G(W)$ and all morphism $\alpha\ :\ V\rightarrow W$ such that $s=G(\alpha)(t)$, $\text{ker}(\alpha)\subset U$.
    \item There exists $W_0$ in $(\mathcal{V}^f)^\text{op}$, $t_0\in G(W_0)$ and $\alpha_0\ :\ V\rightarrow W_0$ such that $s=G(\alpha_0)(t)$ and $\text{ker}(\alpha_0)=U$.
    \item There exists $t_0\in G(V/U)$ such that $s=G(\pi)(t_0)$, where $\pi$ is the projection of $V$ onto $V/U$.
\end{enumerate}
\end{proposition}

\begin{definition}
Let $G\in\mathcal{S}\text{et}^{(\mathcal{V}^f)^\text{op}}$, $V\in\E$ and $s\in G(V)$. We say that $s$ is regular if $\ker(s)=0$. 

Let $\text{reg}(G)(V):=\{x\in G(V)\ ;\ \ker(x)=0\}$.
\end{definition}

We recall the definition of a noetherian functor from \cite{HLS2}.

\begin{definition}
Let $F$ be in $\PFVF$, we say that $F$ is noetherian if it satisfies the following:

\begin{enumerate}
    \item $F\in\FVF$,
    \item there exists an integer $d$ such that $F(V)=\emptyset$ for $\dim(V)> d$,
    \item for all $V\in\E$ and $s\in F(V)$ and for all morphisms $\alpha$ which takes values in $V$, $\text{ker}(F(\alpha)s)=\alpha^{-1}(\text{ker}(s))$.
\end{enumerate} 
\end{definition}

\begin{proposition} \cite[Theorem 7.1]{HLS2}\label{noethnoeth}
\begin{enumerate}
    \item If $K\in\K$ is noetherian, $\text{Hom}_\mathcal{K}(K,H^*(\_))$ is noetherian. 
    \item If $F\in\SVF$ is noetherian, then $m\circ\mathcal{L}(F)\in\K$ is noetherian.
\end{enumerate}

\end{proposition}

\begin{lemma}
Let $F$ be a noetherian functor, then, $\text{reg}(F)$ is an object in $\mathcal{F}\text{in}^{(\mathcal{VI})^\text{op}}$. 
\end{lemma}

\begin{proof}
We only have to prove that if $x\in \text{reg}(F)(V)$, and if $\alpha$ from $E$ to $V$ is an injection, $\alpha^*x$ is regular. But, since $F$ is supposed noetherian, $\ker(\alpha^*x)=\alpha^{-1}(\ker(x))$ which is equal to $\{0\}$ since $x$ is regular and $\alpha$ is an injection.
\end{proof}

\begin{remark}
$\reg$ is not a functor, since the image under a natural transformation of a regular element is not necessarily regular. For example, fix $V$ and $U$ where $U\neq 0$ is a sub-vector space of $V$, and consider the natural transformation from $\HomF(\_,V)\overset{\pi}{\rightarrow}\HomF(\_,V/U)$ induced by the projection from $V$ to $V/U$. Then, $\id_V$ is regular in $\HomF(V,V)$ but $\ker(\pi\circ\id_V)=U$ in $\HomF(V,V/U)$.
\end{remark}

\begin{proposition}\cite[Proposition 4.8]{DW1}\label{rugugu}
For $K$ a noetherian unstable algebra, we have a natural isomorphism $\text{reg}(\text{Hom}_\mathcal{K}(K,H^*(\_)))\cong\text{Hom}_{\mathcal{K}\text{f.g.}}(K,H^*(\_))$.
\end{proposition}

Let us denote by $\mathcal{O}$ the forgetful functor from $\SVF$ to $\SVI$ induced by restriction to $(\mathcal{VI})^\text{op}$. We define the left Kan extension of $\mathcal{O}$ along the identity.

\begin{proposition}\label{kanus}
The functor $\mathcal{O}$ has a left adjoint, which maps a functor $F\in \SVI$ to $\Tilde{F}$, which is defined by $\Tilde{F}(V):=\bigsqcup\limits_{U\in\textbf{S}(V)}F(V/U)$, where $\textbf{S}(V)$ is the set of sub-vector spaces of $V$.

\end{proposition}

\begin{proof}
Let us first prove that $\Tilde{F}$ defines an object in $\SVF$. To a morphism $\alpha\ :\ V\rightarrow W$ we associate $\Tilde{F}(\alpha)$ in the following way: for $x\in F(W/U)$ we consider $\pi\circ\alpha$ where $\pi$ is the projection from $W$ to $W/U$, we factorise $\pi\circ\alpha$ as $\Tilde{\alpha}\circ\psi$ with $\psi$ the projection from $V$ to $V/\ker(\pi\circ\alpha)$, then $\Tilde{\alpha}$ is injective and we define $\Tilde{F}(\alpha)(x)=F(\Tilde{\alpha})(x)\in F(V/\ker(\pi\circ\alpha))$.\\
 
Let $\phi$ be a natural transformation from $\Tilde{F}$ to $A$, where $A$ is an object in $\SVF$. For all $V\in\E$, we have a morphism in $\mathcal{S}et$, $$\phi_V\ : \ \bigsqcup\limits_{U\in\textbf{S}(V)}F(V/U)\rightarrow A(V) .$$ Then,
$$\phi_V|_{F(V)}$$ induces a morphism $\iota(\phi)_V$ from $F(V)$ to $\mathcal{O}(A)(V)$. $\iota(\phi)$ is a natural transformation in $\SVI$. \\
Conversely, for $\phi$ a natural transformation from $F$ to $\mathcal{O}(A)$, we define $$\gamma(\phi)_V\ :\ \Tilde{F}(V)\rightarrow A(V)$$ in the following way: for $x\in F(V/U)$ we define $\gamma(\phi)_V(x):=A(\pi)(\phi_{V/U}(x))$, for $\pi$ the projection from $V$ to $V/U$. Then, $\gamma(\phi)$ is a natural transformation in $\SVF$, and $\gamma$ and $\iota$ are mutually inverse.
\end{proof}

\begin{remark}
For $F\in\SVI$, $\Tilde{F}$ always satisfies the second condition in the definition of a noetherian functor. Moreover, $\Tilde{F}$ is in $\FVF$ if and only if $F\in\FVI$ and there is an integer $d$ such that $F(V)=\emptyset$ for $\dim(V)\geq d$, hence $\Tilde{F}$ is noetherian if and only if $F\in\FVI$ and there exists $d\in\N$ such that $F(V)$ is empty for $\dim(V)$ greater than $d$.
\end{remark}

\begin{proposition}\label{regkan}
For $F$ a noetherian functor, $\widetilde{\text{reg}(F)}\cong F$.
\end{proposition}

\begin{proof}
Let $F$ be a noetherian functor. For $V$ a finite dimensional vector space, we define the following morphism from $\widetilde{reg(F)}(V)$ to $F(V)$:
$\begin{array}{ccc}
   \widetilde{reg(F)}(V)  &  \rightarrow & F(V)\\
   x\in reg(F)(V/U)  & \mapsto & F(\pi)(x),
\end{array}$
where $\pi$ is the projection from $V$ to $V/U$. Then, Proposition \ref{imp} implies that this morphism is an isomorphism, and the fact that $F$ is noetherian implies that it is natural in $V$.
\end{proof}

\begin{corollary}\label{Hfg}
For $K$ a noetherian unstable algebra, there is a natural bijection with respect to $V$, $$\text{Hom}_\mathcal{K}(K,H^*(V))\cong \bigsqcup\limits_{U\in\textbf{S}(V)}\text{Hom}_{\mathcal{K}f.g.}(K,H^*(V/U)).$$
\end{corollary}

\subsection{Central elements of a noetherian algebra}

In the last subsection, we explained why, for $K$ a noetherian unstable algebra, $\HomK(K,H^*(\_))$ is fully determined by its regular elements, given by $\HomKf(K,H^*(\_))$. In this section, we prove that for $K$ noetherian, one can deduce $\textbf{C}(K)$ from the central elements $(V,\phi)$ of $K$ such that $\phi$ is regular.\\

We recall the definition of the centre of a noetherian unstable algebra from \cite{heard2020topological} which is in some sense the maximal regular central element of $K$.\\

We state the following results for noetherian algebras, since it is in this case that we will use those, but all of them are true if we only suppose that $\mathcal{Q}(K)$ is locally finite.

\begin{proposition}
Let $K$ be a noetherian unstable algebra  and let $\phi$ be in $\HomK(K,H^*(E))$, then if $(E,\phi)$ is central, let $\phi_0$ be the only element in $\HomK(K,H^*(E/\ker(\phi)))$ such that $\phi=\pi^*\phi_0$. Then, $\phi_0$ is regular and $(E/\text{ker}(\phi),\phi_0)$ is central.\\

Conversely, if $(E/U,\phi)$ is central with $\phi$ regular, $(E,\pi^*\circ\phi)$ is central. Thus, $(E,\phi)$ is central if and only if $(E/\ker(\phi),\phi_0)$ is central.
\end{proposition}

\begin{proof}

We consider $s$ a section of $\pi$, the projection from $E$ to $E/\text{ker}(\phi)$, then $\phi_0=s^*\circ\phi$. Then, by Proposition \ref{incl}, if $(E,\phi)$ is central, so is $(E/\text{ker}(\phi),\phi_0)$. The converse is a direct conseequence of Proposition \ref{incl}.

\end{proof}

\begin{lemma} \cite[Lemma 4.6]{DW1}\label{boxplus}
Let $K$ be a noetherian unstable algebra, let $f$ be in $\HomK(K,H^*(E))$ and $(C,g)\in\textbf{C}(K)$, then there exists a unique pair $(E\oplus C,f\boxplus g)$ with $f\boxplus g\in\HomK(K,H^*(E\oplus C))$ such that $f\boxplus g$ composed respectively with the projections on $H^*(E)$ and $H^*(C)$ gives $f$ and $g$. 
$$\xymatrix{ & & H^*(C)\\
K\ar[rru]^{g}\ar[rr]^{f\boxplus g}\ar[rrd]_f & & H^*(C)\otimes H^*(E)\ar[u]_-{\id\otimes\epsilon_{H^*(E)}}\ar[d]^-{\epsilon_{H^*(C)}\otimes\id}  \\
& & H^*(E).}$$
\end{lemma}

\begin{remark}
In Lemma \ref{boxplus}, even if $f$ and $g$ are regular, $f\boxplus g$ is not regular in general.
\end{remark}

We recall now the definition of the centre of a noetherian algebra.

\begin{definition}\cite[Definition 3.9]{heard2020topological}
Let $K$ be a noetherian unstable algebra, let $f$ be in $\HomK(K,H^*(E))$ and $(C,g)\in\textbf{C}(K)$, we define $(E\circ C,\sigma(f,g))$ by $E\circ C=(E\oplus C)/\text{ker}(f\boxplus g)$, and $\sigma(f,g)\ :\ K\rightarrow E\circ C$ such that the composition of $\sigma(f,g)$ by $H^*(\pi_{(E\oplus C)/\text{ker}(f\boxplus g)})$ is equal to $f\boxplus g$, for $\pi_{(E\oplus C)/\text{ker}(f\boxplus g)}$ the projection on $(E\oplus C)/\text{ker}(f\boxplus g)$.
\end{definition}

\begin{proposition} \cite[Corollary 3.11]{heard2020topological}\label{pluscentr}
Let $K$ be a noetherian unstable algebra and $(C,g)$ and $(E,f)$ be two central elements of $K$. Then $(E\oplus C,f\boxplus g)$ and $(E\circ C,\sigma(f,g))$ are central.
\end{proposition}

\begin{theorem}\cite[Theorem 3.13]{heard2020topological}\label{defcentre}
Let $K$ be a connected, noetherian, unstable algebra. Then, up to isomorphism, there is a unique couple $(C,\gamma)\in\textbf{C}(K)$ which is regular and satisfy the following: for all central element $(E,f)$ with $f$ regular, there exists an injection $\iota$ from $E$ to $C$, such that $f=\iota^*\gamma$. We will call $(C,\gamma)$ the centre of $K$.
\end{theorem}

The idea of the proof is to consider $(C,\gamma)$ a regular central element of maximal dimension. Then, show that for $(E,f)$ regular and central, either $(E,f)=(E,\iota^*\gamma)$ for $\iota$ an injection from $E$ to $C$, or the dimension of $E\circ C$ is greater than that of $C$ which contradicts the assumption.

\begin{corollary}
For $E\in\E$ and $f\in\HomK(K,H^*(V))$, $(E,f)$ is central if and only if there is a morphism $\alpha$ from $E$ to $C$, with $(C,\gamma)$ the centre of $K$, such that $f=\alpha^*\gamma$.
\end{corollary}

\begin{proof}
It is a direct consequence of Proposition \ref{incl} and Theorem \ref{defcentre}.
\end{proof}

\subsection{The shift functor}

We want to be able to discuss centrality away from $\Nil$, and define accordingly a notion of central elements for objects in $\SVI$. The aim of this subsection is to define, for $F\in\SVI$, $V\in\E$ and $\phi\in F(V)$ a functor $\sigma_{(V,\phi)}F$, in such a way that for $F=\HomKf(K,H^*(\_))$, $\sigma_{(V,\phi)}F\cong\HomKf(T_{(V,\phi)}(K),H^*(\_))$.\\

We start by defining functors $\Sigma_V$ from $\SVF$ to $\SVF$ such that, for all unstable algebras $K$, $\Sigma_V\HomK(K,H^*(\_))\cong\HomK(T_V(K),H^*(\_))$, and such that $\Sigma_{V}\HomK(K,H^*(\_))$ naturally decomposes in $\Sigma_{(V,\phi)}\HomK(K,H^*(\_))\cong\HomK(T_{(V,\phi)}(K),H^*(\_))$, with $\phi$ running through $\HomK(K,H^*(V))$.\\
 
Then, we define a shift functor $\sigma_V $ in the category $\SVI$ and we identify the desired functors $\sigma_{(V,\phi)}F$ as the connected components of $\sigma_V F$. It is worth pointing out that since the decomposition of $\HomKf(T_V(K),H^*(\_))$ is induced by elements $(V,\phi)$ with $\phi\in\HomK(K,H^*(V))$ (and not only elements in $\HomKf(K,H^*(V))$), $\sigma_V\HomKf(K,H^*(\_))$ is only a sub functor of $\HomKf(T_V(K),H^*(\_))$.

\begin{definition}
For $V\in\E$, define the functor $\Sigma_V$ from $\SVF$ to $\SVF$ by $$\Sigma_V F(W)=F(V\oplus W),$$ and $\Sigma_VF(\alpha)=F(\id_V\oplus\alpha)$, for $F\in\SVF$, $W\in\E$ and $\alpha$ a morphism in $\E$.
\end{definition}

\begin{lemma}\label{cococo}
For $K\in\K$ and $V\in\E$, there is an isomorphism natural both in $K$ and $V$, $\HomK(T_V(K),H^*(\_))\cong \Sigma_V \HomK(K,H^*(\_))$.
\end{lemma}

\begin{proof}
This is a direct consequence of the definition of $T_V$ as the adjoint of the tensor product by $H^*(V)$ and of the natural isomorphism $H^*(V)\otimes H^*(W)\cong H^*(V\oplus W)$.
\end{proof}

\begin{definition}
For $F\in\SVF$, $V\in\E$ and $\phi\in F(V)$, we consider $\Sigma_{(V,\phi)}F(W)$ to be the fibre over $\{\phi\}$ of the morphism $\Sigma_V F(W)\rightarrow \Sigma_VF(0)\cong F(V)$ induced by the injection from $0$ to $W$.
\end{definition}

\begin{lemma}\label{gozagip}
Let $K$ be an unstable algebra, $V\in\E$ and $\phi\in\HomK(K,H^*(V))$. There is a natural isomorphism $\Sigma_{(V,\phi)}\HomK(K,H^*(W))\cong\HomK(T_{(V,\phi)}(K),H^*(W))$.
\end{lemma}

\begin{proof}
We have a commutative diagram $$\xymatrix{\Sigma_V\HomK(K,H^*(W))\ar[r]^\cong\ar[d] & \HomK(T_V(K),H^*(W))\ar[d]\\
\HomK(K,H^*(V))\ar[r]^\cong &\HomK(T_V(K),\Fp),}$$ where the vertical maps are induced by the injection from $\{0\}$ to $W$, and the horizontal ones are given by the natural isomorphism of Lemma \ref{cococo}. By definition, $\Sigma_{(V,\phi)}\HomK(K,H^*(W))$ is the fibre over $\phi$ of the left map, and by construction $\HomK(T_{(V,\phi)}(K),H^*(W))$ is the fibre over the adjoint of $\phi$ of the right one. This concludes the proof. 
\end{proof}

\begin{proposition}\label{noeeth}
If $F$ is noetherian, for $x$ an element in $F(V)$, $\Sigma_V F$ and $\Sigma_{(V,x)}F$ are also noetherian.
\end{proposition}

\begin{proof}
If $F$ is finite, $\Sigma_V F$ and $\Sigma_{(V,x)}F$ are obviously finite, we only have to prove the second condition. 
 Let $a$ be an element in $\sigma_{V}F(W)$, then $a\in F(V\oplus W)$. There is an ambiguity in considering $\ker(a)$, since $a$ can be viewed either as an element in $F$ or $\Sigma_VF$. Let $\ker^1(a)$ denote the kernel of $a$ considered as an element in $F$, and $\ker^2(a)$ its kernel in $\Sigma_V F$. We notice that $\ker^1(a)$ is a sub-vector space of $V\oplus W$, whereas $\ker^2(a)$ is a sub-space of $W$. Then, $\ker^2(a)=\ker^1(a)\cap W$. Indeed, for $t$ regular such that $a=\pi^*t$ for $\pi$ the projection from $V\oplus W$ to $(V\oplus W)/\ker^1(a)$, $\pi$ factorises as $\pi_2\circ \pi_1$, for $\pi_1$ the projection from $V\oplus W$ to $V\oplus (W/(\ker^1(a)\cap W))$ and $\pi_2$ the projection from $V\oplus (W/(\ker^1(a)\cap W))$ to $(V\oplus W)/\ker^1(a)$, then $a=\Sigma_VF(\pi')(\pi_2^*t)$, for $\pi'$ the projection from $W$ to $W/(\ker^1(a)\cap W)$. Therefore, $\ker^1(a)\cap W\subset\ker^2(a)$.\\
 
 Conversely, if $\pi$ is now the projection from $W$ to $W/\ker^2(a)$, let $t\in F(V\oplus (W/\ker^2(a)))$ such that $a=\Sigma_V F(\pi)(t)$. Then, $a=(\id_V\oplus\pi)^*t$, hence $\ker^2(a)\subset \ker^1(a)$. Since, we also have $\ker^2(a)\subset W$, $\ker^2(a)=\ker^1(a)\cap W$. Then, for $\alpha$ a morphism from a finite dimensional vector space $U$ to $W$, $\ker^2(\Sigma_VF(\alpha)(a))=\ker^1((\id_V\oplus \alpha)^*a)\cap W$ which is equal to $(\id_V\oplus\alpha)^{-1}(\ker^1(a))\cap W)$, since $F$ is noetherian. Hence, it is equal to $\alpha^{-1}(\ker^1(a)\cap W)=\alpha^{-1}(\ker^2(a))$. This proves that $\Sigma_V F$ is noetherian. The proof for $\Sigma_{(V,x)}F $ is similar.

\end{proof}

\begin{corollary}
If $K$ is noetherian and $nil$-closed, $T_V(K)$ and $T_{(V,\phi)}K$ are noetherian, for $V\in\E$ and $\phi\in\HomK(K,H^*(V))$. 
\end{corollary}

\begin{proof}
By Proposition \ref{noeeth}, $g(T_V(K))$ and $g(T_{(V,\phi)}(K))$ are noetherian, then by Proposition \ref{noethnoeth} $l_1(T_V(K))\cong (m\circ\mathcal{L}\circ g)(T_V(K))$ and $l_1(T_{(V,\phi)}(K))$ are noetherian. Therefore, since by Corollary \ref{tutut} and Lemma \ref{upussu} $T_V(K)$ and $T_{(V,\phi)}(K)$ are $nil$-closed, they are noetherian.
\end{proof}

We want to identify a shift functor $\sigma_V$ from $\SVI$ to $\SVI$, in so that $\sigma_V\HomKf(K,H^*(\_))$ captures the behaviour of $T_V(K)$ away from nilpotent objects for $K$ noetherian, and such that for $F$ in $\SVI$, $\sigma_{V}F$ comes with a decomposition in $\sigma_{(V,\phi)}F$ with $\phi$ running through $F(V)$. \\

In order to define this shift functor, we have to discuss pushouts in $\E$ and $\mathcal{VI}$. It is worth noticing that the pull-back in $\E$ of a diagram whose morphisms are in $\mathcal{VI}$ is also a pull-back in $\mathcal{VI}$.

\begin{remark}
Pushouts usually don't exist in $\mathcal{VI}$, for example, if one consider the following diagram $$\xymatrix{0\ar[r]\ar[d] & V\\V & \ ,}$$
the pushout ``should be" $V\oplus V$, but since non injective morphisms are not in $\mathcal{VI}$ the commutative square $$\xymatrix{0\ar[r]\ar[d] & V\ar[d]\\V \ar[r] & V}$$ does not give rise to a morphism from $V\oplus V$ to $V$ satisfying the universal property of the pushout. So, in the following, we will consider pushouts in $\E$ of diagrams in $\mathcal{VI}$. 
\end{remark}

The following simple lemma deals with this.

\begin{lemma}\label{pushofpull}
Consider a pushout square in $\E$: $$\xymatrix{M\ar[r]^\omega\ar[d]_\nu & W\ar[d]\\ V\ar[r] & P.}$$ Suppose that $\nu$ and $\omega$ are injections. Then, for all pullback squares $$\xymatrix{M\ar[r]^\omega\ar[d]_\nu & W\ar[d]\\ V\ar[r] & N,}$$ such that the morphisms from $V$ and $W$ to $N$ are injections, the morphism from $P$ to $N$ induced by the universal property of the pushout in $\E$, is a morphism in $\mathcal{VI}$.
\end{lemma}

Let us also recall the pasting law for pullbacks: 

\begin{lemma}\cite[Proposition 11.10]{cat}\label{pastlaw}
We consider a diagram of the following shape in a category $\mathcal{C}$:
$$\xymatrix{A\ar[r]\ar[d] & B\ar[r]\ar[d] & C\ar[d]\\
D\ar[r] & E\ar[r] & F.}$$ If the right square is a pullback square, then the outer square is a pullback square if and only if the left square is.
\end{lemma}

\begin{definition}
For $V$ and $W$ some fixed objects of $\E$, let $B(V,W)$ be the set of triples $(M,\nu,\omega)$, where $M\in\E$ and $\nu$ and $\omega$ are morphisms from $M$ to $V$ and $W$ in $\mathcal{VI}$. Let also $\equiv$ be the relation on $B(V,W)$ defined by $(M,\nu,\omega)\equiv(M',\nu',\omega')$ if there exists an isomorphism $\mu$ from $M'$ to $M$ such that $\nu'=\nu\circ\mu$ and $\omega'=\omega\circ\mu$.
\end{definition}

The following lemma is obvious.

\begin{lemma}
$\equiv$ is an equivalence relation on $B(V,W)$.
\end{lemma}

\begin{definition}
For $V$ and $W$ two objects of $\mathcal{VI}$, let $\textbf{B}(V,W)$ be the set of equivalence classes for $\equiv$ in $\text{B}(V,W)$. For $\left[M,\nu,\omega\right]\in\textbf{B}(V,W)$, let $V\oplus_{\nu,\omega} W$ denote the pushout of the following diagram in $\E$: $$\xymatrix{M\ar[r]^{\omega}\ar[d]_{\nu} & W\\
V &\ ,  }$$ let also denote by $\iota_V^{V\oplus_{\nu,\omega} W}$ and $\iota_W^{V\oplus_{\nu,\omega} W}$ the induced injections from $V$ and $W$ to $V\oplus_{\nu,\omega} W$. When there is no ambiguity, we will denote them by $\iota_V$ and $\iota_W$.
\end{definition}

\begin{remark}
 $V\oplus_{\nu,\omega} W$ does not depend on the choice of $(M,\nu,\omega)\in\left[M,\nu,\omega\right]$.
\end{remark}

\begin{proposition}
$\textbf{B}(\_,\_)$ is a bifunctor on $\mathcal{VI}^{op}$.
\end{proposition}

\begin{proof}
Let $\alpha\ :\ V'\hookrightarrow V$ and $\beta\ : W'\hookrightarrow W$ be morphisms in $\mathcal{VI}$, for $\left[M,\nu,\omega\right]\in\textbf{B}(V,W)$ we define $\textbf{B}(\alpha,\beta)(\left[M,\nu,\omega\right])\in\textbf{B}(V',W')$ in the following way. We consider the following diagram, where $M_V$, $M_W$ and $M'$ are defined by pullback. $$\xymatrix{M'\ar[r]\ar[d] & M_W\ar[r]\ar[d] & W'\ar[d]^\beta\\
M_V\ar[r]\ar[d] & M\ar[r]^\omega\ar[d]_\nu & W\\
V'\ar[r]_\alpha & V & \ .}$$ Then, $\nu'$ and $\omega'$ are given by the compositions $M'\rightarrow M_V\rightarrow V'$ and $M'\rightarrow M_W\rightarrow W'$.
\end{proof}

\begin{remark}
By the pasting law of pullbacks (lemma \ref{pastlaw}, that we used two times) we can show that the outer square of the following diagram is a pullback square: $$\xymatrix{M'\ar[r]\ar[d] & M_W\ar[r]\ar[d] & W'\ar[d]\\
M_V\ar[r]\ar[d] & M\ar[r]^\omega\ar[d]_\nu & W\ar[d]\\
V'\ar[r] & V\ar[r] & V\oplus_{\nu,\omega} W.}$$
\end{remark}

\begin{definition}
Let $\alpha\ :\ V'\hookrightarrow V$ and $\beta\ : W'\hookrightarrow W$ be injections in $\mathcal{VI}$, for $\left[M,\nu,\omega\right]\in \textbf{B}(V,W)$ and for $\left[M',\nu',\omega'\right]:=\textbf{B}(\alpha,\beta)(\left[M,\nu,\omega\right])\in\textbf{B}(V',W')$ we define $$\alpha\oplus_{\nu,\omega}\beta\ :\ V'\oplus_{\nu',\omega'}W'\hookrightarrow V\oplus_{\nu,\omega}W,$$ the injective (cf lemma \ref{pushofpull}) morphism induced by the universal property of the pushout in $\E$ of the following pullback square: $$\xymatrix{M'\ar[r]^{\nu'}\ar[d]_{\omega '} & W'\ar[d]^-{\omega\circ\beta}\\
V'\ar[r]_-{\nu\circ\alpha} & V\oplus_{\nu,\omega} W.}$$
\end{definition}

\begin{definition}
 
For $F\in\SVI$ and $V$ and $W$ two objects in $\mathcal{VI}$, we define $$\sigma_V F(W):=\bigsqcup\limits_{\left[M,\nu,\omega\right]\in\textbf{B}(V,W)}F(V\oplus_{\nu,\omega} W).$$
\end{definition}

Even though the roles of $V$ and $W$ are symmetric in the definition of $\sigma_VF(W)$, we use an asymmetric notation to reflect that of $T_V(K)$.

\begin{proposition}\label{sommetronqu}
For $F\in\SVI$, $\sigma_{\_}F$ is a bifunctor on $\mathcal{VI}^{op}$.
\end{proposition}

\begin{proof}
For $\alpha\ :\ V'\hookrightarrow V$ and $\beta\ : W'\hookrightarrow W$ two morphisms in $\mathcal{VI}$ we have to define $\sigma_\alpha(\beta)$ from $\sigma_VF(W)$ to $\sigma_{V'}F(W')$. For $x\in F(V\oplus_{\nu,\omega}W)$ let $\sigma_\alpha F(\beta)(x)=F(\alpha\oplus_{\nu,\omega}\beta)(x)\in F(V'\oplus_{\nu',\omega'}W')$. Then, the fact that $\sigma_\_F$ is a bifunctor comes from the fact that $\id_V\oplus_{\nu,\omega}\id_W=\id_{V\oplus_{\nu,\omega}W}$ and that for $\alpha'\ :\ V"\rightarrow V'$ and $\beta'\ :\ W"\rightarrow W'$, $(\alpha\oplus_{\nu,\omega}\beta)\circ(\alpha'\oplus_{\nu',\omega'}\beta')=(\alpha\circ\alpha')\oplus_{\nu,\omega}(\beta\circ\beta')$.

\end{proof}

\begin{corollary}
For $F\in\SVI$ and $V$ and $W$ some fixed objects in $\mathcal{VI}$, $\sigma_VF$ and $\sigma_\_F(W)$ are objects in $\SVI$.
\end{corollary}

Let us now define the analogue of $T_{(V,\phi)}$ in the category $\SVI$.

\begin{lemma}
For $F\in\SVI$ and $V$ and $W$ in $\mathcal{VI}$, we have natural isomorphisms $\sigma_0 F(W)\cong F(W)$ and $\sigma_V F(0)\cong F(V)$.
\end{lemma}

\begin{definition}
For $F\in\SVI$ and $V$ and $W$ two objects of $\mathcal{VI}$, the morphism from $0$ to $W$ induces a morphism $\sigma_V F(W)\rightarrow \sigma_VF(0)\cong F(V)$. For $x\in F(V)$, let $\sigma_{(V,x)}F(W)$ be the inverse image of $\{x\}$ under this morphism.
\end{definition}

\begin{lemma}
For $F\in\SVI$, $V$ and $W$ objects in $\mathcal{VI}$ and for $x\in F(V)$, $\sigma_{(V,x)}F$ is a sub-functor of $\sigma_V F$. 
\end{lemma}

\begin{proof}
We only have to prove that for $\alpha\ :\ U\rightarrow W$ a morphism in $\mathcal{VI}$ and for $y\in\sigma_{(V,x)}F(W)$, $\sigma_V F(\alpha)(y)\in\sigma_{(V,x)}F(U)$. This follow directly from the fact that the morphism $0\rightarrow W$ factorizes through $0\rightarrow U\rightarrow W$.
\end{proof}

As we stated in the beginning of this sub-section, $$\sigma_V\HomKf(K,H^*(\_))\not\cong\HomKf(T_V(K),H^*(\_)).$$ Nonetheless, as we will see in corollary \ref{opata}, for $\phi\in\HomKf(K,H^*(V))$, we have that $$\sigma_{(V,\phi)}\HomKf(K,H^*(\_))\cong \HomKf(T_{(V,\phi)}(K),H^*(\_)).$$ 

\begin{proposition}\label{optou}
For $F\in\SVI$, $V\in\E$ and $x\in F(V)$, $\widetilde{\sigma_{(V,x)}F}\cong\Sigma_{(V,x)}\Tilde{F}$.
\end{proposition}

\begin{proof}
Let $W$ be a finite dimensional vector space. Let us first notice that both $\widetilde{\sigma_{(V,x)}F}(W)$ and $\Sigma_{(V,x)}\Tilde{F}(W)$ can be seen as sub-sets of $\Tilde{F}(V\oplus W)=\bigsqcup\limits_{H\in\textbf{S}(V\oplus W)}F((V\oplus W)/H)$. It is obvious for $\Sigma_{(V,x)}\Tilde{F}(W)$. For $\widetilde{\sigma_{(V,x)}F}(W)\subset\widetilde{\sigma_V F}(W)$, for every $V\oplus_{\nu,\omega} (W/U)$ the universal property of the product induces a projection $\pi_{\nu,\omega}$ from $V\oplus W$ to $V\oplus_{\nu,\omega} (W/U)$. It induces an isomorphism from $(V\oplus W)/\ker(\pi_{\nu,\omega})$ to $V\oplus_{\nu,\omega} (W/U)$, we identify $F(V\oplus_{\nu,\omega} (W/U))$ with\newline 
$F((V\oplus W)/\ker(\pi_{\nu,\omega}))\subset\Tilde{F}(V\oplus W)$ through this isomorphism. Then we only have to prove that those sub-sets are equals.\\

$\Sigma_{(V,x)}\Tilde{F}(W)$ identifies with the set of elements $\gamma$ in $\Tilde{F}(V\oplus W)$ such that $\Tilde{F}(\iota_V^{V\oplus W})(\gamma)=x\in F(V)$. Since $\Tilde{F}$ is noetherian and $x$ is regular, $\ker(\gamma)\cap V=\{0\}$. Then, there exists $U$ a subspace of $W$ and a class $\left[M,\nu,\omega\right]\in\textbf{B}(V,W/U)$ such that the canonical projection from $V\oplus W$ to $V\oplus_{\nu,\omega} (W/U)$ induces an isomorphism $(V\oplus W)/\ker(\gamma)\cong V\oplus_{\nu,\omega} (W/U)$. Then, up to this isomorphism $\gamma\in F(V\oplus_{\nu,\omega} (W/U))$ and $(\iota_V^{V\oplus_{\nu,\omega}(W/U)})^*\gamma=x$, so $\gamma\in \sigma_{(V,x)}F(W/U)\subset\widetilde{\sigma_{(V,x)}F}(W)$.\\

Conversely, for $\gamma\in\sigma_{(V,x)}F(W/U)\subset\widetilde{\sigma_{(V,x)}F}(W)$, there is $\left[M,\nu,\omega\right]\in\textbf{B}(V,W/U)$ such that $\gamma\in F(V\oplus_{\nu,\omega} (W/U))$ and $(\iota_V^{V\oplus_{\nu,\omega}(W/U)})^*\gamma=x$, then for $H$ the kernel of the projection from $V\oplus W$ to $V\oplus_{\nu,\omega} (W/U)$, we have the following isomorphism induced by the first isomorphism theorem $(V\oplus W)/H\cong V\oplus_{\nu,\omega} (W/U)$. Up to this isomorphism, $\gamma\in F((V\oplus W)/H)\subset\Tilde{F}(V\oplus W)$. Then, by construction and up to the isomorphism from $(V\oplus W)/H$ to $V\oplus_{\nu,\omega}(W/U)$, $\Tilde{F}(\iota_V^{V\oplus W})(\gamma)=(\iota_V^{V\oplus_{\nu,\omega}(W/U)})^*\gamma=x$. Then, $\gamma\in\Sigma_{(V,x)}\Tilde{F}(W)$.

\end{proof}

\begin{corollary}\label{opata}
For $K$ a noetherian algebra, $V\in\E$ and $\phi\in\HomKf(K,H^*(V))$, $$\sigma_{(V,\phi)}\HomKf(K,H^*(\_))\cong\HomKf(T_{(V,\phi)}(K),H^*(\_)).$$
\end{corollary}

\begin{proof}
This is a direct consequence of Lemma \ref{gozagip}, of Proposition \ref{optou} and of the fact that, if $K$ is noetherian, $T_{(V,\phi)}(K)$ is also noetherian. 
\end{proof}

\begin{remark}
Since the morphism $\sigma_VF(W)\rightarrow F(V)$ is a surjection, $\sigma_VF(W)$ is the disjoint union of the fibres over singletons in $F(V)$, we then have an isomorphism which is natural in both $K$ and $W$, $\sigma_V\HomKf(K,H^*(W))\cong\bigsqcup\limits_{\phi\in\HomKf(K,H^*(V))}\HomKf(T_{(V,\phi)}(K),H^*(W))$.
\end{remark}

\subsection{Central elements of an object of $\SVI$}

In Definition \ref{centri}, we defined the notion of a central elements of an unstable algebra $K$. In the following, we define a notion of centrality "away from $\Nil$".

\begin{definition}
For $K\in\K$, $V\in\E$ and $\phi\in\HomK(K,H^*(V))$, we will say that $(V,\phi)$ is central away from $\Nil$, if $g(\rho_{K,(V,\phi)})$ (or equivalently $f(\rho_{K,(V,\phi)})$) is an isomorphism.
\end{definition}

\begin{remark}
Since, by lemma \ref{upussu}, for $K$ $nil$-closed, $T_{(V,\phi)}(K)$ is also $nil$-closed, if $K$ is $nil$-closed $(V,\phi)$ is central away from $\Nil$ if and only if it is central.
\end{remark}

The centrality away from $\Nil$ can be characterised by properties of the functor\newline $\HomKf(K,H^*(\_))\in\SVI$. For $F\in\SVI$, $V\in\E$ and $\phi\in F(V)$, we want to define $\rho_{F,(V,\phi)}\ :\ \sigma_{(V,\phi)}F\rightarrow F$ in such a way that, when $F\cong\HomKf(K,H^*(\_))$ where $K$ is a noetherian unstable algebra, $\rho_{F,(V,\phi)}$ is an isomorphism if and only if $g(\rho_{K,(V,\phi)})$ is.

\begin{definition}
For $F\in\SVI$, $V$ and $W$ two objects of $\mathcal{VI}$, let $\rho_{F,V}\ :\ \sigma_VF\rightarrow F$ be the natural transformation induced by the morphism $0\rightarrow V$ and the natural isomorphism $\sigma_0 F\cong F$.\\

For $x\in F(V)$, we also define $$\rho_{F,(V,x)}\ :\ \sigma_{(V,x)}F\rightarrow F$$ as the restriction of $\rho_{F,V}$ to $\sigma_{(V,x)}F$.\\

If $\rho_{F,(V,x)}$ is an isomorphism we will say that $(V,x)$ is a central element of $F$. We will denote by $\textbf{C}(F)$ the set of central elements of $F$.
\end{definition}

\begin{theorem}\label{principal0}
For $K$ a noetherian unstable algebra and $\phi\in \HomKf(K,H^*(V))$, $(V,\phi)$ is central away from $\Nil$ for $K$ if and only if it is central for $\HomKf(K,H^*(\_))$.
\end{theorem}

\begin{proof}
By construction, $$g(\rho_{K,(V,\phi)})\ :\ \HomK(T_{(V,\phi)}(K),H^*(\_))\cong\Sigma_{(V,\phi)}g(K)\rightarrow g(K),$$ is the morphism which sends $\delta\in\Sigma_{(V,\phi)}\HomK(K,H^*(W))$ to $(\iota_{W}^{V\oplus W})^*\delta\in\HomK(K,H^*(W))$. Then, for $F_K=\HomKf(K,H^*(\_))$ and for $\phi$ regular, by construction the isomorphisms $\widetilde{\sigma_{(V,\phi)}F_K}(W)\cong\Sigma_{(V,\phi)}\HomK(K,H^*(W))$ and $\Tilde{F}_K(W)\cong\HomK(K,H^*(W))$ (Propositions \ref{regkan} and \ref{rugugu}), fit into the following commutative diagram, which is natural with respect to $K$:
$$\xymatrix{\Sigma_{(V,\phi)}\HomK(K,H^*(\_))\ar[rr]^-{g(\rho_{K,(V,\phi)})} & & \HomK(K,H^*(\_)) \\
\widetilde{\sigma_{(V,\phi)}F_K}\ar[u]^\cong\ar[rr]_-{\widetilde{\rho_{F,(V,\phi)}}} & &\Tilde{F}_K\ar[u]^\cong.}$$
This concludes the proof.
\end{proof}

\begin{definition}
For $F\in\SVI$, we will say that $F$ is connected if $F(0)$ is reduced to one element. 

\end{definition}

\begin{remark}
For $F\in\SVI$, $V\in\mathcal{VI}$ and $x\in F(V)$, $\sigma_{(V,x)}F(0)=\{x\}$, therefore it is connected. For $F$ not connected, $\textbf{C}(F)=\emptyset$.
\end{remark}

We give an alternative criterion for the centrality of $(V,x)$. This alternative criterion will often prove to be easier to check, when we have to prove the centrality of a given element.

\begin{lemma}\label{altedef}
Let $F\in\SVI$, $V\in\E$ and $x\in F(V)$. Then, $(V,x)$ is central if for all $W\in\E$ and $y\in F(W)$ there is a unique class $\left[M,\nu,\omega\right]\in\textbf{B}(V,W)$ and a unique $b\in F(V\oplus_{\nu,\omega} W)$ such that $\iota_V^*b=x$ and $\iota_W^*y$.
\end{lemma}

\begin{proof}
By construction, $\rho_{F,V}$ is the morphism which maps $b\in F(V\oplus_{\nu,\omega} W)$ to $\iota^*_Wb$, then $b\in\sigma_{F,(V,x)}F(W)$ if and only if $\iota^*_Vb=x$. So, $y$ has a unique element in its inverse image under $\rho_{F,(V,x)}$ if there is a unique $\left[M,\nu,\omega\right]\in\textbf{B}(V,W)$ and a single $b\in F(V\oplus_{\nu,\omega} W)$ such that $\iota_V^*b=x$ and $\iota_W^*b=y$. 
\end{proof}

In the rest of this subsection, we want to prove an analogue of Theorem \ref{defcentre} for objects in $\SVI$. Namely, we want to prove that, for $F\in\SVI$ such that $F(V)=\emptyset$ when $\dim(V)$ is greater than some integer $d$, there is, up to isomorphism, a unique maximal central element $(C,c)$ in the following sense, for all central element $(V,x)$ there is an injection $\iota$ from $V$ to $C$ such that $x=\iota^*c$.

\begin{lemma}\label{injcentr}
For $F\in\SVI$, $V$ an object in $\E$, $x\in F(V)$ and $\alpha\ :\ T\rightarrow V$ a morphism in $\mathcal{VI}$, $\alpha^*x$ is always in the image of $\rho_{F,(V,x)}$. Furthermore, the morphism $\iota_V^{T\oplus_{\id_T,\alpha} V}$ is inversible and $((\iota_V^{T\oplus_{\id_T,\alpha} V})^{-1})^*x$ is an element in the inverse image of $\alpha^*x$. 
\end{lemma}

\begin{proof}
The morphism $\iota_V^{T\oplus_{\id_T,\alpha}V}\ :\ V\rightarrow T\oplus_{\id_T,\alpha} V$ is an isomorphism. Indeed, the following commutative diagram is a pushout in $\E$:
$$\xymatrix{T\ar[r]^\alpha\ar[d]_{\id_T} & V\ar[d]^{\id_V}\\
T\ar[r]^\alpha & V.}$$
Then, $((\iota_V^{T\oplus_{\id_T,\alpha} V})^{-1})^*x$ satisfies the two following conditions, $\iota_V^{T\oplus_{\id_T,\alpha} V\ *}((\iota_V^{T\oplus_{\id_T,\alpha} V})^{-1})^*x=x$ and $\iota_T^{T\oplus_{\id_T,\alpha} V\ *}((\iota_V^{T\oplus_{\id_T,\alpha} V})^{-1})^*x=\alpha^*x$, where the second is a consequence of the following commutative diagram: $$\xymatrix{T\ar[r]^-\alpha\ar[d]_-{\id_T} & V\ar[d]^-{\iota_V}\ar@/^2.0pc/[rdd]^-{\id_V}&\\
T\ar[r]^-{\iota_T}\ar@/_2.0pc/[rrd]_-\alpha & T\oplus_{\id_T,\alpha}V\ar[rd]^{\iota_V^{-1}} &\\
& & V.}$$
Therefore it is in the inverse image of $\alpha^*x$ under $\rho_{F,(V,x)}$, if $(V,x)$ is central it is the only element in the inverse image.

\end{proof}

In the following, for $V,$ $W$ and $T$ finite dimensional vector spaces, we think of spaces of the form $(V\oplus_{\nu,\omega} W)\oplus_{g,\tau}T$ as $V+W+T$, where we identified $V$, $W$ and $T$ with their images in $(V\oplus_{\nu,\omega} W)\oplus_{g,\tau}T$, in order to use the associativity $(V+W)+T=V+(W+T)$. The following technical construction aims to identify a canonical isomorphism $\zeta_{\nu,\omega,g,\tau}$ from $(V\oplus_{\nu,\omega} W)\oplus_{g,\tau} T$ to the appropriate $V\oplus_{\nu',g'}(W\oplus_{\omega',\tau'} T)$.

\begin{notation}
For $V$, $W$ and $T$ in $\E$, we denote $$T_1(V,W;T):=\bigsqcup\limits_{\left[M,\nu,\omega\right]\in\textbf{B}(V,W)}\textbf{B}(V\oplus_{\nu,\omega}W,T),$$ and
$$T_2(V;W,T):=\bigsqcup\limits_{\left[M',\omega',\tau'\right]\in\textbf{B}(W,T)}\textbf{B}(V,W\oplus_{\omega',\tau'}T).$$
\end{notation}

\begin{lemma}\label{techtechtech}

\begin{enumerate}
    \item $T_1(\_,\_;\_)$ and $T_2(\_;\_,\_)$ are trifunctors on $\mathcal{VI}^{op}$.
    \item There is a natural transformation $\zeta$ from $T_1(\_,\_;\_)$ and $T_2(\_;\_,\_)$,
    \item for $\left[M,g,\tau\right]\in\textbf{B}(V\oplus_{\nu,\omega}W,T)\subset T_1(V,W;T)$ and for $(\nu',g',\omega',\tau')$ such that $\zeta(\left[M,g,\tau\right])=\left[M',\nu',g'\right]\in\textbf{B}(V,W\oplus_{\omega',\tau'}T)$, there is a natural isomorphism $$\zeta_{\nu,\omega,g,\tau}\ :\ (V\oplus_{\nu,\omega}W)\oplus_{g,\tau}T\rightarrow V\oplus_{\nu',g'}(W\oplus_{\omega',\tau'}T),$$ such that the canonical injections from $V$, $W$ and $T$ to $V\oplus_{\nu',g'}(W\oplus_{\omega',\tau'}T)$ factorizes through $\zeta_{\nu,\omega,g,\tau}$ and their injections in $(V\oplus_{\nu,\omega}W)\oplus_{g,\tau}T$.
\end{enumerate}

\end{lemma}

\begin{proof}
The fact that $T_1(\_,\_;\_)$ and $T_2(\_;\_,\_)$ are trifunctors is a direct consequence of the fact that $\textbf{B}(\_,\_)$ is a bifunctor. Let us construct the natural transformation $\zeta$.\\

For $\left[M,g,\tau\right]\in\textbf{B}(V\oplus_{\nu,\omega}W,T)\subset T_1(V,W;T)$, we consider $\iota_W+\iota_T$ from $W\oplus T$ to\newline $(V\oplus_{\nu,\omega}W)\oplus_{g,\tau} T$. Since $\iota_W$ and $\iota_T$ are injections from $W$ and $T$ to $(V\oplus_{\nu,\omega}W)\oplus_{g,\tau} T$,\newline $\ker(\iota_W+\iota_T)\cap W=\ker(\iota_W+\iota_T)\cap T=\{0\}$. Then, for all $x\in\ker(\iota_W\oplus\iota_T)$, $x$ has a non trivial component both in $W$ and $T$, hence there exists $\omega'(x)\in W\backslash\{0\}$ and $\tau'(x)\in T\backslash\{0\}$ such that $x=\omega'(x)-\tau'(x)$. Then, by the first isomorphism theorem, $\iota_W+\iota_T$ factorises through an injection $\iota_{W\oplus_{\omega',\tau'}T}$ from $W\oplus_{\omega',\tau'}T$ to $(V\oplus_{\nu,\omega} W)\oplus_{g,\tau}T$.\\

We do the same construction a second time. We consider $\iota_V+\iota_{W\oplus_{\omega',\tau'}T}$ from $V\oplus(W\oplus_{\omega',\tau'}T)$ to $(V\oplus_{\nu,\omega}W)\oplus_{g,\tau}T$. Since the two are injections, there are injections $\nu'$ and $g'$ from $\ker(\iota_V+\iota_{W\oplus_{\omega',\tau'}T})$ to $V$ and $W\oplus_{\omega',\tau'}T$ such that, for all $x\in\ker(\iota_V+\iota_{W\oplus_{\omega',\tau'}T})$, $x=\nu'(x)-g'(x)$. Then, by the first isomorphism theorem, $\iota_V+\iota_{W\oplus_{\omega',\tau'}T}$ factorises through an injection from $V\oplus_{\nu',g'}(W\oplus_{\omega',\tau'}T)$ to $(V\oplus_{\nu,\omega} W)\oplus_{g,\tau}T$. But, by construction, $\iota_V+\iota_{W\oplus_{\omega',\tau'}T}$ is surjective, therefore this injection is an isomorphism that we denote $\zeta_{\nu,\omega,g,\tau}$. We define $\zeta(\left[M,g,\tau\right]):=\left[\ker(\iota_V+\iota_{W\oplus_{\omega',\tau'}T}),\nu',g'\right]$.

By construction, $\zeta_{\nu,\omega,g,\tau}$ satisfies the desired factorizations.
\end{proof}

\begin{lemma}\label{centrinj2}
Let $F$ be an object in $\SVI$, $V$ be an object of $\mathcal{VI}$ and $x$ be an element of $F(V)$ such that $(V,x)$ is central. For $\alpha\ :\ T\rightarrow V$ in $\mathcal{VI}$, $(T,\alpha^*x)$ is central.
\end{lemma}

\begin{proof}
By Lemma \ref{altedef}, we have to prove that for every $y\in F(W)$ there is a unique $b\in F(W\oplus_{\omega,\tau} T)$ for some $\left[M,\omega,\tau\right]\in\textbf{B}(W,T)$, such that $\iota_W^*b=y$ and $\iota_T^*b=\alpha^*x$.\\

We will prove the existence first. For $a\in F(V\oplus_{\nu,\omega}W)$ satisfying $\iota_V^*a=x$, we have that $$\sigma_\alpha F(W)(a)=(\alpha\oplus_{\nu,\omega}\id_W)^*a$$ is in some $F(T\oplus_{\tau,\epsilon}W)$ and satisfies $\iota_T^*(\alpha\oplus_{\nu,\omega}\id_W)^*a=\alpha^*\iota_V^*a=\alpha^*x$, so $\sigma_\alpha F(W)(a)\in\Sigma_{(T,\alpha^*x)}F(W)$, then if $a$ is the unique element such that $\iota_W^*a=y$ (which exists since $(V,x)$ is central), then $b:=\sigma_\alpha F(W)(a)$ satisfies $\iota_T^*b=\alpha^*x$ and $\iota_W^*b=\iota_W^*a=y$.\\

Let us now prove the uniqueness. Let $b$ be in the preimage of $y$ in $\sigma_{F,(T,\alpha^*x)}F(W)$, with $$b\in F(W\oplus_{\omega,\tau} T).$$
Then, since $(V,x)$ is central, there is a unique $\left[M,g,\nu\right]\in\textbf{B}(W\oplus_{\omega,\tau} T,V)$ and a $$c\in F((W\oplus_{\omega,\tau} T)\oplus_{g,\nu} V)$$ such that $\iota_V^*c=x$ and $\iota_{W\oplus_{\omega,\tau} T}^*c=b$. We consider the isomorphism $\zeta_{\omega,\tau,g,\nu}$ given by Lemma \ref{techtechtech}, which take values in some space $W\oplus_{\omega',g'}(T\oplus_{\tau',\nu'}V)$.  Then, by construction, $\iota_{T\oplus_{\tau',\nu'} V}^*\zeta_{\omega,\tau,g,\nu}^*c$ satisfies $\iota_T^*\iota_{T\oplus_{\tau',\nu'} V}^*\zeta_{\omega,\tau,g,\nu}^*c=\alpha^*x$ and $\iota_V^*\iota_{T\oplus_{\tau',\nu'} V}^*\zeta_{\omega,\tau,g,\nu}^*c=x$. Therefore, since $(V,x)$ is central, Lemma \ref{injcentr} implies that $\iota_V\ :\ V\rightarrow T\oplus_{\tau',\nu'} V$ is an isomorphism and $(\iota_V^{-1})^*\iota_{T\oplus_{\tau',\nu'} V}^*\zeta_{\omega,\tau,g,\nu}^*c=x$. Then, identifying $W\oplus_{\omega',g'}(T\oplus_{\tau',\nu'}V)$ with the appropriate $W\oplus_{\omega'',\nu''}V$, we get that $\zeta_{\omega,\tau,g,\nu}^*c\in F(W\oplus_{\omega'',\nu''}V)$ satisfies $\iota_W^*\zeta_{\omega,\tau,g,\nu}^*c=y$ and $\iota_V^*\zeta_{\omega,\tau,g,\nu}^*c=x$. $\zeta_{\omega,\tau,g,\nu}^*c$ is then the unique element in the inverse image of $y$ by $\rho_{F,(V,x)}$ and $b=\sigma_{\alpha}F(W)(\zeta_{\omega,\tau,g,\nu}^*c)$, which proves the uniqueness of $b$.

\end{proof}

\begin{lemma}\label{somme}
Let $F$ be in $\SVI$ and $(V,x)$ and $(T,y)$ be central elements of $F$, then there exists, up to isomorphism, a unique pair $(R,z)$, with $R=V\oplus_{\nu,\tau}T$ for some $\left[M,\nu,\tau\right]\in\textbf{B}(V,T)$ and with $z\in F(R)$, such that $x=\iota_V^*z$ and $y=\iota_T^*z$. Moreover, $(R,z)$ is central.
\end{lemma}

\begin{proof}
By centrality of $(V,x)$ and by Lemma \ref{altedef}, $z$ is necessarily the only element in the inverse image of $y$ under $\rho_{F,(V,x)}$, which prove the existence and the uniqueness of $(R,z)$.\\

Let us prove that it is central. Let $e$ be in $F(W)$, we take $e'$ in $F(W\oplus_{\omega,\tau} T)$ the only element in the inverse image of $e$ under $\rho_{F,(T,y)}$ and $e''\in F((W\oplus_{\omega,\tau} T)\oplus_{g,\nu} V)$ the only element in the inverse image of $e'$ under $\rho_{F,(V,x)}$. We consider $\zeta_{\omega,\tau,g,\nu}$ as in lemma \ref{techtechtech}, an isomorphism from $(W\oplus_{\omega,\tau} T)\oplus_{g,\nu} V$ to some $W\oplus_{\omega',g'}(T\oplus_{\tau',\nu'} V)$. Then, $\iota_{T\oplus_{\tau',\nu'} V}^*(\zeta_{\omega,\tau,g,\nu}^{-1})^*e''$ is in the inverse image of $y$ under $\rho_{F,(V,x)}$, since $(V,x)$ is central $\iota_{T\oplus_{\tau',\nu'} V}^*(\zeta_{\omega,\tau,g,\nu}^{-1})^*e''=z$. Moreover, $\iota_W^*(\zeta_{\omega,\tau,g,\nu}^{-1})^*e''=e$, this imply that $(\zeta_{\omega,\tau,g,\nu}^{-1})^*e''$ is in the inverse image of $e$ under $\rho_{F,(R,z)}$.\\

Let us now show the uniqueness of the element in the inverse image of $e$. We take $R$ to be of the form $V\oplus_{\nu,\tau} T$, and we take $a\in F(W\oplus_{\omega,g} (V\oplus_{\nu,\tau}T))$ and $b\in F(W\oplus_{\epsilon,\gamma} (V\oplus_{\nu,\tau}T))$ two elements in the inverse image of $e$ under $\rho_{F,(R,z)}$. We consider the isomorphisms $\zeta_{\omega',\nu',g',\tau'}$ and $\zeta_{\epsilon',\nu'',\gamma',\tau''}$ given by Lemma \ref{techtechtech}, from some vector spaces $(W\oplus_{\omega',\nu'}V')\oplus_{g',\tau'}T$ and $(W\oplus_{\epsilon',\nu''}V')\oplus_{\gamma',\tau''}T$ to $W\oplus_{\omega,g} (V\oplus_{\nu,\tau}T)$ and $W\oplus_{\epsilon,\gamma} (V\oplus_{\nu,\tau}T)$. Set $a'=\iota_{W\oplus_{\omega',\nu'}V}^*\zeta_{\omega',\nu',g',\tau'}^*a\in F(W\oplus_{\omega',\nu'}V)$ and $b'=\iota_{W\oplus_{\epsilon',\nu''}V}^*\zeta_{\epsilon',\nu'',\gamma',\tau''}^*b\in F(W\oplus_{\epsilon',\nu''}V)$. They both satisfy that their image by $\iota_W^*$ is $e$ and by $\iota_V^*$ is $x$. Then, by the centrality of $(V,x)$ we have that $a'=b'$. Moreover, applying $\iota_T^*$ to $\zeta_{\omega',\nu',g',\tau'}^*a$ and  $\zeta_{\epsilon',\nu'',\gamma',\tau''}^*b$, we obtain $y$. Therefore, $\zeta_{\omega',\nu',g',\tau'}^*a$ and $\zeta_{\epsilon',\nu'',\gamma',\tau''}^*b$ are both in the inverse image of $a'=b'$ under $\rho_{F,(T,y)}$. Since $(T,y)$ is also central, $\zeta_{\omega',\nu',g',\tau'}^*a=\zeta_{\epsilon',\nu'',\gamma',\tau''}^*b$, hence $a=b$.
\end{proof}

\begin{theorem}
Let $F$ be an object in $\SVI$ such that $\textbf{C}(F)$ is not empty and such that $F(V)=\emptyset$ for $\dim(V)$ greater than some integer $d$, then there is a unique central element $(C,c)$ up to isomorphism satisfying that $(V,x)$ is central if and only if there is an injective morphism $\alpha$ from $V$ to $C$ such that $x=\alpha^*c$. We call $(C,c)$ the centre of $F$.
\end{theorem}

\begin{proof}
It is a direct consequence of lemma \ref{centrinj2} and \ref{somme}. Indeed, for $(C,c)$ a central element with $\dim(C)$ maximal, if $(V,x)$ is a central element such that there is no injection from $V$ to $C$ satisfying $\alpha^*c=x$, then the dimension of $(R,z)$ the unique pair satisfying the assumptions of lemma \ref{somme} for $(T,y)=(C,c)$ is greater than $\dim(C)$ which is absurd, since $\dim(C)$ is supposed to be maximal among central elements.
\end{proof}

\begin{corollary}
For $K$ a noetherian unstable algebra, if $K$ is $nil$-closed, the centre of $K$ is equal to the centre of $\HomKf(K,H^*(\_))$.
\end{corollary}

\section{Definition of the groupoid $\mathcal{G}_F$ and application to the computation of the centre of $F$}

For $F\in\SVI$, Yoneda's Lemma implies the existence of surjections from functors of the form $\bigsqcup\limits_{i\in I}\HomC(\_,W_i)$ to $F$. This statement is not specific to $\SVI$, we could make a similar one for $\SVF$. The interest we have about such surjections in $\SVI$ is that the simplicity of the category $\mathcal{VI}$ makes it easy to classify functors in $\SVI$ with a given surjections from $\bigsqcup\limits_{i\in I}\HomC(\_,W_i)$.
We start this section by explaining how to associate to objects $F$ in $\SVI$ with a given surjection from some $\bigsqcup\limits_{i\in I}\HomC(\_,W_i)$ to $F$, a groupoid $\mathcal{G}_F$ with objects the sub-vector spaces of the $W_i$. We will prove that the isomorphism class of $F$ in the coslice category with respect to $\bigsqcup\limits_{i\in I}\HomC(\_,W_i)$ is determined by $\mathcal{G}_F$. \\

In the second sub-section, we will prove that the centre of $F$ is explicitly determined by the groupoid $\mathcal{G}_F$ of the first sub-section.

\subsection{Definition of the groupoid}

By Yoneda's lemma, $F(W)$ is naturally isomorphic to $\Hom_{\SVI}(\HomC(\_,W),F)$ under the morphism which sends a natural transformation $\eta$ to $\eta(\id_V)$. Then, we can always exhibit a surjection from some $\bigsqcup\limits_{i\in I}\HomC(\_,W_i)$ to $F$, where the cardinal of $I$ may be very big. The most obvious one being defined by \newline
$I=\{x\in F(\Fp^d)\ ;\ d\in\N\}$, $W_x=\Fp^d$ for $x\in F(\Fp^d)$, and $q$ is the only natural transformation which sends $\id_{W_x}$ to $x$ for all $x\in I$. 

So $F$ is always isomorphic to some functor of the form $\bigsqcup\limits_{i\in I}\HomC(\_,W_i)/\sim_F$ with $\sim_F$ an equivalence relation on $\bigsqcup\limits_{i\in I}\HomC(\_,W_i).$ In this part, we will show how the equivalence relation $\sim_F$ associated with a surjection $$\bigsqcup\limits_{i\in I}\HomC(\_,W_i)\twoheadrightarrow F,$$ is encoded by a groupoid whose objects are the sub-spaces of the $W_i$. 

\begin{definition}
Let $F$ be an object in $\SVI$, a generating family of $F$ is a family $(W_i)_{i\in I}$ of objects of $\mathcal{VI}$, together with a surjection $q\ :\ \bigsqcup\limits_{i\in I}\HomC(\_,W_i)\twoheadrightarrow F$. 
\end{definition}

Even though the results of this sub-section do not require any finiteness condition on the cardinal of $I$, for $((W_i)_{i\in I},q)$ a generating family, the constructions might not be exploitable in the case where $I$ is not finite.

\begin{definition}
An object $F$ in $\SVI$ is finitely generated if $F$ admits a generating family $((W_i)_{i\in I},q)$ with $|I|$ finite.
\end{definition}

Let $F$ be an object in $\SVI$, and let $((W_i)_{i\in I},q)$ be a generating family for $F$. We denote by $F_i$ the image of $\HomC(\_,W_i)$ under $q$. Then, we have $F=\bigcup\limits_{i\in I}F_i$. We also denote by $\phi_i^F$ (or just $\phi_i$ when there is no ambiguity) the image of $\id_{W_i}$ under $q$.\\

Let us now define the main ingredient of this section.

\begin{definition}
For $F\in\SVI$ and $((W_i)_{i\in I},q)$ a generating family of $F$, let $\mathcal{G}_F^{((W_i)_{i\in I},q)}$ be the groupoid whose set of objects is the disjoint union of the sets of sub-spaces of the $W_i$ and whose morphisms from $U\subset W_i$ to $U'\subset W_j$ are the isomorphisms $\alpha$ from $U$ to $U'$ such that $\alpha^*\iota_{U'}^*\phi_j=\iota_{U}^*\phi_i$. 
\end{definition}

The groupoid $\mathcal{G}_F^{((W_i)_{i\in I},q)}$ depends heavily on the choice of a generating family $((W_i)_{i\in I},q)$, so cannot be made functorial on the category $\SVI$. It can nonetheless be made functorial on the category whose objects are provided with a surjection from a fixed $\bigsqcup\limits_{i\in I}\HomC(\_,W_i)$.

\begin{definition}\label{houlhoop}
Let $_{(W_i)_{i\in I}}\SVI$ be the slice category whose objects are pairs $(F,q_F)$ with $F\in\SVI$ and $q_F$ a natural surjection from $\bigsqcup\limits_{i\in I}\HomC(\_,W_i)$ to $F$, and whose morphisms from $(F,q_F)$ to $(G,q_G)$ are natural transformations $\eta$ from $F$ to $G$ such that the following diagram commutes:
$$\xymatrix{ & F\ar[dd]^-{\eta}\\
\bigsqcup\limits_{i\in I}\HomC(\_,W_i)\ar[ru]^-{q_F}\ar[rd]_-{q_G} &\\
& G.}$$
\end{definition}

\begin{proposition}
$(F,q_F)\mapsto \mathcal{G}_{(F,q_F)}:=\mathcal{G}_F^{(((W_i)_{i\in I},q_F)}$ is a functor from the category $_{(W_i)_{i\in I}}\SVI$ to the category of groupoids.
\end{proposition}

\begin{proof}
For $\eta$ a morphism from $(F,q_F)$ to $(G,q_G)$ two objects in $_{(W_i)_{i\in I}}\SVI$, we have to define a morphism $\mathcal{G}_\eta$ from $\mathcal{G}_{(F,q_F)}$ to $\mathcal{G}_{(G,q_G)}$. $\mathcal{G}_{(F,q_F)}$ has the same objects as $\mathcal{G}_{(G,q_G)}$, furthermore for all $i\in I$, $\eta(\phi_i^F)=\phi_i^G$, otherwise the diagram of Definition \ref{houlhoop} would not commute. Then, for all $\alpha\in\mathcal{G}_{(F,q_F)}(U\subset W_i,U'\subset W_j)$, the commutativity of the preceding diagram implies that $$\alpha^*\iota_{U'}^*\phi_j^G=\eta(\alpha^*\iota_{U'}\phi_j^F)=\iota_U^*\phi_i^G,$$ therefore $\alpha\in\mathcal{G}_{(G,q_G)}$. Hence, $\mathcal{G}_{(F,q_F)}$ is a sub-groupoid of $\mathcal{G}_{(G,q_G)}$ and we define $\mathcal{G}_\eta$ to be the inclusion of this sub-groupoid.
\end{proof}

In the following, we fix a family $(W_i)_{i\in I}$. When there can be no ambiguity, for $(F,q_F)\in _{(W_i)_{i\in I}}\SVI$, we will use the notation $\mathcal{G}_F$ instead of $\mathcal{G}_{(F,q_F)}$. We want to prove that the isomorphism class of $(F,q_F)$ in $_{(W_i)_{i\in I}}\SVI$ is determined by the groupoid $\mathcal{G}_{(F,q_F)}$.

\begin{lemma}\label{condidi}
Let $(F,q_F)$ be an object in $_{(W_i)_{i\in I}}\SVI$. Then, $\mathcal{G}_F$ satisfies the following property. For $\alpha\in\mathcal{G}_F(U\subset W_i,U'\subset W_j)$, for $M$ a sub-space of $U$, and for $\alpha_M\ :\ M\rightarrow \alpha(M)$ the restriction of $\alpha$ to $M$ corestricted to $\alpha(M)$,  $\alpha_M\in\mathcal{G}_F(M\subset W_i,\alpha(M)\subset W_j)$.
\end{lemma}

\begin{proof}
For $\iota_M^U$ the inclusion of $M$ in $U$ (we reserve the notation $\iota_M$ for the inclusion of $M$ in $W_i$) we have, $\alpha\circ\iota_M^U=\alpha_M$. Then, $\alpha^*\iota_{U'}^*\phi_j=\iota_U^*\phi_i$ implies that $\alpha_M^*\iota_{\alpha(M)}^*\phi_j=(\iota_M^U)^*\alpha^*\iota_{U'}^*\phi_j=(\iota_M^U)^*\iota_U^*\phi_i=\iota_M^*\phi_i$. Then, $\alpha_M\in\mathcal{G}_F(M\subset W_i,\alpha(M)\subset W_j)$.
\end{proof}

\begin{example}\label{exempleparad}
Let $W$ be a finite dimensional vector-space and $G$ be a sub-group of $\text{Gl}(W)$. We consider $F(V):=\HomC(V,W)/G$ where $G$ acts on $\HomC(V,W)$ by composition. $F$ is an object in $\SVI$. Then, for $q$ the canonical projection $\HomC(V,W)\twoheadrightarrow\HomC(V,W)/G$, $(F,q)\in\ _W\SVI$. $\mathcal{G}_F$ is the groupoid whose objects are the sub-vector spaces of $W$, and whose morphisms from $U$ to $U'$ are the isomorphisms $\alpha$ from $U$ to $U'$ such that $\alpha^*q(\iota_{U'})=q(\iota_U)$. But, $\alpha^*q(\iota_{U'})=q(\iota_U)$ implies that $q(\iota_{U'}\circ\alpha)=q(\iota_U)$, by definition of $q$, it implies that there is $g\in G$ such that $\iota_{U'}\circ\alpha=g\circ\iota_U$, in other terms $\alpha=g_U$. So $\mathcal{G}_F(U,U')$ is the set of isomorphisms of the form $g_U$, for $g\in G$ such that $g(U)=U'$.
\end{example}

\begin{definition}\label{restrictionprop}
For $\mathcal{G}$ a groupoid whose objects are the sub-vector spaces of the $W_i$ with $i\in I$, and whose morphisms are isomorphisms of vector spaces, we say that $\mathcal{G}$ has the restriction property if, for all $U\subset W_i$, for all $U'\subset W_j$ and for all $\alpha\in\mathcal{G}(U\subset W_i,U'\subset W_j)$, $\alpha_M$ is in $\mathcal{G}(M\subset W_i,\alpha(M)\subset W_j)$.
\end{definition}

\begin{remark}
Thus, Lemma \ref{condidi} asserts that, for $F\in\ _{(W_i)_{i\in I}}\SVI$, $\mathcal{G}_F$ has the restriction property.
\end{remark}

%\begin{hyp}\label{hyphyp}
%In the following, if not specified otherwise, $\mathcal{G}$ denotes a groupoid whose objects are the sub-spaces of the $W_i$ and whose morphisms are isomorphisms of vector spaces such that $\mathcal{G}$ has the restriction property.
%\end{hyp}

\begin{remark}
The set of groupoids whose set of objects is the disjoint union of the sets of sub-vector spaces of the $W_i$ and that satisfies the restriction property, is ordered by inclusion, where $\mathcal{G}\subset\mathcal{H}$ if, for each pair $(U\subset W_i,U'\subset W_j)$, $\mathcal{G}(U\subset W_i,U'\subset W_j)\subset\mathcal{H}(U\subset W_i,U'\subset W_j)$.  
\end{remark}

\begin{notation}\label{frakfrak}
We denote by $\textit{Groupoid}((W_i)_{i\in I})$ the poset of groupoids whose objects are the sub-spaces of the $W_i$ and that satisfies the restriction property.

For $W$ a finite dimensional vector space, let us also denote by $\textit{Group}(W)$ the poset of sub-groups of $\text{Gl}(W)$.
\end{notation}

\begin{definition}\label{profrak}
Let $\mathfrak{g}$ be the poset preserving map from $\textit{Group}(W)$ to $\textit{Groupoid}(W)$ which sends $G$ a subgroup of $\text{Gl}(W)$ to the groupoid $\mathfrak{g}(G)$ such that, for $U$ and $U'$ subspaces of $W$ and $\alpha$ from $U$ to $U'$, $\alpha\in\mathfrak{g}(G)(U,U')$ if and only if there is $g\in G$ sucht that $\alpha=g_U$.
\end{definition}

\begin{remark}
We have seen in Example \ref{exempleparad} that $\mathcal{G}_{(\HomC(\_,W)/G,q_G)}=\mathfrak{g}(G)$, for $q_G$ the canonical projection from $\HomC(\_,W)$ to $\HomC(\_,W)/G$.
\end{remark}

In the following, we want to show that for $\mathcal{G}$ a groupoid whose objects are the sub-vector spaces of the $W_i$ and whose morphisms are isomorphisms of vector spaces, there exists $F\in\ _{(W_i)_{i\in I}}\SVI$ such that $\mathcal{G}_F=\mathcal{G}$ if and only if $\mathcal{G}$ has the restriction property. We will then show that those groupoids are in on-to-one correspondence with isomorphism classes of objects in $_{(W_i)_{i\in I}}\SVI$.

\begin{definition}
For $\mathcal{G}\in\textit{Groupoid}((W_i)_{i\in I})$, let $\sim_\mathcal{G}$ be the relation on \newline$\bigsqcup\limits_{i\in I}\HomC(V,W_i)$ defined by $\zeta\sim_\mathcal{G}\epsilon$ if there is $\alpha\in\mathcal{G}(\Ima(\zeta)\subset W_i,\Ima(\epsilon)\subset W_j)$ such that $\Tilde{\epsilon}=\alpha\circ\Tilde{\zeta}$, where $\Tilde{\zeta}$ and $\Tilde{\epsilon}$ denote the corestriction of $\zeta$ and $\epsilon$ to their images.
\end{definition}

\begin{lemma}
For $\mathcal{G}\in\textit{Groupoid}((W_i)_{i\in I})$, $\sim_\mathcal{G}$ is an equivalence relation.

We will denote by $[\epsilon]_\mathcal{G}$ the equivalence class of $\epsilon$ for this equivalence relation.
\end{lemma}

\begin{proof}
Since for all $U\subset W_i$, $\id_{U}\in\mathcal{G}(U\subset W_i,U\subset W_i)$, the relation is reflexive. The transitivity comes from the composition of morphisms in $\mathcal{G}$ and the symmetry from the fact that $\mathcal{G}$ is a groupoid and hence that every morphism in $\mathcal{G}$ has an inverse.
\end{proof}

\begin{proposition}
For $\mathcal{G}\in\textit{Groupoid}((W_i)_{i\in I})$, $\bigsqcup\limits_{i\in I}\HomC(\_,W_i)/\sim_\mathcal{G}$ defines an object of $\SVI$.\\

Moreover, $\bigsqcup\limits_{i\in I}\HomC(\_,W_i)/\sim_\mathcal{G}$ is connected if and only if, for all $(i,j)\in I^2$,\newline 
$\mathcal{G}(0\subset W_i,0\subset W_j)$ contains the only morphism from $0$ to $0$.

\end{proposition}

\begin{proof}

To prove that  $\bigsqcup\limits_{i\in I}\HomC(\_,W_i)/\sim_\mathcal{G}$ is an object in $\SVI$, let us prove that if $\rho\sim_\mathcal{G}\zeta,$ for $\rho\in\HomC(V,W_i)$ and $\zeta\in\HomC(V,W_j)$ with $V\in\mathcal{VI}$, and if $\beta\ :\ U\rightarrow V$ is a morphism in $\mathcal{VI}$, then $\beta^*\rho\sim_\mathcal{G}\beta^*\zeta$. Let $\alpha$ be in $\mathcal{G}(\Ima(\rho)\subset W_i,\Ima(\zeta)\subset W_j)$ such that $\Tilde{\zeta}=\alpha\circ\Tilde{\rho}$. Then since $\mathcal{G}$ has the restriction property, $\alpha_{\rho\circ\beta(U)}\ :\ \rho\circ\beta(U)\rightarrow \zeta\circ\beta (U)$ is in $\mathcal{G}(\Ima(\rho\circ\beta)\subset W_i,\Ima(\zeta\circ\beta)\subset W_j)$ and $\widetilde{\zeta\circ\beta}=\alpha_{\rho\circ\beta(U)}\circ\widetilde{\rho\circ\beta}$ therefore $\beta^*\rho\sim_\mathcal{G}\beta^*\zeta$.
\end{proof}

\begin{theorem}\label{principal1}
\begin{enumerate}
    \item\label{oon} For $\mathcal{G}\in\textit{Groupoid}((W_i)_{i\in I})$, $$\mathcal{G}_{(\bigsqcup\limits_{i\in I}\HomC(\_,W_i)/\sim_\mathcal{G},q)}=\mathcal{G},$$ for $q$ the canonical surjection from $\bigsqcup\limits_{i\in I}\HomC(\_,W_i)$ to $\bigsqcup\limits_{i\in I}\HomC(\_,W_i)/\sim_\mathcal{G}$.
    \item Conversely, let $F\in\ _{(W_i)_{i\in I}}\SVI$. Then, $F$ is isomorphic to $\bigsqcup\limits_{i\in I}\HomC(\_,W_i)/\sim_{\mathcal{G}_F}$. 
\end{enumerate}

\end{theorem}

\begin{proof}

Let us prove the first point. Let $\alpha$ be an isomorphism from $U\subset W_i$ to $U'\subset W_j$. $\alpha\in \mathcal{G}_{\bigsqcup\limits_{i\in I}\HomC(\_,W_i)/\sim_\mathcal{G}}(U\subset W_i,U'\subset W_j)$ if and only if $\alpha^*\left[\iota_{U'}\right]_\mathcal{G}=\left[\iota_{U}\right]_\mathcal{G}$. This the case if and only if $\alpha\in\mathcal{G}(U\subset W_i,U'\subset W_j)$.\\

By construction, for $\alpha$ an isomorphism from $U$ to $U'$ with $U$ and $U'$ sub-spaces of $W_i$ and $W_j$, $\alpha^*\iota_{U'}^*[id_{W_j}]_{\mathcal{G}}=\iota_{U}^*[id_{W_i}]_{\mathcal{G}}$ if and only if $\alpha\in\mathcal{G}(U\subset W_i,U'\subset W_j)$.  
We only have to prove that $[\epsilon]_{\mathcal{G}_F}\mapsto \epsilon^*\phi_i$, for $\epsilon\in\HomC(V,W_i)$, is well defined and defines a bijection from $\bigsqcup\limits_{i\in I}\HomC(V,W_i)/\sim_{\mathcal{G}_F}$ to $F(V)$ for $V\in\mathcal{VI}$.\\ 

For the second point, if $\epsilon\sim_{\mathcal{G}_F}\zeta$, there exists $\alpha\ :\ \epsilon(V)\rightarrow \zeta(V)$ such that $\alpha\in\mathcal{G}_F(\Ima(\epsilon)\subset W_i,\Ima(\zeta)\subset W_j)$ and $\Tilde{\zeta}=\alpha\circ\Tilde{\epsilon}$. Then, 
\begin{align*}
     \zeta^*\phi_j & =\Tilde{\zeta}^*\iota_{\zeta(V)}^*\phi_j,\\
     & =\Tilde{\epsilon}^*\alpha^*\iota_{\zeta(V)}^*\phi_j,\\
     & =\Tilde{\epsilon}^*\iota_{\epsilon(V)}^*\phi_i,\\
     & =\epsilon^*\phi_i.
\end{align*}

The map which maps $[\epsilon]_{\mathcal{G}_F}$ to $\epsilon^*\phi_i$ is then well defined. It is obviously surjective.\\

Let us prove that it is injective. Let $\epsilon$ and $\zeta$ be morphisms from $V$ to $W_i$ and $W_j$ such that $\epsilon^*\phi_i=\zeta^*\phi_j$. Then, $\Tilde{\epsilon}^*\iota_{\epsilon(V)}^*\phi_i=\Tilde{\zeta}^*\iota_{\zeta(V)}^*\phi_j$. Since $\Tilde{\epsilon}$ is an isomorphism from $V$ to the image of $\epsilon$, we have $\iota_{\epsilon(V)}^*\phi_i=(\Tilde{\zeta}\circ\Tilde{\epsilon}^{-1})^*\iota_{\zeta(V)}^*\phi_j$. Therefore $\Tilde{\zeta}\circ\Tilde{\epsilon}^{-1}\in\mathcal{G}_F(\Ima(\epsilon)\subset W_i,\Ima(\zeta)\subset W_j)$ and $\epsilon\sim_{\mathcal{G}_F}\zeta$.
\end{proof}

The following is a ``converse" of example \ref{exempleparad}.

\begin{example}
Let $F$ be an object in $_W\SVI$ and consider the group $G=\mathcal{G}_F(W,W)$. We know, since $\mathcal{G}_F$ has the restriction property, that $\mathcal{G}_F\subset\mathfrak{g}(G)$ (see Definition \ref{profrak}). Suppose that this inclusion is an equality. Then, $\sim_\mathcal{G}$ is the equivalence relation defined by $\rho\sim_\mathcal{G}\zeta$ if and only if there is $g\in G$ such that $g\circ\rho=\zeta$. By Theorem \ref{principal1}, $F\cong\HomC(\_,W)/G$.
\end{example}

\subsection{Computation of $\textbf{C}(F)$ using $\mathcal{G}_F$}

In this sub-section, we take $F\in\ _{(W_i)_{i\in I}}\SVI$. We want to prove that, under some assumptions on the generating family of $F$, the central elements of $F$ are determined by the groupoid $\mathcal{G}_F$. \\

We start by proving that, if $F$ has a generating family with one element, $\rho_{F,(V,x)}$ is surjective for all $V\in\mathcal{VI}$ and for all $x\in F(V)$. 

\begin{notation}
For $F\in\SVI$ and for $(W,q)$ a generating family with one element, we denote by $\phi^F\in F(W)$ the image of $\id_W$ under $q$.
\end{notation}

On one hand, for $F\in\ _W\SVI,$ the functor $F$ is determined up to isomorphism by the groupoid $\mathcal{G}_F$. On the other hand, Lemma \ref{altedef} gives a criterion for the centrality of a pair $(V,x)$ through the sets $F(V\oplus_{\nu,\mu}U)$ with the $V\oplus_{\nu,\mu} U\in\textbf{B}(V,U)$ which are not sub-vector spaces of $W$. To characterise the centre of $F$ using $\mathcal{G}_F$, we need to reformulate the centrality condition by comparing the $V\oplus_{\nu,\mu}U$ with sub-vector spaces of $W$.

\begin{remark}\label{lastremark}
Consider $\delta\ :\ V\rightarrow W$ and $\epsilon\ :\ U\rightarrow W$ two morphisms in $\mathcal{VI}$. For $\delta+\epsilon$ the morphism from $V\oplus U$ to $W$ which sends $v$ to $\delta(v)+\epsilon(v)$, since $\delta$ and $\epsilon$ are injective, there exists $\nu$ and $\mu$ from $\ker(\delta+\epsilon)$ to $V$ and $U$ such that for all $v\in\ker(\delta+\epsilon)$, $v=\nu(v)-\mu(v)$. Then, by the universal property of the pushout, $(V\oplus U)/\ker(\delta+\epsilon)$ is isomorphic to $V\oplus_{\nu, \mu} U$. The first isomorphism theorem implies that $\delta+\epsilon$ induces an isomorphism between $V\oplus_{\nu, \mu} U$ and $\delta(V)+\epsilon(U)$ (which is the image of $\delta+\epsilon$).
\end{remark}

\begin{notation}
For $\delta\ :\ V\rightarrow W$ and $\epsilon\ :\ U\rightarrow W$ two morphisms in $\mathcal{VI}$, let $\nu$ and $\mu$ be as in Remark \ref{lastremark}. Then, let $\epsilon\uparrow\delta$ denote the isomorphism from $V\oplus_{\nu,\mu} U$ to $\delta(V)+\epsilon(U)$ induced by the first isomorphism theorem. It satisfies $\delta+\epsilon=\iota_{\delta(V)+\epsilon(U)}\circ (\epsilon\uparrow\delta)\circ\pi$, for $\pi$ the canonical projection from $V\oplus U$ to $V\oplus_{\nu,\mu}U$ and $\iota_{\delta(V)+\epsilon(U)}$ the inclusion from $\delta(V)+\epsilon(U)$ in $W$.
\end{notation}

\begin{lemma}\label{centritri}
For $F\in\ _W\SVI$, $\delta$ from $V$ to $W$ and $\epsilon$ from $U$ to $W$, $(\epsilon\uparrow\delta)^*\iota_{\delta(V)+\epsilon(U)}^*\phi^F\in F(V\oplus_{\nu,\mu} U)$ is in the inverse image of $\epsilon^*\phi^F$ under $\rho_{F,(V,\delta^*\phi^F)}$.
\end{lemma}

\begin{proof}
By definition of $\epsilon\uparrow\delta$, we have $\iota_{\delta(V)+\epsilon(U)}\circ(\epsilon\uparrow\delta)\circ\iota_V=\delta$ and $\iota_{\delta(V)+\epsilon(U)}\circ(\epsilon\uparrow\delta)\circ\iota_U=\epsilon$, for $\iota_V$ and $\iota_U$ the inclusions of $U$ and $V$ in $V\oplus_{\nu,\mu}U$. Therefore, $\iota_U^*((\epsilon\uparrow\delta)^*\iota_{\delta(V)+\epsilon(U)}^*\phi^F)=\epsilon^*\phi^F$ and $\iota_V^*((\epsilon\uparrow\delta)^*\iota_{\delta(V)+\epsilon(U)}^*\phi^F)=\delta^*\phi^F$.
\end{proof}

This lemma implies two things, firstly that $\rho_{F,(V,\delta^*\phi^F)}$ is always surjective, for $F\in\ _W\SVI$ and $\delta$ from $V$ to $W$. Secondly, that $\delta^*\phi^F\in F(V)$ is central if and only if, for all morphisms $\epsilon\ :\ U\rightarrow W$ in $\mathcal{VI}$, $(\epsilon\uparrow\delta)^*\iota_{\delta(V)+\epsilon (U)}^*\phi^F\in F(V\oplus_{\nu,\mu}U)$ is the only element in the inverse image of $\epsilon^*\phi^F$ under $\rho_{F,(V,\delta^*\phi^F)}$.\\ 

In the case $F\in\ _{(W_i)_{i\in I}}\SVI$ with $|I|$ greater than one, we recall that $F_i$ denotes the image of $\HomC(\_,W_i)$ in $F$. We want to use the surjectivity of each $\rho_{F_i,(V,x)}$ with $i\in I$ and $x\in F_i(V)$. We start by proving that, under some condition on the surjection $q_F$ from $\bigsqcup\limits_{i\in I}\HomC(\_,W_i)$ to $F$, if $(V,x)$ is central, then $x\in F_i(V)$ for all $i\in I$. 

\begin{definition}
Let $F$ be an object in $\SVI$, with generating family $((W_i)_{i\in I},q)$. We say that $((W_i)_{i\in I},q)$ is a minimal generating family, if for all $i\neq j$, $F_i\not\subset F_j$.
\end{definition}

\begin{remark}
We can always extract a minimal sub-generating family from a given generating family.
\end{remark}

\begin{lemma}
Let $F$ be in $\SVI$ and $((W_i)_{i\in I},q)$ a minimal generating family of $F$. Let also $(V,x)$ be a central element of $F$. Then, $x\in\bigcap\limits_{i\in I} F_i(V)$ and $(V,x)$ is a central element for each $F_i$.
\end{lemma}

\begin{proof}
For each $i$, let $x_i\in F(W_i\oplus_{\omega_i,\nu_i} V)$ be in the inverse image of $\phi_i$ under $\rho_{F,(V,x)}$. If $x_i\in F_j(W_i\oplus_{\omega_i,\nu_i} V)$, $\phi_i=\iota_{W_i}^*x_i$ hence $F_i$ is a sub-functor of $F_j$. Since we supposed the generating family to be minimal, this implies that $i=j$. Then, $x=\iota_{V}^*x_i$ and therefore $x\in F_i(V)$.\\
By Lemma \ref{centritri}, we get that an element in $F_i(E)$, for $E\in\mathcal{VI}$, has at least one element in its inverse image under $\rho_{F_i,(V,x)}$. If $(V,x)$ is central for $F$, this element is unique and $(V,x)$ is also central for $F_i$.
\end{proof}

Now, let us consider $x\in\bigcap\limits_{i\in I}F_i(V)$. We want to find a criterion on $\mathcal{G}_F$ to determine whether $(V,x)$ is central. Since, $x$ is in every $F_i(V)$ and that each $F_i$ are generated by one element, Lemma \ref{centritri} implies that every $\rho_{F_i,(V,x)}$ is surjective, and so $\rho_{F,(V,x)}$ is surjective. So the only way that $(V,x)$ can fail to be central is if a given element has two elements in its inverse image under $\rho_{F,(V,x)}$.\\

Since $x\in\bigcap\limits_{i\in I}F_i(V)$, for all $i\in I$, there is some inclusion $\delta_i$ from $V$ to $W_i$ such that $x=\delta_i^*\phi_i$ and we have the following lemma.

\begin{lemma}\label{nouvtech}
For $F\in _{(W_i)_{i\in I}}\SVI$ and $(V,x)$ a central element of $F$, we consider for all $i\in I$, $\delta_i$ from $V$ to $W_i$ such that $x=\delta_i^*\phi_i$. Then, for $\epsilon_k\ :\ U\rightarrow W_k$ and $\epsilon_j\ :\ U\rightarrow W_j$ such that $\epsilon_k^*\phi_k=\epsilon_j^*\phi_j$, we have $(\epsilon_k\uparrow{\delta_k})^*\iota_{\delta_k(V)+\epsilon_k(U)}^*\phi_k=(\epsilon_j\uparrow{\delta_j})^*\iota_{\delta_j(V)+\epsilon_j(U)}^*\phi_j$.
\end{lemma}

\begin{proof}
 Lemma \ref{centritri} states that $(\epsilon_k\uparrow{\delta_k})^*\iota_{\delta_k(V)+\epsilon_k(U)}^*\phi_k$ and $(\epsilon_j\uparrow{\delta_j})^*\iota_{\delta_j(V)+\epsilon_j(U)}^*\phi_j$ are in the inverse image of $\epsilon_k^*\phi_k=\epsilon_j^*\phi_j$, respectively under $\rho_{F_j,(V,x)}$ and $\rho_{F_k,(V,x)}$. Hence, each of them are in its inverse image under $\rho_{F,(V,x)}$. Thus, if $(V,x)$ is central, we necessarily have that $(\epsilon_k\uparrow{\delta_k})^*\iota_{\delta_k(V)+\epsilon_k(U)}^*\phi_k=(\epsilon_j\uparrow{\delta_j})^*\iota_{\delta_j(V)+\epsilon_j(U)}^*\phi_j$. 
\end{proof}
 
 This gives us a necessary condition on $\mathcal{G}_F$, for $(V,x)$ to be central. This condition will also be sufficient.\\

For $\alpha\in\mathcal{G}(U\subset W_i,U'\subset W_j)$, $\alpha^*\iota_{U'}^*\phi_j=\iota_{U}^*\phi_i$ has two, possibly equal, elements in its inverse image under $\rho_{F,(V,x)}$, where $x=\delta_i^*\phi_i$ for all $i$. Namely $((\iota_{U'}\circ\alpha)\uparrow{\delta_j})^*\iota_{U'+\delta_j(V)}^*\phi_j$ and $(\iota_U\uparrow{\delta_i})^*\iota_{U+\delta_i(V)}\phi_i$. If $(V,x)$ is central, those two have to be equal, in particular they have to live in the same $F(V\oplus_{\nu,\mu} U)$, which gives us a first condition on $\mathcal{G}_F$.

\begin{lemma}\label{centritru2}
Let $F$ in $\SVI$ and $((W_i)_{i\in I},q)$ a minimal generating family of $F$. Let $(V,x)$ be a central element of $F$ and for all $i\in I$, let $$\delta_i\ :\ V\rightarrow W_i,$$ such that $x=\delta_i^*\phi_i$. Then, for $\alpha\in\mathcal{G}_F(U\subset W_i,U'\subset W_j)$, if $v\in U\cap V_i$, $\alpha(v)=\delta_j\circ(\delta_i|^{V_i})^{-1}(v)$, for $V_i$ the image of $V$ under $\delta_i$.
\end{lemma}

\begin{proof}
We consider $$(\iota_U\uparrow{\delta_i})^*(\iota_{U+ V_i}^{W_i})^*\phi_i\in F_i(V\oplus_{\nu,\mu}U)$$ and $$((\iota_{U'}\circ\alpha)\uparrow{\delta_j})^*(\iota_{U'+ V_j}^{W_j})^*\phi_j\in F_j(V\oplus_{\nu',\mu'}U),$$ where $\mu$ is the injection of $U\cap V_i$ into $U$, $\nu$ is the restriction to $U\cap V_i$ of $(\delta_i|^{V_i})^{-1}$, and $\nu'$ and $\mu'$ are the restriction to $U'\cap V_j$ of $\alpha^{-1}$ and $(\delta_j|^{V_j})^{-1}$. They are equal by Lemma \ref{nouvtech}. In particular they need to live in the same set. We get that $V\oplus_{\nu,\mu} U=V\oplus_{\nu',\mu'} U$, which implies that, for $v\in U\cap V_i$, $\alpha(v)\in V_j$ and $(\delta_j|^{V_j})^{-1}(\alpha(v))=(\delta_i|^{V_i})^{-1}(v)$, which means that $\alpha(v)=\delta_j\circ(\delta_i|^{V_i})^{-1}(v)$.
\end{proof}

If the condition of Lemma \ref{centritru2} is satisfied, both $((\iota_{U'}\circ\alpha)\uparrow{\delta_j})^*\iota_{U'+\delta_j(V)}^*\phi_j$ and \newline$(\iota_U\uparrow{\delta_i})^*\iota_{U+\delta_i(V)}^*\phi_i$ live in $F(V\oplus_{\nu,\mu} U)$ where $\mu$ is the inclusion from $U\cap V_i$ in $U$ and $\nu$ is the restriction to $U\cap V_i$ of $(\delta_i|^{V_i})^{-1}$. 
Since by definition $\iota_U\uparrow{\delta_i}$ and $(\iota_{U'}\circ\alpha)\uparrow{\delta_j}$ are isomorphisms from $V\oplus_{\nu,\mu} U$ respectively to $V_i+U$ and $V_j+U'$, $((\iota_{U'}\circ\alpha)\uparrow{\delta_j})\circ(\iota_U\uparrow{\delta_i})^{-1}$ is an isomorphism from $U+V_i$ to $U'+V_j$ that we can compute.

\begin{lemma}\label{couterop}
Let $\alpha$ be an isomorphism from $U\subset W_i$ to $U'\subset W_j$ and let $\delta_i$ and $\delta_j$ be two injections from $V$ respectively to $W_i$ and $W_j$ whose images are $V_i$ and $V_j$. Then, if for every $v\in U\cap V_i$, $\alpha(v)=\delta_j\circ(\delta_i|^{V_i})^{-1}(v)$,  $\iota_U\uparrow{\delta_i}$ and $(\iota_{U'}\circ\alpha)\uparrow{\delta_j}$ are isomorphism from $V\oplus_{\nu,\mu} U$ respectively to $V_i+U$ and $V_j+U'$, where $\nu$ is the inclusion of $U\cap V_i$ in $U$ and $\mu$ is the restriction to $U\cap V_i$ of $(\delta_i|^{V_i})^{-1}$.\\

Furthermore, $((\iota_{U'}\circ\alpha)\uparrow{\delta_j})\circ(\iota_U\uparrow{\delta_i})^{-1}$ is the isomorphism from $U+V_i$ to $U'+V_j$ which sends $u\in U$ to $\alpha(u)$ and $v\in V_i$ to $\delta_j((\delta_i|^{V_i})^{-1}(v))$.
\end{lemma}

\begin{proof}
The first part of the statement is a direct consequence of the definition of $\iota_U\uparrow{\delta_i}$ and $(\iota_{U'}\circ\alpha)\uparrow{\delta_j}$. 

For $u\in U$, $(\iota_U\uparrow{\delta_i})^{-1}(u)=\iota_U^{V\oplus_{\nu,\mu} U}(u)$, then $((\iota_{U'}\circ\alpha)\uparrow{\delta_j})(\iota_U^{V\oplus_{\nu,\mu} U}(u))=\alpha(u)$, and for $v\in V_i$, $(\iota_U\uparrow{\delta_i})^{-1}(v)=\iota_V^{V\oplus_{\nu,\mu}U}((\delta_i|^{V_i})^{-1}(v))$ and $((\iota_{U'}\circ\alpha)\uparrow{\delta_j})(\iota_V^{V\oplus_{\nu,\mu}U}((\delta_i|^{V_i})^{-1}(v)))=\delta_j\circ(\delta_i|^{V_i})^{-1}(v)$. 
\end{proof}

\begin{notation}
Let $\alpha$ be an isomorphism from $U\subset W_i$ to $U'\subset W_j$ and let $\delta_i$ and $\delta_j$ be two injections from $V$ respectively to $W_i$ and $W_j$ whose images are $V_i$ and $V_j$ and such that for every $v\in U\cap V_i$, $\alpha(v)=\delta_j\circ(\delta_i|^{V_i})^{-1}(v)$. We denote by $\Bar{\alpha}$ the morphism from $U+V_i$ to $U'+V_j$ which sends $u\in U$ to $\alpha(u)$ and $v\in V_i$ to $\delta_j((\delta_i|^{V_i})^{-1}(v))$. 
\end{notation}

If $(V,x)$ is central, we have $((\iota_{U'}\circ\alpha)\uparrow{\delta_j})^*\iota_{U'+\delta_j(V)}^*\phi_j=(\iota_U\uparrow{\delta_i})^*\iota_{U+\delta_i(V)}^*\phi_i$, which implies that $\Bar{\alpha}$ must be in $\mathcal{G}_F(U+V_i\subset W_i,U'+V_j\subset W_j)$. We state the principal theorem of this sub-section.

\begin{theorem}\label{principal2}
Let $F$ be in $\SVI$ and $((W_i)_{i\in I},q)$ a minimal generating family of $F$. Let also $x\in F(V)$ and $(\delta_i\ :\ V\rightarrow W_i)_{i\in I}$ a family of injective morphisms from $V$ to the $W_i$, such that, for all $i\in I$, $x=\delta_i^*\phi_i$. Then, $(V,x)$ is central if and only if the two following conditions are satisfied: 
\begin{enumerate}
    \item\label{possible} for all $\alpha\in\mathcal{G}_F(U\subset W_i,U'\subset W_j)$, and for all $v\in U\cap V_i$, $\alpha(v)=\delta_j\circ(\delta_i|^{V_i})^{-1}(v)$,
    \item  for all sub-spaces $U$ of $W_i$ and $U'$ of $W_j$ and for all isomorphism $\alpha$ from $U$ to $U'$ satisfying $\alpha(v)=\delta_j\circ(\delta_i|^{V_i})^{-1}(v)$, for $v\in U\cap V_i$, $\alpha\in\mathcal{G}_F(U\subset W_i,U'\subset W_j)$ if and only if $\Bar{\alpha}\in\mathcal{G}_F(U+ V_i\subset W_i,U'+ V_j\subset W_j)$.
\end{enumerate}
\end{theorem}

\begin{proof}
Let us show first necessity. We suppose that $(V,x)$ is central, then by Lemma \ref{centritru2}, the first condition is satisfied. We consider $\alpha\in\mathcal{G}_F(U\subset W_i,U'\subset W_j)$, then by Lemma \ref{nouvtech}, $(\iota_U\uparrow\delta_i)^*(\iota_{U+ V_i}^{W_i})^*\phi_i=((\iota_{U'}\circ\alpha)\uparrow\delta_j)^*(\iota_{\alpha(U)+ V_j}^{W_j})^*\phi_j$. By Lemma \ref{couterop}, this implies that $\Bar{\alpha}^*\iota_{U'+V_j}^*\phi_j=\iota_{U+V_i}^*\phi_i$. Then, $\Bar{\alpha}\in\mathcal{G}_F(U+ V_i\subset W_i,U'+ V_j\subset W_j)$. Conversely, if $\Bar{\alpha}\in\mathcal{G}_F(U+ V_i\subset W_i,U'+ V_j\subset W_j)$, since $\mathcal{G}_F$ has the restriction property, $\alpha=\Bar{\alpha}_U\in\mathcal{G}_F(U\subset W_i,U'\subset W_j)$.\\

We now prove sufficiency. Suppose that $F$ and $(V,x)$ satisfy the two conditions. We consider $\zeta^*\phi_i\in F_i(V\oplus_{\nu,\mu} U)$ in the inverse image of an element $\beta^*\phi_j\in F(U)$, where $j$ is not necessarily equal to $i$. The first condition implies that, if $\eta\ :\ V\rightarrow W_i$ satisfies $\eta\sim_{\mathcal{G}_F}\delta_i$, then $\eta=\delta_i$. But $\zeta|_V^*\phi_i=\delta_i^*\phi_i$, hence $\zeta|_V=\delta_i$. Then, the equality $\zeta|_U^*\phi_i=\beta^*\phi_j$ implies the existence of a morphism $\alpha\in\mathcal{G}_F(\beta(U), \zeta(U))$ such that $\widetilde{\zeta|_U}=\alpha\circ\Tilde{\beta}$, where $\widetilde{\zeta|_U}$ and $\Tilde{\beta}$ are the corestrictions of $\zeta|_U$ and $\beta$ to their images. Therefore, $\zeta=(\iota_{\zeta(U)+ V_i}^{W_i})\circ\Bar{\alpha}\circ(\beta\uparrow{\delta_j})$ and by the second condition, we get that $\zeta^*\phi_i=(\beta\uparrow{\delta_j})^*(\iota_{\beta(U)+ V_j}^{W_j})^*\phi_j$. Hence, $(\beta\uparrow{\delta_j})^*(\iota_{\beta(U)+ V_j}^{W_j})^*\phi_j$ is the only element in the inverse image of $\beta^*\phi_j$ under $\rho_{F,(V,x)}$.
\end{proof}

\begin{example}
When $F\in\ _W\SVI$, the two conditions become simpler. We consider $\delta$ from $V$ to $W$, such that $(V,\delta^*\phi^F)$ is central.

The first condition becomes that, for $\alpha\in\mathcal{G}(U,U')$ and for $v\in U\cap \delta(V)$, $\alpha(v)=v$. Furthermore, in this case, $\Bar{\alpha}$ is the morphism from $U+\delta(V)$ to $U'+\delta(V)$ which sends $u\in U$ to $\alpha(u)$ and $v\in \delta(V)$ to $v$.
\end{example}

We return to Example \ref{exempleparad}.

\begin{proposition}
For $G$ a sub-group of $\text{Gl}(W)$, the centre of $F=\HomC(\_,W)/G$ is given by $(C,\iota_C)$ with $C$ the maximal sub-vector space of $W$ such that for all $v\in C$ and $g\in G$, $g(v)=v$ and with $\iota_C$ the inclusion of $C$ in $W$.
\end{proposition}

\begin{proof}
As we have seen in Example \ref{exempleparad}, $\mathcal{G}_F=\mathfrak{g}(G)$ (see Definition \ref{profrak}). Let $\delta$ be an injection from some vector space $V$ to $W$ such that $(V,q(\delta))$ is central, for $q$ the canonical projection from $\HomC(\_,W)$ to $F$. We will denote by $V'$ the sub-space $\delta(V)$ of $W$. The first condition of Theorem \ref{principal2} implies that, for all $g\in G$, $g\in\mathcal{G}_F(W,W)$, therefore for all $x\in V'$, $g(x)=x$. Conversely, if, for all $g\in G$ and for $x\in V'$, $g(x)=x$, since for all sub-spaces $U$ of $W$, and for all $\alpha\in\mathcal{G}_F(U,U')$, there exists $g\in G$ such that $\alpha=g_U$, $(V,q(\delta))$ satisfies the first condition of the theorem.\\

In this case, if $(V,q(\delta))$ satisfies the first condition of the theorem, the second condition of Theorem \ref{principal2} is immediately satisfied. Indeed, for $\alpha\in\mathcal{G}_F(U,U')$, $\Bar{\alpha}$ is the morphism from $U+V'$ to $U'+V'$ which sends $u\in U$ to $\alpha(u)$ and $v\in V'$ to $v$. Then, let $g\in G$ such that $\alpha=g_U$. Since for all $v\in V'$, $g(v)=v$, $g_{U+V'}=\Bar{\alpha}$. Hence, $\Bar{\alpha}\in\mathcal{G}_F(U+V',U'+V')$.\\

Therefore, $q(\delta)$ is central if and only if $\delta$ factorizes through $\iota_C$. Thus, $(C,\iota_C)$ is the centre of $F$.
\end{proof}

\section{The algebras $H^*(W)^\mathcal{G}$}

This section shows how to apply the groupoid $\mathcal{G}_F$ of the last section to classification problems of $nil$-closed, integral, noetherian, unstable algebras. Before we explain the focus of this section, let us recall the theorem of Adams-Wilkerson.

\begin{definition}\cite[Part II.2]{HLS2}
For $K\in\K$, the transcendence degree of $K$ is $d\in\N\cup\{\infty\}$, the supremum of the cardinals of finite sets of homogeneous elements in $K$ which are algebraically independent .
\end{definition}

\begin{remark}
If $K$ is noetherian, the transcendence degree of $K$ is finite.
\end{remark}

Let us recall the theorem of Adams-Wilkerson. 

\begin{theorem}\cite[Theorem 3]{HLS2}
Let $K$ be an integral, unstable algebra of transcendence degree lesser or equal to $\dim(W)$, then there exists an injection $\phi$ from $K$ to $H^*(W)$. Furthermore, this injection is regular if and only if the transcendence degree of $K$ equals $\dim(W)$.
\end{theorem}

Therefore, every integral, $nil$-closed, noetherian, unstable algebra is isomorphic to a $nil$-closed, noetherian sub unstable algebra of some $H^*(W)$. In the first sub-section we define $H^*(W)^\mathcal{G}$  for $\mathcal{G}$ a groupoid on sub spaces of $W$ that satisfies the restriction property. Then, $\mathcal{G}\mapsto H^*(W)^\mathcal{G}$ defines an explicit one-to-one correspondence between the groupoids on sub spaces of $W$ satisfying the restriction property and the noetherian, $nil$-closed, unstable sub algebras of $H^*(W)$ of transcendence degree $\dim(W)$.

For simplicity, we restrict to the integral case. The constructions generalizes to the non integral case. \\

Let us now recall the definition of the primitive elements of a comodule.
 
 \begin{definition}
 For $K\in\K$ provided with a $H^*(V)$-comodule structure $\kappa$ in $\K$, the algebra of primitive elements of $K$ is the sub-algebra of $K$ whose elements are those satisfying that $\kappa(x)=x\otimes 1$, for $1$ the unit of $H^*(V)$. We will denote by $P(K,\kappa)$ the algebra of primitive elements of $K$ for the $H^*(V)$-comodule structure $\kappa$.
 \end{definition}
 
 \begin{remark}
 By Corollary \ref{DWCT}, for all $(V,\phi)\in\textbf{C}(K)$, there is a unique structure $\kappa_\phi$ of $H^*(V)$-comodule on $K$ such that $(\epsilon_K\otimes\id_{H^*(V)})\circ\kappa_\phi=\phi$.
 \end{remark}
 
 \begin{notation}
 We will also denote $P(K,\kappa_\phi)$ by $P(K,\phi)$.
 \end{notation}
 
 The problem that we are interested in is the following. If we fix $V$ some finite dimensional vector space and $P$ some unstable algebra, can we classify, under suitable hypothesis, the connected, noetherian, $nil$-closed unstable algebras $K$, satisfying that $K$ admit a $H^*(V)$-comodule structure $\kappa$ in $\K$, whose algebra of primitive elements is isomorphic to $P$. Since, every $nil$-closed, noetherian, integral, unstable algebra of transcendence degree $\dim(W)$ is isomorphic to some $H^*(W)^\mathcal{G}$, we need to be able to identify the primitive elements associated with a regular central element $(V,\phi)$ of $H^*(W)^\mathcal{G}$.\\

In the second subsection, we consider $H^*(W)^\mathcal{G}$ and an inclusion $\delta$ from some vector space $V$ to $W$, such that $(V,\delta^*\phi)\in\textbf{C}(H^*(W)^\mathcal{G})$ for $\phi$ the inclusion of $H^*(W)^\mathcal{G}$ in $H^*(W)$. Then, we prove that $P(H^*(W)^\mathcal{G},\delta^*\phi)$ is a $nil$-closed and noetherian sub-algebra of $\pi^*(H^*(W/\Ima(\delta)))$ for $\pi$ the projection from $W$ to $W/\Ima(\delta)$. Since $\pi^*$ is injective, there exists $H^*(W/\Ima(\delta))^{\mathcal{G}'}\subset H^*(W/\Ima(\delta))$ such that $P(H^*(W)^\mathcal{G},\delta^*\phi)=\pi^*(H^*(W/\Ima(\delta))^{\mathcal{G}'})$. We conclude this sub-section by explaining how to determine $\mathcal{G}'$ from $\mathcal{G}$.\\

Finally, in the last sub-section, we give examples of how to answer the following question: given $\mathcal{G}'$ a groupoid with the restriction property on sub-spaces of $W/V$, with $W$ a finite dimensional vector space and $V$ a sub-space of $W$, what are the groupoids $\mathcal{G}$ with the restriction property on sub-spaces of $W$ such that $(V,\iota_V^*\phi)\in\textbf{C}(H^*(W)^\mathcal{G})$, for $\iota_V$ the inclusion of $V$ in $W$, and such that the primitive elements of $H^*(W)^\mathcal{G}$ for the $H^*(V)$-comodule structure induced by $\iota_V^*\phi$ are $\pi^*(H^*(W/V)^{\mathcal{G}'})$.
 
 \subsection{Noetherian, $nil$-closed, unstable sub-algebras of $H^*(W)$}
 
In this sub-section, we give an explicit one-to-one correspondence between the groupoids on sub spaces of $W$ satisfying the restriction property and the noetherian, $nil$-closed, unstable sub algebra of $H^*(W)$ of transcendence degree $\dim(W)$.

\begin{theorem}\label{principal4}
For all $W\in\E$, there is a one-to-one correspondence between the set of $nil$-closed and noetherian sub-algebras of $H^*(W)$ whose transcendence degree is $\dim(W)$ and the set of groupoids with the restriction property, whose objects are the sub-vector spaces of $W$.
\end{theorem}

\begin{proof}
By Theorem \ref{principal1}, there is a one-to-one correspondence between isomorphism classes in $_W\SVI$ and the set of groupoids satisfying the assumption. Thus, we have to justify that the set of $nil$-closed and noetherian sub-algebras of $H^*(W)$ of transcendence degree $\dim(W)$ are in one-to-one correspondence with isomorphism classes in $_W\SVI$. Let $K$ be a $nil$-closed, noetherian, sub-algebra of $H^*(W)$ whose transcendence degree is $\dim(W)$. Then, for $\phi_K$ the inclusion of $K$ in $H^*(W)$, since the transcendence degree of $K$ is $\dim(W)$, by the Theorem of Adams-Wilkerson, $\phi_K$ is regular. Since $K$ is noetherian, this implies that we can restrict $\phi_K^*\ :\ \Hom_{\E}(\_,W)\twoheadrightarrow \HomK(K,H^*(\_))$ to regular elements. For $F:=\HomKf(K,H^*(\_))$ and $q_F\ :\ \HomC(\_,W)\twoheadrightarrow F$ the restriction of $\phi_K^*$ to regular elements, we have that $\phi_K^*=\widetilde{q_F}$ and that $(F,q_F)$ is an element in $_W\SVI$. This defines a map $h$ from the set of $nil$-closed and noetherian sub-algebras of $H^*(W)$ whose transcendence degree is $\dim(W)$ to the set of isomorphism classes in $_W\SVI$.\\

Let us prove that $h$ is injective. We consider, $K$ and $K'$ two  $nil$-closed and noetherian sub-algebras of $H^*(W)$ such that, for $h(K)=(F,q_F)$ and $h(K')=(F',q_{F'})$ are isomorphic in $_W\SVI$. Then, there is an isomorphism $\eta$ from $F$ to $F'$ such that $q_{F'}=\eta\circ q_F$. By applying the functor $\widetilde{(\_)}$, we get the following commutative diagram in $\SVF$: $$\xymatrix{ & \HomK(K,H^*(\_))\ar[dd]^-{\Tilde{\eta}}\\
\HomK(H^*(W),H^*(\_)) \ar[ru]^-{\phi_K^*}\ar[rd]_-{\phi_{K'}^*} &  \\
 & \HomK(K',H^*(\_)).}$$ Finally, by applying the functor $m\circ\mathcal{L}$, we get the following commutative diagram:
$$\xymatrix{l_1(K)\ar@{^{(}->}[rd]^-{l_1(\phi_K)} &\\
& l_1(H^*(W))\\
l_1(K')\ar[uu]^-{m\circ\mathcal{L}(\Tilde{\eta})}\ar@{^{(}->}[ru]_-{l_1(\phi_{K'})} &\ ,}$$ where, $m\circ\mathcal{L}(\Tilde{\eta})$ is an isomorphism. Since $K$, $K'$ and $H^*(W)$ are $nil$-closed, this implies that there is an isomorphism from $K$ to $K'$ such that the following diagram commutes: 
$$\xymatrix{K\ar@{^{(}->}[rd]^-{\phi_K}\ar[dd]_-{\cong} &\\
& H^*(W)\\
K'\ar@{^(->}[ru]^{\phi_{K'}} &\ ,}$$ so that $K$ and $K'$ are the same sub-algebra of $H^*(W)$, hence $h$ is injective.\\

We conclude, by exhibiting a right inverse of $h$, for $(F,q_F)\in\ _W\SVI$, we consider the injection $m\circ\mathcal{L}(\Tilde{F})\overset{m\circ\mathcal{L}(\widetilde{q_F})}{\hookrightarrow} H^*(W)$ and we define $j(F,q_F)$ the image of $m\circ\mathcal{L}(\widetilde{q_F})$ in $H^*(W)$, which does not depend on the choice of $(F,q_F)$ in its isomorphism class. By construction, $j(F,q_F)$ is $nil$-closed, and $g(j(F,q_F))\cong\Tilde{F}$ is noetherian so, by Proposition \ref{noethnoeth}, $j(F,q_F)$ is noetherian. $j$ is obviously a right inverse of $h$; since $h$ is injective, it is a bijection.  
\end{proof}

\begin{definition}
For $\mathcal{G}$ a groupoid whose objects are the sub-vector spaces of $W$ and which has the restriction property, for $F:=\HomC(\_,W)/\sim_\mathcal{G}$ and for $q_F$ the canonical surjection from $\HomC(\_,W)$ to $F$, $H^*(W)^\mathcal{G}$ is the image of the map $$m\circ\mathcal{L}(\Tilde{q}_F)\ :\ m\circ\mathcal{L}(\Tilde{F})\hookrightarrow H^*(W).$$
\end{definition}

\begin{remark}
 $H^*(W)^\mathcal{G}$ is the unique $nil$-closed, noetherian, sub-algebra of $H^*(W)$, whose transcendence degree is $\dim(W)$ and such that $\mathcal{G}_{(F,q_F)}=\mathcal{G}$, for $F:=\HomKf(H^*(W)^\mathcal{G},H^*(\_))$ and $q_F$ the natural surjection from $\HomC(\_,W)$ to $\HomKf(H^*(W)^\mathcal{G},H^*(\_))$ induced by the inclusion from $H^*(W)^\mathcal{G}$ to $H^*(W)$.
 
 Furthermore, $\mathcal{G}\mapsto H^*(W)^\mathcal{G}$ defines a contravariant functor between $\textit{Groupoid}(W)$ (see Notation \ref{frakfrak}) and the poset of $nil$-closed, noetherian, sub-algebras of $H^*(W)$, whose transcendence degrees are $\dim(W)$, ordered by inclusion.
\end{remark}

\begin{corollary}\label{principal5}
Any $nil$-closed, integral, noetherian, unstable, algebra whose transcendence degree is equal to $\dim(W)$ is isomorphic to $H^*(W)^{\mathcal{G}}$ for some $\mathcal{G}$.
\end{corollary}

\begin{proof}
It is a reformulation of the theorem of Adams-Wilkerson using Theorem \ref{principal4}.
\end{proof}

\begin{example}
For $G$ a sub-group of $\text{Gl}(W)$, $H^*(W)^{\mathfrak{g}(G)}=H^*(W)^G$, for $H^*(W)^G$ the algebra of invariant element of $H^*(W)$ under the action of $G$.
\end{example}

Let us identify precisely the sub-algebra $H^*(W)^\mathcal{G}$ of $H^*(W)$.

\begin{proposition}\label{defdif}
Let $\mathcal{G}$ be a groupoid whose objects are the sub-vector spaces of $W$ with the restriction property. Then, $$H^*(W)^\mathcal{G}=\{x\in H^*(W)\ ;\ \alpha^*\iota_{U'}^*(x)=\iota_U^*(x)\text{ for all }\alpha\in\mathcal{G}(U,U')\}.$$
\end{proposition}

\begin{proof}
Let $\phi$ be the inclusion of $H^*(W)^\mathcal{G}$ in $H^*(W)$ and let $K(\mathcal{G})=\{x\in H^*(W)\ ;\ \alpha^*\iota_{U'}^*(x)=\iota_U^*(x)\text{ for all }\alpha\in\mathcal{G}(U,U')\}$. By definition, $\alpha^*\iota_{U'}^*\phi=\iota_U^*\phi\text{ for all }\alpha\in\mathcal{G}(U,U')$ and for all sub-spaces $U$ and $U'$ of $W$. Then, $$H^*(W)^\mathcal{G}\subset K(\mathcal{G}).$$ Furthermore, the inclusion from $K(\mathcal{G})$ to $H^*(W)$ induces a surjection from $\HomC(\_,W)$ to $g(K(\mathcal{G}))$ which factorises through an isomorphism from $\HomC(\_,W)/\sim_\mathcal{G}$ to $g(K(\mathcal{G}))$. The fact that the surjection factorises through $\HomC(\_,W)/\sim_\mathcal{G}\twoheadrightarrow g(K(\mathcal{G}))$ is a direct consequence of the definitions of $K(\mathcal{G})$ and $\sim_\mathcal{G}$, and the fact that this morphism is injective, comes from the fact that the isomorphism from $\HomC(\_,W)/\sim_\mathcal{G}$ to $g(H^*(W)^\mathcal{G})$ factorises as the following diagram: $$\HomC(\_,W)/\sim_\mathcal{G}\rightarrow g(K(\mathcal{G}))\rightarrow g(H^*(W)^\mathcal{G}).$$\\

Then, by applying the functor $m\circ\mathcal{L}$ to the last diagram, we get the bottom line of the following:

$$\xymatrix{H^*(W)^\mathcal{G}\ar@{^(->}[r]\ar[d]^-{\eta_{H^*(W)^\mathcal{G}}} & K(\mathcal{G})\ar@{^(->}[r]\ar[d]^-{\eta_{K(\mathcal{G})}} & H^*(W)\ar[d]^-{\eta_{H^*(W)}}\\
l_1(H^*(W)^\mathcal{G})\ar[r]^-\cong & l_1(K(\mathcal{G}))\ar@{^(->}[r]& l_1(H^*(W)),}$$ where $\eta$ denotes the unit of the adjunction between $f$ and $m$. Then, since  $H^*(W)^\mathcal{G}$ and $H^*(W)$ are $nil$-closed, $\eta_{H^*(W)^\mathcal{G}}$ and $\eta_{H^*(W)}$ are isomorphisms. Furthermore, $K(\mathcal{G})$ is a sub unstable algebra of $H^*(W)$, hence it does not contains any nilpotent sub module, and $\eta_{K(\mathcal{G})}$ is injective. Then, the commutativity of the diagram implies that $\eta_{K(\mathcal{G})}$ is an isomorphism, and therefore that $H^*(W)^\mathcal{G}=K(\mathcal{G})$.
\end{proof}

\begin{definition}
 For $g\in \text{Gl}(W)$, and $\mathcal{G}\in\textit{Groupoid}(W)$, $g\cdot\mathcal{G}$ is the groupoid in $\textit{Groupoid}(W)$ defined by $\beta\in g\cdot\mathcal{G}(R,R')$, for $\beta$ an isomorphisms between subspaces $R$ and $R'$ of $W$, if there exist $\alpha\in\mathcal{G}(U,U')$ for $U=g^{-1}(R)$ and $U'=g^{-1}(R')$, such that the following diagram commutes: $$\xymatrix{U\ar[r]^-{g|_U^R}\ar[d]_\alpha & R\ar[d]^\beta\\
 U'\ar[r]_-{g|_{U'}^{R'}}& R'.}$$
 
 This defines a poset preserving action of $\text{Gl}(W)$ on $\textit{Groupoid}(W)$.
\end{definition}

\begin{remark}
 This action generalises the action by conjugation on $\textit{Group}(W)$. Indeed, for $G$ a subgroup of $\text{Gl}(W)$ and $g\in \text{Gl}(W)$, $g\cdot\mathfrak{g}(G)=\mathfrak{g}(gGg^{-1})$.
\end{remark}

\begin{proposition}\label{goutzi}
For $g\in\text{Gl}(W)$ and $\mathcal{G}\in\textit{Groupoid}(W)$, $H^*(W)^{g\cdot\mathcal{G}}=(g^{-1})^*(H^*(W)^\mathcal{G})$.
\end{proposition}

\begin{proof}
This is a direct consequence of Proposition \ref{defdif}.
\end{proof}

\begin{remark}\label{gatouzou}
We want to notice that the $(H^*(W)^{\mathcal{G}})_{\mathcal{G}\in \textit{Groupoid}(W)}$ does not constitute a minimal list for representing elements of isomorphism classes of $nil$-closed, integral and noetherian unstable algebras of transcendence degree $\dim(W)$. For $g\in\text{Gl}(W)$ and $\mathcal{G}\in\textit{Groupoid}(W)$, $g\cdot\mathcal{G}$ needs not to be equal to $\mathcal{G}$, but, by Proposition \ref{goutzi}, $H^*(W)^\mathcal{G}\cong H^*(W)^{g\cdot\mathcal{G}}$.\\

Conversely, since the inclusion of $H^*(W)^\mathcal{G}$ in $H^*(W)$ induces a surjection from\newline $\HomK(H^*(W),H^*(W))$ to $\HomK(H^*(W)^\mathcal{G},H^*(W))$, and since $g\mapsto g^*$ induces an isomorphism between $\HomK(H^*(W),H^*(W))$ and $\text{Gl}(W)$, we have that if $H^*(W)^\mathcal{G}\cong H^*(W)^\mathcal{H}$, there exists $g\in\text{Gl}(W)$ such that $H^*(W)^\mathcal{H}=(g^{-1})^*(H^*(W)^\mathcal{G})$. By Proposition \ref{goutzi}, $\mathcal{H}=g\cdot \mathcal{G}$.
\end{remark}

\subsection{Centrality and primitive elements of $H^*(W)^\mathcal{G}$}

Throughout this sub-section, we fix $V$ and $W$ two objects in $\E$, as well as an injection $\delta$ from $V$ to $W$. \\

We consider $K$ a $nil$-closed, noetherian unstable sub algebra of $H^*(W)$ of transcendence degree $\dim(W)$, such that $(V,\delta^*\phi)\in\textbf{C}(K)$, for $\phi$ the inclusion of $K$ in $H^*(W)$. We start by explaining why the $H^*(V)$-comodule structure on $K$ induced by $\delta^*\phi$ is induced from the $H^*(V)$-comodule structure on $H^*(W)$ given by $(\id_W+\delta)^*\ :\ H^*(W)\rightarrow H^*(W)\otimes H^*(V)$. \\

Then, for $K=H^*(W)^\mathcal{G}$, we explain how to determine the primitive elements of this comodule structure from $\mathcal{G}$.

\begin{proposition}\label{commutemonbrave}
Let $K$ be a noetherian unstable sub algebra of $H^*(W)$ of finite transcendence degree $\dim(W)$ such that $(V,\delta^*\phi)\in\textbf{C}(K)$, for $\phi$ the inclusion of $K$ in $H^*(W)$. The $H^*(V)$-comodule structure $\kappa$ on $K$, induced by $\delta^*\phi$ and Corollary \ref{DWCT}, fits into the following commutative diagram:
$$\xymatrix{K\ar[d]_-\phi\ar[rr]^-\kappa & & K\otimes H^*(V)\ar[d]^-{\phi\otimes\id_{H^*(V)}}\\
H^*(W)\ar[rr]_-{(\id_W+\delta)^*}& & H^*(W)\otimes H^*(V).}$$
\end{proposition}

\begin{proof}
We consider the following diagram:
$$\xymatrix{K\ar[d]_-\phi\ar[rr]^-\kappa & & K\otimes H^*(V)\ar[d]^-{\phi\otimes\id_{H^*(V)}}\\
H^*(W)& & H^*(W)\otimes H^*(V).}$$
The existence of a morphism $\psi^*$ from $H^*(W)$ to $H^*(W)\otimes H^*(V)$ which turns it into a commutative diagram is a consequence of the surjectivity of $\phi^*$ from $\HomK(H^*(W),H^*(W\oplus V))$ to $\HomK(K,H^*(W\oplus V))$. We only have to justify why we can take $\psi=\id_W+\delta$. We have that the composition of $(\phi\otimes\id_{H^*(V)})\circ\kappa$ with $\epsilon_K\otimes\id_{H^*(V)}$ is equal to $\delta^*\phi$ and that with $\id_{H^*(W)}\otimes\epsilon_{H^*(V)}$ is equal to $\phi$. Hence, since $\delta^*\phi$ is central, $(\phi\otimes\id_{H^*(V)})\circ\kappa$ is the unique element in the inverse image of $\phi$ under $\rho_{\HomK(K,H^*(\_)),(V,\delta^*\phi)}$. But $(\id_W+\delta)^*\phi$ is also in this inverse image of $\phi$, hence the diagram commutes.
\end{proof}

We consider $(\id_W+\delta)^*\ : H^*(W)\rightarrow H^*(W)\otimes H^*(V)$ which is the $H^*(V)$-comodule structure on $H^*(W)$ associated with $(V,\delta)\in\textbf{C}(H^*(W))$.

\begin{proposition}\label{gazouzou}
Let $K$ be a noetherian unstable sub algebra of $H^*(W)$ of finite transcendence degree $\dim(W)$ such that $(V,\delta^*\phi)\in\textbf{C}(K)$, for $\phi$ the inclusion of $K$ in $H^*(W)$. Then, we have a pullback diagram of the following form:
$$\xymatrix{P(K,\delta^*\phi)\ar@{^(->}[r] \ar@{^(->}[d] & K\ar@{^(->}[d]^-{\phi}\\
H^*(W/\Ima(\delta))\ar@{^(->}[r]_-{\pi^*} & H^*(W).}$$
\end{proposition}

\begin{proof}
Proposition \ref{commutemonbrave} says that the following diagram commutes:
$$\xymatrix{K\ar[d]_-\phi\ar[rr]^-\kappa & & K\otimes H^*(V)\ar[d]^-{\phi\otimes\id_{H^*(V)}}\\
H^*(W)\ar[rr]_-{(\id_W+\delta)^*} & & H^*(W)\otimes H^*(V).}$$ This means that the $H^*(V)$-comodule structure on $K$ is induced by that on $H^*(W)$. Hence, the primitive elements of $K$ are simply the primitive elements of $H^*(W)$ that are in $K$. But the comodule structure on $H^*(W)$ is the morphism $(\id_W+\delta)^*$ whose algebra of primitive elements is the image of $H^*(W/\Ima(\delta))$ under $\pi^*$, for $\pi$ the projection from $W$ to $W/\Ima(\delta)$.
\end{proof}

\begin{corollary}\label{cuitcuit}
Let $K$ be a noetherian unstable sub algebra of $H^*(W)$ of finite transcendence degree $\dim(W)$ such that $(V,\delta^*\phi)\in\textbf{C}(K)$, for $\phi$ the inclusion of $K$ in $H^*(W)$. Then, the following is a pushout diagram: $$\xymatrix{\Hom_{\E}(\_,W)\ar@{->>}[r]\ar@{->>}[d] & \HomK(K,H^*(\_))\ar@{->>}[d]\\
\Hom_{\E}(\_,W/\Ima(\delta))\ar@{->>}[r] & \HomK(P(K,\delta^*\phi),H^*(\_)).}$$
\end{corollary}

\begin{proof}
It is a direct consequence of Lemma \ref{coco limite} and of Proposition \ref{gazouzou}.
\end{proof}

We can thus identify $\HomK(P(K,\delta^*\phi),H^*(\_))$ in this context. In particular, we show that $P$ is always noetherian.

\begin{lemma}\label{simsym}
For $S$ a set, and $\sim_1$ and $\sim_2$ two equivalence relations on $S$, we denote by $\sim$ the smallest equivalence relation on $S$ (in the sense that $\{(a,b)\in S\times S\ ;\ a\sim b\}\subset S\times S$ is the smallest) such that, for all $a$ and $b$ in $S$ such that $a\sim_1 b$ or $a\sim_2 b$, $a\sim b$. Then, the following is a pushout in $\mathcal{S}et$:
$$\xymatrix{S\ar@{->>}[r]\ar@{->>}[d] & S/\sim_1\ar@{->>}[d]\\
S/\sim_2\ar@{->>}[r] & S/\sim.}$$
\end{lemma}

\begin{proof}
Let $\Sigma$ denote the pushout of $$\xymatrix{S\ar@{->>}[r]\ar@{->>}[d] & S/\sim_1\\
S/\sim_2 &\ . }$$ Then, for $s\ :\ S\rightarrow \Sigma$ the composition of the projection from $S$ to $S/\sim_1$ with the surjective application $S/\sim_1\rightarrow \Sigma$, $s$ is surjective. We define $\sim'$ the equivalence relation on $S$ defined by $a\sim' b$ if and only if $s(a)=s(b)$. $\Sigma$ is isomorphic in $\mathcal{S}et$ with $S/\sim'$ and we will show that $\sim'=\sim$.\\

By commutativity of the pushout diagram, for $a$ and $b$ in $S$ such that $a\sim_1 b$ or $a\sim_2 b$, $s(a)=s(b)$. Suppose that $\sim'$ is not the smallest such equivalence relation. Then, there exists $x$ and $y$ with $x\sim' y$ and an equivalence relation $\sim''$, satisfying that for $a$ and $b$ such that $a\sim_1 b$ or $a\sim_2 b$, $a\sim''b$, and such that $x$ is not equivalent to $y$ for $\sim''$. Then, the following diagram is commutative:
$$\xymatrix{S\ar@{->>}[r]\ar@{->>}[d] & S/\sim_1\ar@{->>}[d]\\
S/\sim_2\ar@{->>}[r] & S/\sim'' ,}$$ and factorise by a morphism $S/\sim'\rightarrow S/\sim''$. This is a contradiction, so $\sim'=\sim$.
\end{proof}

\begin{remark}
For $\sim_1$ and $\sim_2$ as in Lemma \ref{simsym}, and for $S$ finite, the smallest equivalence relation $\sim$ on $S$  such that, for all $a$ and $b$ in $S$ such that $a\sim_1 b$ or $a\sim_2 b$, $a\sim b$, is the equivalence relation defined by $a\sim b$ if there is a finite family $(s_i)_{i\in\left[|1,n|\right]}$ of objects in $S$ such that: \begin{enumerate}
    \item $s_1=a$,
    \item $s_n=b$,
    \item for all $1\leq i\leq n$, if $i$ is odd $s_i\sim_1s_{i+1}$ and if $i$ is even $s_i\sim_2s_{i+1}$.
\end{enumerate} 
\end{remark}

We deduce the following proposition.

\begin{proposition}\label{caduc}
Let $K$ be a noetherian unstable sub algebra of $H^*(W)$ of finite transcendence degree $\dim(W)$ such that $(V,\delta^*\phi)\in\textbf{C}(K)$, for $\phi$ the inclusion of $K$ in $H^*(W)$. Then, for $\zeta$ and $\gamma$ in $\HomC(U,V)$, $\gamma^*\phi|_{P(K,\delta^*\phi)}=\zeta^*\phi|_{P(K,\delta^*\phi)}\in\HomK(P(K,\delta^*\phi),H^*(U))$ if and only if there exists a family $(\epsilon_i)_{i\in\left[|1,n\right|]}\in\Hom_{\E}(U,W)^n$ with $n\in\mathbb{N}$ greater than $1$, such that: \begin{enumerate}
    \item $\gamma=\epsilon_1$, 
    \item $\zeta=\epsilon_n$,
    \item for all $1\leq i\leq n-1$, $\epsilon_i^*\phi=\epsilon_{i+1}^*\phi$ if $i$ is odd and $\pi\circ\epsilon_i=\pi\circ\epsilon_{i+1}$ if $i$ is even.
\end{enumerate}
\end{proposition}

\begin{proof}
By Corollary \ref{cuitcuit}, the following is a pushout: 
$$\xymatrix{\Hom_{\E}(\_,W)\ar@{->>}[r]^-{k}\ar@{->>}[d]_-{\pi_*} & \HomK(K,H^*(\_))\ar@{->>}[d]^-p\\
\Hom_{\E}(\_,W/\Ima(\delta))\ar@{->>}[r] & \HomK(P(K,\delta^*\phi),H^*(\_)),}$$ where $k$ maps $\zeta\ :\ U\rightarrow W$ to $\zeta^*\phi\ :\ K\rightarrow H^*(U)$, $p$ maps $\psi\ :\ K\rightarrow H^*(U)$ to $\psi|_{P(K,\delta^*\phi)}$ and $\pi_*$ maps $\zeta\ :\ U\rightarrow W$ to $\pi\circ\zeta\ :\ U\rightarrow W/\Ima(\delta)$. Then, by Lemma \ref{simsym}, $p\circ k(\zeta)=p\circ k(\gamma)$ if and only if there exists a family $(\epsilon_i)_{i\in\left[|1,n\right|]}\in\Hom_{\E}(U,W)^n$ with $n\in\mathbb{N}$ greater than $1$, such that $\gamma=\epsilon_1$, $\zeta=\epsilon_n$ and for all $1\leq i\leq n-1$, $k(\epsilon_i)=k(\epsilon_{i+1})$ if $i$ is odd and $\pi_*(\epsilon_i)=\pi_*(\epsilon_{i+1})$ if $i$ is even. 
\end{proof}

\begin{corollary}\label{blabla}
Let $K$ be a noetherian unstable sub algebra of $H^*(W)$ of finite transcendence degree $\dim(W)$ such that $(V,\delta^*\phi)\in\textbf{C}(K)$, for $\phi$ the inclusion of $K$ in $H^*(W)$. Then, for $\zeta\in\Hom_{\E}(U,W)$, $\ker(\zeta^*\phi|_{P(K,\delta^*\phi)})=\ker(\pi\circ\zeta)$.
\end{corollary}

\begin{proof}
Let $\zeta_0\in\Hom_{\E}(U/\ker(\zeta^*\phi|_{P(K,\delta^*\phi)},W)$ such that $\zeta^*\phi|_{P(K,\delta^*\phi)}=\pi_U^*\zeta_0^*\phi|_{P(K,\delta^*\phi)}$, with $\pi_U$ the projection from $U$ to $U/\ker(\zeta^*\phi|_{P(K,\delta^*\phi)})$. Let $\epsilon_1=\zeta_0\circ\pi_U$, $\epsilon_n=\zeta$ and for all $i$ $\epsilon_i^*\phi=\epsilon_{i+1}^*\phi$ if $i$ is odd and $\pi\circ\epsilon_i=\pi\circ\epsilon_{i+1}$ if $i$ is even. Then, since $\HomK(K,H^*(\_))$ is noetherian, $\ker(\pi\circ\epsilon_i)=\ker(\pi\circ\epsilon_{i+1})$ for all $1\leq i\leq n-1$. Hence, $\ker(\pi\circ\zeta)=\ker(\pi\circ\zeta_0\circ\pi_U)=\ker(\zeta^*\phi|_{P(K,\delta^*\phi)})$.
\end{proof}

\begin{corollary}
Let $K$ be a noetherian unstable sub algebra of $H^*(W)$ of finite transcendence degree $\dim(W)$ such that $(V,\delta^*\phi)\in\textbf{C}(K)$, for $\phi$ the inclusion of $K$ in $H^*(W)$. Then, \newline$\HomK(P(K,\delta^*\phi),H^*(\_))$ is noetherian.
\end{corollary}

\begin{proof}
The fact that $\HomK(P(K,\delta^*\phi),H^*(\_))$ is finite is obvious. Let $\zeta^*\phi|_{P(K,\delta^*\phi)}$ in \newline$\HomK(P(K,\delta^*\phi),H^*(U))$ and let $\alpha$ be a morphism from a vector space $Y$ to $U$. Then, by Corollary \ref{blabla} $$\ker(\alpha^*\zeta^*\phi|_{P(K,\delta^*\phi)})=\ker(\pi\circ\zeta\circ\alpha).$$ Since $\alpha$ is injective, this is equal to $$\alpha^{-1}(\ker(\pi\circ\zeta))=\alpha^{-1}(\ker(\zeta^*\phi|_{P(K,\delta^*\phi)})).$$
\end{proof}

\begin{theorem}\label{principal7}
Let $K$ be a noetherian unstable sub algebra of $H^*(W)$ of finite transcendence degree $\dim(W)$ such that $(V,\delta^*\phi)\in\textbf{C}(K)$, for $\phi$ the inclusion of $K$ in $H^*(W)$. Then, $P(K,\delta^*\phi)$ is $nil$-closed and noetherian.
\end{theorem}

\begin{proof}
Since, $P(K,\delta^*\phi)$ is the kernel of $\kappa-\id_K\otimes 1$ from $K$ to $K\otimes H^*(V)$ which are $nil$-closed, for $\kappa$ the comodule structure of $K$ associated with $\delta^*\phi$, and since $f$ is exact and $m$ is left-exact, the following is an exact sequence:
 $$0\rightarrow l_1(P(K,\delta^*\phi))\rightarrow l_1(K)\overset{l_1(\kappa-\id_K\otimes 1)}{\longrightarrow} l_1(K\otimes H^*(V)).$$ Therefore, since $K$ is $nil$-closed, $P(K,\delta^*\phi)$ is also $nil$-closed. Then, the noetherianity of $P(K,\delta^*\phi)$ is a consequence of the noetherianity of $\HomK(P(K,\delta^*\phi),H^*(\_))$ and of Proposition \ref{noethnoeth}.
\end{proof}

\begin{remark}
We have identified $P(K,\delta^*\phi)$ with a sub-algebra of $H^*(W/\Ima(\delta))$. Furthermore, we proved that $P(K,\delta^*\phi)$ is $nil$-closed and noetherian, and (because we took $\delta$ to be an injection) Corollary \ref{blabla} implies that the inclusion from $P(K,\delta^*\phi)$ into $H^*(W/\Ima(\delta))$ is regular. Therefore, by Theorem \ref{principal4}, $P(K,\delta^*\phi)$ has the form $H^*(W/\Ima(\delta))^{\mathcal{G}'}$, for some $\mathcal{G}'\in\textit{Groupoid}(W/\Ima(\delta))$.
\end{remark}

This leads to the following question: for $W$ and $V$ in $\E$, for $\delta$ an inclusion from $V$ to $W$ and for $\mathcal{G}'$ a groupoid with the restriction property and whose objects are the sub-spaces of $(W/\Ima(\delta))$, which are the groupoids $\mathcal{G}$ with the restriction property and whose objects are the sub-spaces of $W$, such that \begin{enumerate}
    \item $H^*(W)^\mathcal{G}$ is a sub $H^*(V)$-comodule of $H^*(W)$ for the comodule struture induced by $\delta$,
    \item the intersection of $H^*(W)^\mathcal{G}$ with $\pi^*(H^*(W/\Ima(\delta)))$ is the image under $\pi^*\ :\ H^*(W/\Ima(\delta))\rightarrow H^*(W)$ of $H^*(W/\Ima(\delta))^{\mathcal{G}'}$.
\end{enumerate}

\begin{remark}
$H^*(W)^\mathcal{G}$ is a sub $H^*(V)$-comodule of $H^*(W)$ for the comodule struture induced by $\delta$ if and only if $\mathcal{G}$ satisfies the two conditions of theorem \ref{principal2}.
\end{remark}

We recall that, from the begining of this sub-section, $V$ and $W$ are fixed objects of $\E$ and $\delta$ a fixed injective morphism from $V$ to $W$.

\begin{theorem}\label{principal6}
Let $\mathcal{G}$ be a groupoid with the restriction property and whose objects are the sub-vector spaces of $W$, such that $H^*(W)^\mathcal{G}$ is a sub $H^*(V)$-comodule of $(H^*(W),(\id_W+\delta)^*)$. For $\mathcal{G'}$ the only groupoid whose objects are sub spaces of $W/\Ima(\delta)$ which satisfies that $\pi^*(H^*(W/\Ima(\delta))^{\mathcal{G}'})$ is the algebra of primitive elements of $H^*(W)^\mathcal{G}$, the two following conditions are equivalent:

\begin{enumerate}
    \item\label{ain} $\alpha\in\mathcal{G}'(U,U')$, where $U$ and $U'$ are sub-vector spaces of $W/\Ima(\delta)$ and $\alpha$ is an isomorphism from $U$ to $U'$,
    \item\label{deu} there exists $N$ and $N'$ sub spaces of $W$ such that $\pi$ induce isomorphisms from $N$ and $N'$ to $U$ and $U'$, as well as an element $\beta\in\mathcal{G}(N,N')$ such that $\alpha=\pi|_{N'}^{U'}\circ\beta\circ (\pi|_N^U)^{-1}$.
\end{enumerate}
\end{theorem}

\begin{proof}

We consider the pushout diagram of Corollary \ref{cuitcuit}: $$\xymatrix{\Hom_{\E}(\_,W)\ar@{->>}[r]^-k\ar@{->>}[d]_-{\pi^*} & \HomK(H^*(W)^\mathcal{G},H^*(\_))\ar@{->>}[d]^-q\\
\Hom_{\E}(\_,W/\Ima(\delta))\ar@{->>}[r]_-p & \HomK(H^*(W/\Ima(\delta))^{\mathcal{G}'},H^*(\_)),}$$
where $\pi_*$ maps $\gamma\ :\ U\rightarrow W$ to $\pi\circ\gamma$, $k$ maps $\gamma$ to $\gamma^*\phi_\mathcal{G}$ for $\phi_\mathcal{G}$ the inclusion from $H^*(W)^\mathcal{G}$, $q$ maps $\psi\ :\ H^*(W)^\mathcal{G}\rightarrow H^*(U)$ to $\psi|_{\pi^*(H^*(W/\Ima(\delta)))}$ the restriction of $\psi$ to $\pi^*(H^*(W/\Ima(\delta)))$ and, finally, $p$ maps $\zeta$ from $U$ to $W/\Ima(\delta)$ to $\zeta^*\phi_{\mathcal{G}'}$, for $\phi_{\mathcal{G}'}$ the inclusion of $H^*(W/\Ima(\delta))^{\mathcal{G}'}$ into $H^*(W/\Ima(\delta))$.\\

We fix a section $s$ from $W/\Ima(\delta)$ to $W$. Since $\pi\circ s=\id_{W/\Ima(\delta)}$, $\pi_*(s)=\id_{W/\Ima(\delta)}$. Then, by commutativity of the pushout diagram, we have $\phi_{\mathcal{G}'}=q(s^*\phi_{\mathcal{G}})=s^*\phi_\mathcal{G}|_{\pi^*(H^*(W/\Ima(\delta)))}$.\\

By construction, there are natural isomorphisms $\HomKf(H^*(W)^\mathcal{G},H^*(\_))\cong \HomC(\_,W)/\sim_\mathcal{G}$ and $\HomKf(H^*(W/\Ima(\delta))^{\mathcal{G}'},H^*(\_))\cong\HomC(\_,W/\Ima(\delta))/\sim_{\mathcal{G}'}$. These are the isomorphisms that map $\phi_{\mathcal{G}}$ to $\left[\id_W\right]_\mathcal{G}$ and $\phi_{\mathcal{G}'}$ to $\left[\id_{W/\Ima(\delta)}\right]_{\mathcal{G}'}$ respectively.\\ 

Let us first prove \ref{deu}) $\Rightarrow$ \ref{ain}). We consider, $\beta\in\mathcal{G}(N,N')$ such that $\pi$ induces isomorphisms $\pi|_N^U$ and $\pi|_{N'}^{U'}$ between $N$ and $U$ and between $N'$ and $U'$. Let $\alpha$ be an isomorphism from $U$ to $U'$ such that $\alpha=\pi|_{N'}^{U'}\circ\beta\circ (\pi|_N^U)^{-1}$. Then, $(\pi|_{N'}^{U'})^{-1}\circ\alpha=\beta\circ (\pi|_N^U)^{-1}$. Therefore, 
\begin{align*}
    \alpha^*((\pi|_{N'}^{U'})^{-1})^*\iota_{N'}^*\phi_\mathcal{G}&=((\pi|_N^U)^{-1})^*\beta^*\iota_{N'}^*\phi_\mathcal{G}\\
    &=((\pi|_N^U)^{-1})^*\iota_N^*\phi_\mathcal{G}.
\end{align*}
We can choose the section $s$ in such a way that $s\circ\iota_{U'}=\iota_{N'}\circ ((\pi|_{N'}^{U'})^{-1})$, then $$\alpha^*\iota_{U'}^*s^*\phi_\mathcal{G}=((\pi|_N^U)^{-1})^*\iota_N^*\phi_\mathcal{G}.$$ This implies that $\alpha^*\iota_{U'}^*s^*q(\phi_\mathcal{G})=((\pi|_N^U)^{-1})^*\iota_N^*q(\phi_\mathcal{G})$. Furthermore, $\pi\circ s\circ\iota_U=\pi\circ\iota_N\circ(\pi|_N^U)^{-1}$. Hence, we also have, by Proposition \ref{caduc}, that $\iota_U^*s^*q(\phi_\mathcal{G})=((\pi|_N^U)^{-1})^*\iota_N^*q(\phi_\mathcal{G})$. Hence, $$\alpha^*\iota_{U'}^*s^*q(\phi_\mathcal{G})=\iota_U^*s^*q(\phi_\mathcal{G}).$$ Since $s^*q(\phi_\mathcal{G})=\left[\id_{W/\Ima(\delta)}\right]_{\mathcal{G}'}$, this implies that $\alpha\in\mathcal{G}'(U,U')$, as required.\\

Now, let us prove the far more challenging \ref{ain}) $\Rightarrow$ \ref{deu}).

We consider $\alpha\in\mathcal{G}'(U,U')$ where $U$ and $U'$ are two sub-spaces of $W/\Ima(\delta)$. Then, $$\alpha^*\iota_{U'}^*\left[\id_{W/\Ima(\delta)}\right]_{\mathcal{G}'}=\iota_U^*\left[\id_{W/\Ima(\delta)}\right]_{\mathcal{G}'},$$ or, equivalently, $$\left[\iota_{U'}\circ\alpha\right]_{\mathcal{G}'}=\left[\iota_U\right]_{\mathcal{G}'}.$$ 
By Proposition \ref{caduc}, since $\left[\id_{W/\Ima(\delta)}\right]_{\mathcal{G}'}$ is the image of $s^*q(\phi_\mathcal{G})$ under the isomorphism $$\HomKf(H^*(W/\Ima(\delta))^{\mathcal{G}'},H^*(\_))\cong\HomC(\_,W/\Ima(\delta))/\sim_{\mathcal{G}'},$$ we have that $\left[\zeta\right]_{\mathcal{G}'}=\left[\gamma\right]_{\mathcal{G}'}$, if and only if there exists a family $(\epsilon_i)_{i\in\left[|1,n\right|]}\in\Hom_{\E}(U,W)^n$ with $n\in\mathbb{N}$ greater than $1$, such that $s\circ\gamma=\epsilon_1$, $s\circ\zeta=\epsilon_n$ and, for all $1\leq i\leq n-1$, $\epsilon_i^*\phi_\mathcal{G}=\epsilon_{i+1}^*\phi_\mathcal{G}$ if $i$ is odd and $\pi\circ\epsilon_i=\pi\circ\epsilon_{i+1}$ if $i$ is even.\\ 

So let $(\epsilon_i)_{i\in\left[|1,n\right|]}\in\Hom_{\E}(U,W)^n$ be such that $\epsilon_1=s\circ\iota_U$, $\epsilon_n=s\circ\iota_{U'}\circ\alpha$ and for all $1\leq i\leq n-1$, $\epsilon_i^*\phi_\mathcal{G}=\epsilon_{i+1}^*\phi_\mathcal{G}$ if $i$ is odd and $\pi\circ\epsilon_i=\pi\circ\epsilon_{i+1}$ if $i$ is even.\\

By induction, for all $i\in\left[|1,n|\right]$, $\epsilon_i^*\phi_\mathcal{G}$ and $\pi\circ\epsilon_i$ are regular elements respectively of \newline$\HomK(H^*(W)^\mathcal{G},H^*(U))$ and $\Hom_{\E}(U,W/\Ima(\delta))$. Hence, $\epsilon_i$ and $\pi\circ\epsilon_i$ are injections. For all $i$, let $N_i$ denote the image of $\epsilon_i$ in $W$, we denote also by $\Tilde{\epsilon}_i$ the corestriction of $\epsilon$ to $N_i$. Then, for $i$ odd, $\epsilon_i^*\phi_\mathcal{G}=\epsilon_{i+1}^*\phi_\mathcal{G}$ implies that there exists $\beta_i$ in $\mathcal{G}(N_i,N_{i+1})$ such that $\Tilde{\epsilon}_{i+1}=\beta_i\circ\Tilde{\epsilon}_i$.\\

We take some moment to explain subtlety in the proof. We would like, for $i$ even, to have $\epsilon_i=\epsilon_{i+1}$. Then, the composition of the $\beta_i$ with $i$ odd would give an isomorphism $\beta$ between $N_1=s(U)$ and $N_n=s(U')$ such that $\beta\in\mathcal{G}(N_1,N_n)$, since $\mathcal{G}$ is a groupoid and we would have $(s\circ\iota_{U'})|^{N_n}\circ\alpha=\beta\circ(s\circ\iota_U)|^{N_1}$. Since $(s\circ\iota_U)|^{N_1}$ and $(s\circ\iota_{U'})|^{N_n}$ are inverse isomorphisms of $\pi|_{N_1}^U$ and $\pi|_{N_n}^{U'}$, we would have $\alpha=\pi|_{N_n}^{U'}\circ\beta\circ(\pi|_{N_1}^{U})^{-1}$. If this were the case, we would have found a $\beta$ for any $N$ and $N'$ such that $\pi$ induces isomorphisms between $U$ and $N$ and between $U'$ and $N'$, and we would have done so without using the assumption that $\delta^*\phi_\mathcal{G}$ is central. Unfortunately, this naive approach fails, and $N$ and $N'$ must be chosen carefully. The hypothesis on the $\epsilon_i$ for $i$ even indicates how to modify our original $N_1$ and $N_n$ to make it work, using the centrality of $\delta^*\phi_\mathcal{G}$.\\

First notice that, since $\pi\circ\epsilon_i$ is injective for all $i$, we always have $N_i\cap\Ima(\delta)=\{0\}$. Then, the assumption that, for $i$ even, $\pi\circ\epsilon_i=\pi\circ\epsilon_{i+1}$ implies that there exists $\rho_i$ from $U$ to $W$ whose image is inside $\Ima(\delta)$ and such that $\epsilon_{i+1}=\epsilon_i+\rho_i$. Now, since $H^*(W)^\mathcal{G}$ is a sub $H^*(V)$-comodule of $(H^*(W),(\id_W+\delta)^*)$ we know that $\left[\delta\right]_\mathcal{G}$ is a central element of $\HomC(\_,W)/\sim_\mathcal{G}$. Then, by Theorem \ref{principal2}, for $i$ odd, we know that the isomorphisms $\Bar{\beta}_i$ from $N_i\oplus\Ima(\delta)$ to $N_{i+1}\oplus\Ima(\delta)$ defined by $\Bar{\beta}_i(n)=\beta_i(n)$ for $n\in N_i$ and $\Bar{\beta}_i(v)=v$ for $v\in\Ima(\delta)$ satisfy $\Bar{\beta}_i\in\mathcal{G}(N_i\oplus\Ima(\delta),N_{i+1}\oplus\Ima(\delta))$. Moreover, for $i$ even, $\pi\circ\epsilon_i=\pi\circ\epsilon_{i+1}$ implies that $N_i\oplus\Ima(\delta)=N_{i+1}\oplus\Ima(\delta)$. Then, at each even step $i$, we can ``correct" $\epsilon_{i-1}$ to get $\beta_{i-1}\circ\Tilde{\epsilon}_{i-1}=\Tilde{\epsilon}_{i+1}$ instead of $\Tilde{\epsilon}_i$. \\

For each $i\in\left[|1,n|\right]$, we define $\epsilon'_i$ (the ``corrected" $\epsilon_i$) by 
$$\epsilon'_i:=\iota_{N_i\oplus\Ima(\delta)}^{W}\circ(\Tilde{\epsilon}_i\oplus\sum\limits_{\{j\text{ even ; }i\leq j<n\}}\beta_{i\rightarrow j}^{-1}\circ\rho_j|^{N_j\oplus\Ima(\delta)}),$$ where $\beta_{i\rightarrow j}$ is the composition of all the $\Bar{\beta}_k$ with $k$ odd and $i\leq k<j$. The family $(\epsilon'_i)_{i\in\left[|1,n|\right]}$ satisfies the following:
\begin{enumerate}
    \item $\pi\circ\epsilon'_1=\iota_U$, $\pi\circ\epsilon'_n=\iota_{U'}\circ\alpha$, 
    \item for all $i$, if we denote by $N'_i$ the image of $\epsilon'_i$, then $N'_i\oplus\Ima(\delta)=N_i\oplus\Ima(\delta)$,
    \item for $i$ odd, if we denote by $\beta'_i$ the restriction of $\Bar{\beta}_i$ to $N'_i$ corestricted to $N'_{i+1}$,  $\Tilde{\epsilon}'_{i+1}=\beta'_i\circ\Tilde{\epsilon}'_i$, with $\beta'_i\in\mathcal{G}(N'_i,N'_{i+1}))$, since $\Bar{\beta}\in\mathcal{G}(N_i\oplus\Ima(\delta),N_{i+1}\oplus\Ima(\delta))$ and $\mathcal{G}$ has the restriction property,
    \item for $i$ even, $\epsilon'_i=\epsilon'_{i+1}$.
\end{enumerate} 

Then, let $N=N'_1$, $N'=N'_n$ and $\beta=(\beta'_k\circ...\circ\beta'_3\circ \beta'_1)$, where $k=n-2$ if $n$ is odd, $k=n-1$ otherwise. Then, $\beta\in\mathcal{G}(N,N')$ and $\beta\circ\Tilde{\epsilon}'_1=\Tilde{\epsilon}'_n$. Finally, $\pi\circ\epsilon'_1=\iota_U$ implies that $\Tilde{\epsilon}'_1=(\pi|_{N}^{U})^{-1}$ and $\pi\circ\epsilon'_n=\iota_{U'}\circ\alpha$ implies that $\Tilde{\epsilon}'_n=(\pi|_{N'}^{U'})^{-1}\circ\alpha$. Hence, $\alpha=\pi|_{N'}^{U'}\circ\beta\circ (\pi|_N^U)^{-1}$.
\end{proof}

\subsection{Applications}

We end this section by presenting some applications of Theorem \ref{principal6}. The first result was already known.

\begin{proposition}
Let $K$ be a noetherian, $nil$-closed, integral, unstable algebra of transcendence degree $d$. We assume that the centre of $K$ is of dimension $d$. Then, $K\cong H^*(W)$ with $\dim(W)=d$.
\end{proposition}

\begin{proof}
By Theorem \ref{principal4} and the theorem of Adams-Wilkerson $K\cong H^*(W)^\mathcal{G}$ for some groupoid $\mathcal{G}$ with the restriction property and whose objects are the sub-spaces of $W$. Then, since the centre of $K$ has dimension $\dim(W)$, up to isomorphism the centre of $K$ is given by $(W,\phi)$ where $\phi$ is the inclusion from $K$ to $H^*(W)$ induced by the theorem of Adams-Wilkerson. Then, by Theorem \ref{principal2}, for all $\alpha\in\mathcal{G}(U,U')$, where $U$ and $U'$ are sub-vector spaces of $W$, and for all $x\in U=U\cap W$, $\alpha(x)=x$. Hence, $\mathcal{G}$ is the groupoid in which the only morphisms are identities of sub-spaces of $W$. Then, $H^*(W)^\mathcal{G}=H^*(W)$.
\end{proof}

Let us now consider the case where the centre is of dimension $d-1$, for $d$ the transcendence degree of $K$.

\begin{theorem}\label{principal8}
Let $K$ be a noetherian, $nil$-closed, integral, unstable algebra of transcendence degree $d$. We assume that the centre of $K$ is of dimension $d-1$. Then, there exists $G$ a sub-group of $\text{Gl}(W)$ such that $K$ is isomorphic to the algebra of invariant elements $H^*(W)^G$ with $\dim(W)=d$. Furthermore, $G$ satisfies that the set of element $x\in W$ such that $g(x)=x$ for all $g\in G$ is a sub-vector space of $W$ of dimension $d-1$.
\end{theorem}

\begin{proof}
By Theorem \ref{principal4} and the theorem of Adams-Wilkerson $K\cong H^*(W)^\mathcal{G}$ for some groupoid $\mathcal{G}$ with the restriction property and whose objects are the sub-spaces of $W$. Then, the centre of $K$ is associated with a $H^*(V)$-comodule structure on $H^*(W)^\mathcal{G}$, with $\dim(V)=d-1$. Up to isomorphism, we can suppose that this $H^*(V)$-comodule structure is induced by $\iota_V$ the inclusion of a sub-vector space $V$ in $W$ (we can assume that the comodule structure is induced by an injection, because it is associated with the centre of $K$ which is regular by definition). By Theorem \ref{principal2}, for all $\alpha\in\mathcal{G}(U,U')$ and for all $x\in U\cap V$, $\alpha(x)=x$. Furthermore, for all sub-spaces $U$ and $U'$ of $W$, $\mathcal{G}(U,U')$ is determined from $\mathcal{G}(V\oplus U,V\oplus U')$ by the restriction property. Since we assumed $\dim(V)=\dim(W)-1$, $\mathcal{G}$ is uniquely determined by $\mathcal{G}(V,V)$ and $\mathcal{G}(W,W)$. Then, since for all $x\in V$ and for all $\alpha\in\mathcal{G}(V,V)$, $\alpha(x)=x$, $\mathcal{G}(V,V)$ is reduced to the identity of $V$. Therefore, for $g\in G:=\mathcal{G}(W,W)$, $g(x)=x$ for all $x\in V$, and $\alpha\in\mathcal{G}(U,U')$ if and only if there is $g\in G$ such that $\alpha=g_U$. We get that $H^*(W)^\mathcal{G}=H^*(W)^G$. 
\end{proof}

\begin{corollary}
Let $K$ be a noetherian, $nil$-closed, integral, unstable algebra of transcendence degree $\dim(W)$, and let $(C,\gamma)$ be its centre. Then if $K$ is not isomorphic to $H^*(W)^G$ for some sub-group $G$ of $\text{Gl}(W)$, $\dim (C)\leq d-2$.
\end{corollary}

In the next section, we will have examples where $\dim(C)=d-2$ and $K$ is not an algebra of invariant elements under some action on $H^*(W)$. 

\section{Two examples}

In the following, we focus on the case where $p=2$, so that $H^*(V)\cong \F\left[x_1,...,x_n\right]$, for $(x_1,...,x_n)$ a basis of $V^\sharp.$\\

We take $V_3=\F e_1\oplus\F e_2\oplus \F e_3$ and denote by $(x,y,z)$ the dual basis of $(e_1,e_2,e_3)$. The morphism $$(\id_{V_3}+\iota_{\F e_1})^*\ :\ H^*(V_3)\rightarrow H^*(V_3)\otimes H^*(\F e_1),$$ which maps $x$ to $x\otimes1+1\otimes x$, $y$ to $y\otimes1$ and $z$ to $z\otimes1$ defines an $H^*(\F e_1)$-comodule structure on $H^*(V_3)$ in $\K$. In each of the following examples we want to find the list of noetherian, $nil$-closed, unstable sub-algebras of $H^*(V_3)$ of transcendence degree $3$, which are sub $H^*(\F e_1)$-comodules of $H^*(V_3)$ for the comodule structure given by $(\id_{V_3}+\iota_{\F e_1})^*$
and whose algebra of primitive elements is $\pi^*(H^*(V_3/\F e_1)^{\mathcal{G}'})$ for $\pi$ the projection from $V_3$ over $V_3/\F e_1$, and for $H^*(V_3/\F e_1)^{\mathcal{G}'}$ a chosen $nil$-closed, noetherian, sub algebra of $H^*(V_3/\F e_1)$. By Theorem \ref{principal6}, $\mathcal{G}$ has to satisfy that $\alpha\in\mathcal{G'}(U,U')$ if and only if there is $\beta\in\mathcal{G}(N,N')$ such that $\pi|_{N'}\circ\beta=\alpha\circ\pi|_N$. \\

Furthermore, by Theorem \ref{principal2}, since $\left[\iota_{\F e_1}\right]_\mathcal{G}\in \HomKf(H^*(W)^\mathcal{G},H^*(\F e_1))$ is central, all morphism $\beta$ in $\mathcal{G}(N,N')$ is the restriction of a morphism $\Bar{\beta}$ in $\mathcal{G}(\F e_1 + N,\F e_1 +N')$ defined by $\Bar{\beta}(e_1)=e_1$ and $\Bar{\beta}(n)=\beta(n)$ for $n\in N$. So $\mathcal{G}$ can be deduced from its full sub-groupoid whose objects are $V_3$, $\F e_1 \oplus\F e_2$, $\F e_1 \oplus\F e_3$, $\F e_1 \oplus\F (e_2+e_3)$ and $\F e_1$. This full sub-groupoid can be represented as follows: \\
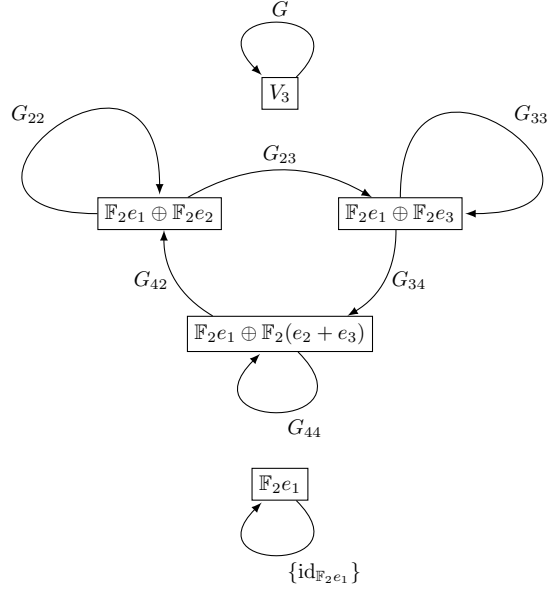
\begin{center}

\begin{tikzpicture}[scale=0.8, every node/.style={transform shape}]
\node[draw] (V) at (5,5) {$V_3$};
\node[draw] (P2) at (3,3) {$\F e_1\oplus \F e_2$};
\node[draw] (P3) at (7,3) {$\F e_1\oplus \F e_3$};
\node[draw] (P4) at (5,1) {$\F e_1\oplus \F (e_2+e_3)$};
\node[draw] (P1) at (5,-1.5) {$\F e_1$};
\draw[->,>=latex] (V) edge[in=135, out=45, loop] node[above]{$G$} (V) ;
\draw[->,>=latex] (P2) edge[in=90, out=180, loop] node[above left]{$G_{22}$} (P2) ;
\draw[->,>=latex] (P3) edge[in=0, out=90, loop] node[above right]{$G_{33}$} (P3) ;
\draw[->,>=latex] (P4) edge[in=-135, out=-45, loop] node[below right]{$G_{44}$} (P4) ;
\draw[->,>=latex] (P2) to[bend left] node[above]{$G_{23}$} (P3) ;
\draw[->,>=latex] (P3) to[bend left] node[right]{$G_{34}$} (P4) ;
\draw[->,>=latex] (P4) to[bend left] node[left]{$G_{42}$} (P2) ;
\draw[->,>=latex] (P1) edge[in=-135, out=-45, loop] node[below right]{$\{\id_{\F e_1}\}$} (P1) ;

\end{tikzpicture}
\end{center}
Where $G$ denotes $\mathcal{G}(V_3,V_3)$ and $G_{ij}$ denotes $\mathcal{G}(\F e_1\oplus\F e_i,\F e_1\oplus\F e_j)$ with $e_4:=e_2+e_3$, and where we omitted $G_{43}$, $G_{24}$, etc, which can be deduced by composition and inversion of the already given sets of morphisms.  \\

We consider first the case where $\mathcal{G}'$ contains only identities.

\begin{proposition}\label{example627}
There are exactly $15$ $nil$-closed, noetherian, unstable sub-algebras $K$ of $H^*(V_3)$ of transcendence degree $3$, which are sub $\F\left[x\right]$-comodules of $H^*(V_3)$, for the comodule structure that maps $x$ to $x\otimes1+1\otimes x$, $y$ to $y\otimes 1$ and $z$ to $z\otimes 1$, and such that the algebra of primitive elements of $K$ is $\F\left[y,z\right]$. These are :

\begin{enumerate}
    \item $\F\left[y,z,x(x+y)(x+z)(x+y+z)\right],$
    \item $\F\left[y,z,x(x+y)\right],$
    \item $\F\left[y,z,x(x+z)\right],$
    \item $\F\left[y,z,x(x+y+z)\right],$
    \item $\F\left[y,z,x(x+y)(x+z)(x+y+z)\right]+\F\left[y,z,x(x+y)\right]y,$
    \item $\F\left[y,z,x(x+y)(x+z)(x+y+z)\right]+\F\left[y,z,x(x+z)\right]z,$
    \item $\F\left[y,z,x(x+y)(x+z)(x+y+z)\right]+\F\left[y,z,x(x+y+z)\right](y+z),$
    \item $\F\left[x,y,z\right],$
    \item $\F\left[x,y,z\right](y+z)\oplus\F\left[z,x(x+z)\right],$
    \item $\F\left[x,y,z\right]y\oplus\F\left[z,x(x+z)\right],$
    \item $\F\left[x,y,z\right]z\oplus\F\left[y,x(x+y)\right],$
    \item $\F\left[z,x(x+z)\right]\oplus \F\left[x,y,z\right](y+z)y\oplus \F\left[y,x(x+y)\right]y,$
    \item $\F\left[y,x(x+y)\right]\oplus \F\left[x,y,z\right](y+z)z\oplus \F\left[z,x(x+z)\right]z,$
    \item $\F\left[(y+z),x(x+y+z)\right]\oplus \F\left[x,y,z\right]yz\oplus \F\left[z,x(x+z)\right]z,$
    \item $\F\left[z,x(x+z)\right]z\oplus\F\left[y,x(x+y)\right]y\oplus\F\left[y,x(x+y)\right](y+z)y\oplus \F\left[x,y,z\right](y+z)yz .$
     
\end{enumerate}
\end{proposition}

\begin{proof} 
In this case, $\F\left[y,z\right]=\pi^*(H^*(V_3/\F e_1))$. Then, the groupoid $\mathcal{G}'$ is the trivial groupoid, which contains only identities. Therefore, by Theorem \ref{principal6}, for all $\beta\in\mathcal{G}(N,N')$, $\pi\circ\beta=\pi$. For $\beta\in G_{ij}$ with $i$ and $j$ in $\{2,3,4\}$, $\pi\circ\beta(e_i)=\pi(e_i)$ which is not possible if $i\neq j$. Hence, $G_{ij}=\emptyset$ if $i\neq j$. We only have to determine $G$, $G_{22}$, $G_{33}$ and $G_{44}$. 
We have two condition on $\beta\in G$ or $\beta\in G_{ii}$, the first is that $\beta(e_1)=e_1$ (by Theorem \ref{principal2} and the centrality of $\left[\iota_{\F e_1}\right]_\mathcal{G}$) and the second that $\pi\circ\beta=\pi$ (by Theorem \ref{principal6} and the hypothesis on $\mathcal{G}'$). Therefore, $\beta$ has a block matrix of the following form:
$$\left(\begin{array}{ccc}
     \id_{\F e_1}& \hat{\beta}  \\
     0 &  \id_{N}
\end{array}\right),$$ where $N\in\{\F e_2,\F e_3,  \F (e_2+e_3),\F e_2\oplus\F e_3 \}$ and $\hat{\beta}$ is a morphism from $N$ to $\F e_1$. Hence, the set of matrices of morphisms in $G$ in the basis $(e_1,e_2,e_3)$ is a sub group of $$H=\{\left(\begin{array}{ccc}
     1 & 0 & 0  \\
     0 & 1 & 0\\
     0 & 0 & 1
\end{array}\right),\left(\begin{array}{ccc}
     1 & 1 & 0  \\
     0 & 1 & 0\\
     0 & 0 & 1
\end{array}\right),\left(\begin{array}{ccc}
     1 & 0 & 1  \\
     0 & 1 & 0\\
     0 & 0 & 1
\end{array}\right),\left(\begin{array}{ccc}
     1 & 1 & 1  \\
     0 & 1 & 0\\
     0 & 0 & 1
\end{array}\right)\}.$$ 

There are five possibilities for $G$, $G_1=\{\id_{V_3}\}$, $G_2=<\left(\begin{array}{ccc}
     1 & 1 & 0  \\
     0 & 1 & 0\\
     0 & 0 & 1
\end{array}\right)>$, $G_3=<\left(\begin{array}{ccc}
     1 & 0 & 1  \\
     0 & 1 & 0\\
     0 & 0 & 1
\end{array}\right)>$, $G_4=<\left(\begin{array}{ccc}
     1 & 1 & 1  \\
     0 & 1 & 0\\
     0 & 0 & 1
\end{array}\right)>$ and $G_5$ the full group. The groups $G_{ii}$ are either trivial or equal to $H_{ii}$ the group generated by the morphism which sends $e_1$ to itself and $e_i$ to $e_1+e_i$. Furthermore for a chosen $G$, the restriction property (Definition \ref{restrictionprop}) requires that some of them are non trivial. Let us summarise the possible values of the $G_{ii}$ for each value of $G$.\\

\begin{tabular}{|c|c|c|c|}
\hline
    $G$ & $G_{22}$ & $G_{33}$ & $G_{44}$ \\
    \hline
    $G_1$ & $\{\id\}$ or $H_{22}$ &  $\{\id\}$ or $H_{33}$ &  $\{\id\}$ or $H_{44}$\\
    \hline
    $G_2$ & $H_{22}$ &  $\{\id\}$ or $H_{33}$ & $H_{44}$\\
     \hline
    $G_3$ &  $\{\id\}$ or $H_{22}$ & $H_{33}$ & $H_{44}$\\
     \hline
    $G_4$ & $H_{22}$ &  $H_{33}$ &  $\{\id\}$ or $H_{44}$\\
     \hline
    $G_5$ & $H_{22}$ & $H_{33}$ & $H_{44}$\\
    \hline
\end{tabular}\\

We find, $15$ possible values for $(G,G_{22},G_{33},G_{44})$. Each one characterising precisely one $\mathcal{G}\in\textit{Groupoid}(V_3)$ such that $H^*(V_3)^\mathcal{G}$ satisfies the required conditions.\\

We give details for only two computations of $H^*(V_3)^\mathcal{G}$. The other ones are left to the reader.\\

\textbf{First case: } If $G$ contains all morphisms whose matrix are in $H$, then by the restriction property $G_{ii}$ contains all the morphisms $\beta$ that satisfy the two conditions, for $i\in\{2,3,4\}$. Hence, all morphisms in $G_{ii}$ are induced by restriction from a morphism in $G$. Therefore, $H^*(V_3)^\mathcal{G}$ is equal to $H^*(V_3)^G$ the algebra of invariant elements under the action of $G$ on $H^*(V_3)$. Let us start by using that an element of $H^*(V_3)^G$ has to be invariant under the action $g_y$ which sends $x$ to $x+y$, and $y$ and $z$ to themselves. Then, every polynomial $P(x,y,z)$ satisfies that $$x^2P(x,y,z)=x(x+y)P(x,y,z)+yxP(x,y,z).$$ Hence, we can express every polynomial $P(x,y,z)$ as $P_1(x(x+y),y,z)+xP_2(x(x+y),y,z)$. Then, since $P_1(x(x+y),y,z)$ is invariant under $g_y$, $xP_2(x(x+y),y,z)$ has also to be invariant, therefore $P_2=0$. So every element in $H^*(V_3)^G$ can be expressed as $P(X,y,z)$ with $X=x(x+y)$.\\ 

But an element of $H^*(V_3)^G$ is also invariant under the action of $g_z$ which send $x$ to $x+z$ and $y$ and $z$ to themselves. $g_z$ sends $X$ to $X+y(y+z)$. We can show as before that $P(X,y,z)$ can be expressed as $P_1(X(X+y(y+z)),y,z)+XP_2(X(X+y(y+z)),y,z)$. Since $P_1(X(X+y(y+z)),y,z)$ is invariant under $g_z$, $P_2=0$. We have $X(X+y(y+z))=x(x+y)(x+z)(x+y+z)$, so every element in $H^*(V_3)^G$ is in $\F\left[y,z,x(x+y)(x+z)(x+y+z)\right]$, but conversely, $y$, $z$ and $x(x+y)(x+z)(x+y+z)$ are invariant under $G$. Hence, $ H^*(V_3)^\mathcal{G}=\F\left[y,z,x(x+y)(x+z)(x+y+z)\right]$. \\

\textbf{Second case: } If $G=G_2$, by the restriction property $G_{22}$ and $G_{44}$ are the groups of isomorphisms of $\F e_1\oplus\F e_i$, with $i=2\text{ or }4$, generated by the morphism whose matrix in the basis $(e_1,e_i)$ is $\left(\begin{array}{ccc}
     1 & 1  \\
     0 & 1 
\end{array}\right).$ The only group which is not fully determined by the restriction property is $G_{33}$ which can either be equal to $\{\id_{\F e_1\oplus\F e_3}\}$, or to the subgroup $B_2$ of $\text{Gl}(\F e_1\oplus\F e_3)$ generated by the morphism $\beta$ whose matrix in the basis $(e_1,e_3)$ is $\left(\begin{array}{ccc}
     1 & 1  \\
     0 & 1 
\end{array}\right).$ 
\begin{enumerate}
    \item If $G_{33}=\{\id_{\F e_1\oplus\F e_3}\}$, then all the sets $\mathcal{G}(U,U')$ are determined by $G_2=\mathcal{G}(V_3,V_3)$, by the restriction property, hence $ H^*(V_3)^\mathcal{G}=H^*(V_3)^{G_2}=\F\left[y,z,x(x+y)\right]$.
    \item If $G_{33}=B_2$, we get our first example which is not an algebra of invariant elements. By Proposition \ref{defdif}, $$H^*(V_3)^\mathcal{G}=\{x\in H^*(V_3)\ ;\ \alpha^*\iota_{U'}^*(x)=\iota_U^*(x)\text{ for all }\alpha\in\mathcal{G}(U,U')\},$$
    which is equal to $$H^*(V_3)^{G_2}\cap\{x\in H^*(V_3)\ ;\ \beta^*\iota_{\F e_1\oplus\F e_3}^*(x)=\iota_{\F e_1\oplus\F e_3}^*(x)\}.$$
    Hence, $H^*(V_3)^\mathcal{G}$ contains the element of $H^*(V_3)^G$ which are in the inverse image of $H^*(\F e_1\oplus\F e_3)^{B_2}=\F\left[x(x+z),z\right]$ under $\iota_{\F e_1\oplus\F e_3}^*$ which maps $x(x+y)$ to $x^2$, $y$ to $0$ and $z$ to $z$. We notice that $\iota_{\F e_1\oplus\F e_3}^*(\F\left[x(x+y),y,z\right])\cap\F\left[x(x+z),z\right]=\F\left[x^2(x+z)^2,z\right]$ and that $\ker(\iota_{\F e_1\oplus\F e_3}^*)=\F\left[x(x+y),y,z\right]y$. Since $\F\left[x(x+y)(x+z)(x+y+z),y,z\right]$ is a sub algebra of $\F\left[x(x+y),y,z\right]$ that surjects onto $\F\left[x^2(x+z)^2,y,z\right]$ under $\iota_{\F e_1\oplus\F e_3}^*$, we get that $$(\iota_{\F e_1\oplus\F e_3}^*)^{-1}(\F\left[x^2(x+z)^2,z\right])=\F\left[y,z,x(x+y)(x+z)(x+y+z)\right]+\F\left[y,z,x(x+y)\right]y.$$ Therefore, $$ H^*(V_3)^\mathcal{G}=\F\left[y,z,x(x+y)(x+z)(x+y+z)\right]+\F\left[y,z,x(x+y)\right]y.$$ 
    
\end{enumerate}

\end{proof}

\begin{remark}
 In Proposition \ref{restrictionprop}, it is worth noticing that, for example, $\F\left[y,z,x(x+y)\right],$ \newline$\F\left[y,z,x(x+z)\right]$ and $\F\left[y,z,x(x+y+z)\right]$ are conjugates. In particular, they are isomorphic. To get a minimal list of isomorphism classes of $nil$-closed, noetherian, integral, unstable algebras of transcendence degree $3$, with a $\F\left[x\right]$-comodule structure whose algebra of primitive elements is isomorphic to $\F\left[y,z\right]$, we should consider the conjugacy classes of the algebras found in Proposition \ref{example627}. (See Remark \ref{gatouzou}.)
\end{remark}

\begin{proposition}
There are exactly $12$ $nil$-closed, noetherian, unstable sub-algebras $K$ of $H^*(V_3)$ of transcendence degree $3$, which are sub $\F\left[x\right]$-comodules of $H^*(V_3)$, for the comodule structure that maps $x$ to $x\otimes1+1\otimes x$, $y$ to $y\otimes 1$ and $z$ to $z\otimes 1$, and such that the algebra of primitive elements of $K$ is $\F\left[z,y(y+z)\right]$. These are :

\begin{enumerate}
    \item $\F\left[z,y(y+z),x(x+y)(x+z)(x+y+z)\right],$
    \item $\F\left[z,y(y+z),x(x+y)(x+z)(x+y+z)\right]\newline+\F\left[z,y(y+z),x(x+y)(x+z)(x+y+z)\right](xz(x+z)+y^2(y+z)),$
    \item $\F\left[z,y(y+z),x(x+z)\right],$
    \item $\F\left[y(y+z),x(x+y)(x+z)(x+y+z)\right]\oplus \F\left[z,y(y+z),x(x+z)\right]z,$
    \item $\F\left[x,z,y(y+z)\right],$
    \item $\F\left[x,z,y(y+z)\right]z\oplus\F\left[y(y+z),x(x+y)(x+z)(x+y+z)\right],$
    \item $\F\left[x,z,y(y+z)\right]y(y+z)\oplus\F\left[z,x(x+z)\right],$
    \item $(\F\left[x,z,y(y+z)\right]y(y+z)\oplus\F\left[z,x(x+z)\right])z\oplus\F\left[y(y+z),x(x+y)(x+z)(x+y+z)\right],$
    \item $\F\left[z,(x+y),y(y+z)\right],$
    \item $\F\left[z,(x+y),y(y+z)\right]z\oplus\F\left[y(y+z),x(x+y)(x+z)(x+y+z)\right],$
    \item $\F\left[(x+y),z,y(y+z)\right]y(y+z)\oplus\F\left[z,(x+y)(x+y+z)\right],$
    \item $(\F\left[z,y(y+z),(x+y)\right]y(y+z)\oplus\F\left[z,(x+y)(x+y+z)\right])z\oplus\F\left[y(y+z),x(x+y)(x+z)(x+y+z)\right].$
    
\end{enumerate}
\end{proposition}

\begin{proof} 
In this case, $\F \left[z,y(y+z)\right]=\pi^*(H^*(V_3/\F e_1)^{B_2})$, for $B_2$ the group generated by the morphism $b_2$ which sends $e'_2$ to itself and $e'_3$ to $e'_2+e'_3$, where $e'_2$ and $e'_3$ are the images of $e_2$ and $e_3$ under the canonical projection. Then, the groupoid $\mathcal{G}'=\mathfrak{g}(B_2)$ is the following:
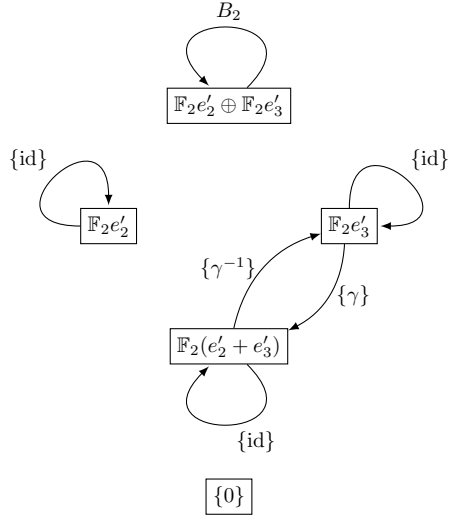
\begin{center}
\begin{tikzpicture}[scale=0.8, every node/.style={transform shape}]
\node[draw] (V) at (5,5) {$\F e'_2\oplus \F e'_3$};
\node[draw] (P2) at (3,3) {$\F e'_2$};
\node[draw] (P3) at (7,3) {$ \F e'_3$};
\node[draw] (P4) at (5,1) {$ \F (e'_2+e'_3)$};
\node[draw] (P1) at (5,-1.5) {$\{0\}$};
\draw[->,>=latex] (V) edge[in=135, out=45, loop] node[above]{$B_2$} (V) ;
\draw[->,>=latex] (P2) edge[in=90, out=180, loop] node[above left]{$\{\id\}$} (P2) ;
\draw[->,>=latex] (P3) edge[in=0, out=90, loop] node[above right]{$\{\id\}$} (P3) ;
\draw[->,>=latex] (P4) edge[in=-135, out=-45, loop] node[below right]{$\{\id\}$} (P4) ;
\draw[->,>=latex] (P3) to[bend left] node[right]{$\{\gamma\}$} (P4) ;
\draw[->,>=latex] (P4) to[bend left] node[left]{$\{\gamma^{-1}\}$} (P3) ;

\end{tikzpicture}
\end{center}

Theorem \ref{principal6} implies that, for $\beta\in\mathcal{G}(N,N')$, when $e_2\in N$, $\pi\circ\beta(e_2)=e'_2$, when $e_3\in N$, $\pi\circ\beta(e_3)=e'_3\text{ or }e'_4$ and when $e_4\in N$, $\pi\circ\beta(e_4)=e'_3\text{ or }e'_4$.  In this case, $G_{34}$ is not necessarily empty. We have to determine $G$, $G_{22}$, $G_{33}$, $G_{44}$ and $G_{34}$ ($G_{43}$ being the set of inverses of morphisms in $G_{34}$). 
By the two conditions on $\beta\in G$, the matrix of $\beta$ is in $H$, for $$H=\{\left(\begin{array}{ccc}
     1 & 0 & 0  \\
     0 & 1 & 0\\
     0 & 0 & 1
\end{array}\right),\left(\begin{array}{ccc}
     1 & 1 & 0  \\
     0 & 1 & 0\\
     0 & 0 & 1
\end{array}\right),\left(\begin{array}{ccc}
     1 & 0 & 1  \\
     0 & 1 & 0\\
     0 & 0 & 1
\end{array}\right),\left(\begin{array}{ccc}
     1 & 1 & 1  \\
     0 & 1 & 0\\
     0 & 0 & 1
\end{array}\right),$$
$$\left(\begin{array}{ccc}
     1 & 0 & 0  \\
     0 & 1 & 1\\
     0 & 0 & 1
\end{array}\right),\left(\begin{array}{ccc}
     1 & 1 & 0  \\
     0 & 1 & 1\\
     0 & 0 & 1
\end{array}\right),\left(\begin{array}{ccc}
     1 & 0 & 1  \\
     0 & 1 & 1\\
     0 & 0 & 1
\end{array}\right),\left(\begin{array}{ccc}
     1 & 1 & 1  \\
     0 & 1 & 1\\
     0 & 0 & 1
\end{array}\right)\}.$$ 
Furthermore, since $b_2\in\mathcal{G}'(\F e'_2\oplus\F e'_3,\F e'_2\oplus\F e'_3)$, by Theorem \ref{principal6} there exist $N$ and $N'$ and $\beta\in\mathcal{G}(N,N')$ such that $\pi$ induces isomorphisms from $N$ and $N'$ to $\F e'_2\oplus\F e'_3$ and $b_2\circ\pi=\pi\circ\beta$. Then, by Theorem \ref{principal2}, $\Bar{\beta}\in G$. Therefore $G$ contains at least one element among $$\{\left(\begin{array}{ccc}
     1 & 0 & 0  \\
     0 & 1 & 1\\
     0 & 0 & 1
\end{array}\right),
\left(\begin{array}{ccc}
     1 & 1 & 0  \\
     0 & 1 & 1\\
     0 & 0 & 1
\end{array}\right),
\left(\begin{array}{ccc}
     1 & 0 & 1  \\
     0 & 1 & 1\\
     0 & 0 & 1
\end{array}\right),\left(\begin{array}{ccc}
     1 & 1 & 1  \\
     0 & 1 & 1\\
     0 & 0 & 1
\end{array}\right)\}.$$

$H$ is isomorphic to the dihedral group $D_4$. It admits ten subgroups, among those the possible values of $G$ are:  $G_1=<\left(\begin{array}{ccc}
     1 & 0 & 0  \\
     0 & 1 & 1\\
     0 & 0 & 1
\end{array}\right)>$, $G_2=<\left(\begin{array}{ccc}
     1 & 0 & 1  \\
     0 & 1 & 1\\
     0 & 0 & 1
\end{array}\right)>$, \newline$G_3=<\left(\begin{array}{ccc}
     1 & 0 & 0  \\
     0 & 1 & 1\\
     0 & 0 & 1
\end{array}\right),\left(\begin{array}{ccc}
     1 & 0 & 1  \\
     0 & 1 & 0\\
     0 & 0 & 1
\end{array}\right)>$,  $G_4=<\left(\begin{array}{ccc}
     1 & 1 & 0  \\
     0 & 1 & 1\\
     0 & 0 & 1
\end{array}\right)>$ and $G_{5}$ the full group.

For $\beta\in G_{ii}$, we necessarily have $\beta(e_1)=e_1$. By Theorem \ref{principal6}, $\beta$ has a block matrix of the following form:
$$\left(\begin{array}{ccc}
     \id_{\F e_1}& \hat{\beta}  \\
     0 &  \id_{N}
\end{array}\right),$$ where $N\in\{\F e_2,\F e_3,  \F (e_2+e_3) \}$ and $\hat{\beta}$ is a morphism from $N$ to $\F e_1$. 

Finally, $G_{34}$ cannot be trivial, otherwise $(b_2)_{\F e'_3}\not\in\mathcal{G}'(\F e'_3,\F e'_4)$, and for $\beta\in G_{34}$, $\beta(e_1)=e_1$ and $\beta(e_3)=e_4\text{ or }e_1+e_4$. We get that the $G_{ii}$ with $i\in\{2,3,4\}$ is a subgroup of $H_{ii}:=\{\id,\beta_{ii}\}$, for $\beta_{ii}$ the morphism which sends $e_1$ to itself and $e_i$ to $e_1+e_i$, and $G_{34}$ is a non trivial subset of the set $H_{34}:=\{\beta_1,\beta_2\}$ where $\beta_1$ and $\beta_2$ are the morphisms from $\F e_1\oplus\F e_3$ to $\F e_1\oplus\F e_4$ whose matrix in the basis $(e_1,e_3)$ and $(e_1,e_4)$ are either $\left(\begin{array}{ccc}
     1 & 1  \\
     0 & 1 
\end{array}\right),$ or $\left(\begin{array}{ccc}
     1 & 0  \\
     0 & 1 
\end{array}\right).$\\

 As in Proposition \ref{example627} (or rather its proof), some values of $G$ imply, by restriction property, the maximality of some of the $G_{ii}$ or the one of $G_{34}$. Furthermore, in this case since the groupoid admits a morphism between $\F e_1\oplus \F e_3$ and $\F e_1\oplus \F e_4$, for $\mathcal{G}$ to be a groupoid we need the compositions of morphisms in $G_{33}$, $G_{34}$, $G_{44}$ and $G_{43}$ to be in $\mathcal{G}$. This is equivalent to requiring that $G_{33}$, $G_{44}$ and $G_{34}$ have the same cardinal.\\

\begin{tabular}{|c|c|c|c|c|}
\hline
    $G$ & $G_{22}$ & $G_{33}$ & $G_{44}$ & $G_{34}$ \\
    \hline
     
     & $\{\id\}$ or & $\{\id\}$ & $\{\id\}$ & $\{\beta_1\}$\\
    \cline{3-5}
    
    $G_1$ & $H_{22}$ & $H_{33}$ & $H_{44}$ & $H_{34}$\\
    \hline
     & $\{\id\}$ or & $\{\id\}$ & $\{\id\}$ & $\{\beta_2\}$\\
    \cline{3-5}
    
    $G_2$ & $H_{22}$ & $H_{33}$ & $H_{44}$ & $H_{34}$\\
    \hline
   
    $G_3$ & $\{\id\}$ or $H_{22}$ & $H_{33}$ & $H_{44}$ & $H_{34}$\\
    \hline

     $G_4$ & $H_{22}$ & $H_{33}$ & $H_{44}$ & $H_{34}$\\
    \hline
     $G_{5}$ & $H_{22}$ & $H_{33}$ & $H_{44}$ & $H_{34}$\\
    \hline
\end{tabular}\\

We find $12$ possible values for $(G,G_{22},G_{33},G_{44},G_{34})$. Each one characterising precisely one $\mathcal{G}\in\textit{Groupoid}(V_3)$ such that $H^*(V_3)^\mathcal{G}$ satisfies the required conditions. We leave the computations of the corresponding $H^*(V_3)^\mathcal{G}$ to the reader.

\end{proof}

\bibliographystyle{alpha}
\bibliography{Biblio.bib}

\begin{center}
    Address : Univ Angers, CNRS, LAREMA, SFR MATHSTIC, F-49000 Angers, France\\
    e-mail : bloede@math.univ-angers.fr
\end{center}

\end{document}